\documentclass[11pt]{article}
\usepackage{epic,latexsym,amssymb,xcolor}
\usepackage{color}
\usepackage{tikz}
\usepackage{amsfonts,epsf,amsmath,leftidx}
\usepackage{float}
\usepackage{pgfplots}
\pgfplotsset{compat=1.15}
\usepackage{mathrsfs}
\usetikzlibrary{arrows}
\usepackage{pdflscape}
\usepackage[bottom=1.5in, top=1.2in, left=1.5in, right=1.5in]{geometry}

\newtheorem{thm}{Theorem}

\newcommand{\qed}{$\Box$}

\let\oldenumerate\enumerate
\renewcommand{\enumerate}{
  \oldenumerate
  \setlength{\itemsep}{0pt}
  \setlength{\parskip}{0pt}
  \setlength{\parsep}{0pt}
}

\begin{document}

\title{Maximal Independent Sets in Polygonal Cacti}

\author{Natawat Klamsakul \, Pantaree Thengarnanchai \, Mattanaporn Suebtangjai,\\
Pailin Kaewperm \, Nuttanon Songsuwan \, and \, Pawaton Kaemawichanurat
\\ \\
Mathematics and Statistics with Applications~(MaSA),\\
Department of Mathematics, Faculty of Science,\\ King Mongkut's University of Technology Thonburi, \\
Bangkok, Thailand \\
\small \tt Email: natawat.kla@kmutt.ac.th, pantaree.mewmew@mail.kmutt.ac.th, \\\small \tt mattanaporn.km@mail.kmutt.ac.th,
pailin.test@mail.kmutt.ac.th, \\
\small \tt nuttanon.19701@gmail.com, pawaton.kae@mail.kmutt.ac.th}

\date{}
\maketitle

\begin{abstract}
Counting the number of maximal independent sets of graphs was started over $50$ years ago by Erd\H{o}s and Mooser. The problem has been continuously studied with a number of variations. Interestingly, when the maximal condition of an independent set is removed, such the concept presents one of topological indices in molecular graphs, the so called Merrifield-Simmons index. In this paper, we applied the concept of bivariate generating function to establish the recurrence relations of the numbers of maximal independent sets of regualr $n$-gonal cacti when $3 \leq n \leq 6$. By the ideas on
meromorphic functions and the growth of power series coefficients, the asymptotic behaviors
through simple functions of these recurrence relations have been established.
\end{abstract}


\section{Introduction}
Let $G = (V(G), E(G))$ be a graph with the vertex set $V(G)$ and the edge set $E(G)$. The \emph{order} of $G$ is $|V(G)|$. For a vertex $v \in V(G)$, the \emph{neighbor set} of $v$ in $G$, $N_{G}(v)$, is the set $\{u | uv \in E(G)\}$. The \emph{degree} $deg_{G}(v)$ of a vertex $v$ in $G$ is $|N_{G}(v)|$. A \emph{cycle} of length $n$ is denoted by $C_{n}$. For a graph $G$ containing a vertex $v$, we say that $v$ is a \emph{cut vertex} of $G$ if $G - v$ has more components than $G$. A maximal subgraph that does not have a cut vertex is called a \emph{block}. A block $B$ is an \emph{end block} if $B$ has exactly one cut vertex of $G$. For a natural number $n \geq 3$, a \emph{regular} $n$\emph{-gonal cactus} is a graph that has exactly two end blocks and all the blocks are cycles of the same length. A regular $n$-gonal cactus $G$ is said to be \emph{ortho} if the two cut vertices of $G$ that belong to a non-end block are adjacent. A regular $n$-gonal cactus $G$ is said to be \emph{meta} if the two cut vertices of $G$ that belong to a non-end block are a pair of vertices at distance two of the polygon. Further, a regular $n$-gonal cactus $G$ is said to be \emph{para} if the two cut vertices of $G$ that belong to a non-end block a pair of vertices at distance three of the polygon.

A vertex subset $I$ of $V(G)$ is \emph{independent} if any pair of vertices in $I$ are not adjacent in $G$. An independent set $I$ is \emph{maximal} if $I \cup \{v\}$ is no longer independent for any vertex $v \in V(G) \setminus I$. A vertex subset $K$ of $V(G)$ is \emph{clique} if any pair of vertices in $K$ are adjacent in $G$. A clique $K$ is \emph{maximal} if $K \cup \{v\}$ is not a clique for any vertex $v \in V(G) \setminus I$. Clearly, if $\alpha(n)$ is the number of maximal independent sets of all graphs of order $n$ and $\omega(n)$ is the number of maximal cliques of all graphs of order $n$, then 
\begin{align*}
    \alpha(n) = \omega(n).
\end{align*}
\indent The study of counting the largest number of maximal independent sets, or cliques, was started over fifty years ago when Erd\H{o}s and Mooser rose the question how large $\alpha(n)$ can be and what is the graphs that satisfy this upper bound. This problem was solved by Moon and Mooser \cite{MM} in 1965. However, the extremal graphs in \cite{MM} are not connected. Hence, it was asked further by Wilf \cite{W}, who found the largest number of maximal independent sets of trees, what is the largest number of maximal independent sets of all connected graphs of order $n$. This problem was solved by Griggs et. al. \cite{GGG}.
\vskip 5 pt

\indent Although the value $\alpha(n)$ for graphs, connected graphs or trees had been found, a number of graph theorists still have paid attention to count the number of maximal independent sets. The study is to give some condition to the graphs or is to focus on some particular classes of graphs. For examples, counting the number of maximal independent sets of graphs with given the number of cycles see \cite{JC,JL,LZZ,SV}, counting the number of maximal independent sets of comparability graphs (the graphs that are constructed from partial order sets) see \cite{DM1}, counting the number of maximal independent sets of unlabelled cycle of length $n$ see \cite{BM}, counting the number of maximal independent sets of random graphs see \cite{DR,DFJ} and for example of studies on maximal independent sets see \cite{AH1,AH,BH}.
\vskip 5 pt

\indent In 1952, Husimi \cite{H} generalized the cluster and irreducible integrals from the book of Statistical Mechanics by Mayer and Mayer \cite{Mayer} to solve problems in the Theory of Condensation. Later the same year, Uhlenbeck \cite{U} pointed out that Husimi's integrals can be interpreted by graphs. The such graphs were called Husimi trees, the graphs whose each edge is in at most one cycle. Husimi trees have been attracted much attention as they can be applied to explain many of condensation phenomena, for example of the studies see \cite{HN,HU,R}. The Husimi trees were known in graph theory literature as cacti after Harary and Palmer \cite{HP} published their classical book on graph enumeration in 1973 which was 20 years since it was first introduced.

\indent For graphs that are applied in chemistry, Merrifield and Simmons \cite{MS1,MS2,MS3,MS4} found relationship between physical properties of hydrocarbons and topological indices of graphs representing their molecular structures such as the number of independent sets, connected sets, irredundant sets and maximal independent sets (kernel). In particular, in \cite{MS1}, Merrifield and Simmons observed that the numbers of independent sets of graphs representing Alkanes have inverse variations to the boiling points and heat of formations of these compounds. By these observations, the number of independent sets of graphs representing molecular structures is known to be Merrifield-Simmons index and has been studied by a number of graph theorists, see \cite{DL,DM,OC} for example. For more studies on topological indices of molecular graphs see \cite{CXL,CGR,DG,D,G,Gu1,Gu2,LYH,Ra,SIA,SLH,TBJK} for example.
\vskip 5 pt

\indent In this paper, we applied the concept of bivariate generating functions, which have been used in \cite{DL,DM,DZ}, to establish recurrence relations of the number of maximal independent sets of ortho, meta and para $n$-gonal cacti when $3 \leq n \leq 6$. Further, the asymptotic behaviors through simple functions of these recurrence relations has been established  by the ideas on meromorphic functions and the growth of power series coefficients which were introduced in \cite{W2}.

%\indent Doslic and Maloy \cite{DM} applied bivariate generating functions to obtain formulae to count the number of independent sets of hexagonal cacti.  Doslic and Maloy proved that

%\begin{thm}
%For a natural number $n \geq 1$, we let $O_{n}, M_{n}$ and $P_{n}$ be the numbers of independent sets of ortho-, meta- and para-hexagonal cactus of $n$ hexagons, respectively. Then
%\begin{align*}
%O_{n} &= \Big(1 + \frac{5}{\sqrt{41}}\Big)(4 + \sqrt{41})^{n} + \Big(1 - \frac{5}{\sqrt{41}}\Big)(4 - \sqrt{41})^{n}\\ \\
%M_{n} &= \frac{2}{5} + \frac{8}{5}\cdot 11^{n}\\ \\
%P_{n} &= \Big(1 + \frac{1}{\sqrt{2}}\Big)(5 + 4\sqrt{2})^{n} + %\Big(1 - \frac{1}{\sqrt{2}}\Big)(5 - 4\sqrt{2})^{n}.
%\end{align*}
%\end{thm}

\section{Main Results}
In this section, we state all our main results as detailed in each of the following subsections. All the results will be summarized in Tables 1 and 2 at the end of this section.
%which are the bivariate generating functions, recurrnce relations and thei asymptotic behavior of $n$-gonal cactus when $3 \leq n \leq 6$.

\subsection{Triangular Cacti}
Let $T(n)$ be a triangular cactus of $n$ triangles as detailed in Figure \ref{Tn0}.

\vskip -0.25 cm
\setlength{\unitlength}{0.5cm}
\begin{figure}[htb]
\begin{center}
\begin{picture}(13, 6.5)

%T1
\put(-2, 2){\circle*{0.2}}
\put(0, 4){\circle*{0.2}}
\put(2, 2){\circle*{0.2}}
\put(-2, 2){\line(1, 1){2}}
\put(-2, 2){\line(1, 0){4}}
\put(0, 4){\line(1, -1){2}}

%\put(-2.2, 1.5){\footnotesize$x^{1}_{1}$}
%\put(1.5, 1.5){\footnotesize$x^{3}_{1} = x^{1}_{2}$}
%\put(-0.2, 4.2){\footnotesize$x^{2}_{1}$}

%T2
\put(2, 2){\circle*{0.2}}
\put(4, 4){\circle*{0.2}}
\put(6, 2){\circle*{0.2}}
\put(2, 2){\line(1, 1){2}}
\put(2, 2){\line(1, 0){4}}
\put(4, 4){\line(1, -1){2}}

%\put(5.5, 1.5){\footnotesize$x^{3}_{2} = x^{1}_{3}$}
%\put(3.8, 4.2){\footnotesize$x^{2}_{2}$}

\put(7.45, 2){\footnotesize$\cdots$}

%Tn
\put(10, 2){\circle*{0.2}}
\put(12, 4){\circle*{0.2}}
\put(14, 2){\circle*{0.2}}
\put(10, 2){\line(1, 1){2}}
\put(10, 2){\line(1, 0){4}}
\put(12, 4){\line(1, -1){2}}

%\put(9.8, 1.5){\footnotesize$x^{1}_{n}$}
%\put(13.8, 1.5){\footnotesize$x^{3}_{n}$}
%\put(11.8, 4.2){\footnotesize$x^{2}_{n}$}

\end{picture}
\end{center}
\vskip -1 cm
\caption{\label{Tn0} The triangular cactus $T(n)$ of $n$ triangles}
\end{figure}
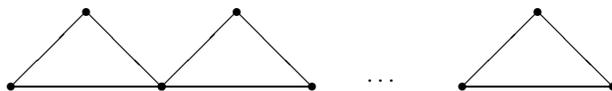
\vskip 5 pt

\begin{thm}\label{bivariate-t}
Let $t(n, k)$ be the number of maximal independent sets containing $k$ vertices of $T(n)$. If $T(x, y)$ is a bi-variate generating function of $t(n, k)$, then 
\begin{align}
T(x, y) = \frac{1 + 2xy + x^{2}y^{2}}{1 - xy - x^{2}y}.\notag
\end{align}
\end{thm}

\begin{thm}\label{recurrence-t}
Let $t(n)$ be the number of maximal independent sets of $T(n)$. Then, $t(1) = 3, t(2) = 5$ and
\begin{align*}
 t(n) = t(n - 1) + t(n - 2) 
\end{align*}
\noindent for $n \geq 3$.
\end{thm}

\noindent From Theorem \ref{recurrence-t}, we obtain the exact formula of $t(n)$ as follows.

\begin{thm}\label{exact-t}
If $t(n)$ is the number of maximal independent sets of $T(n)$. Then,
\begin{align*}
t(n) = 0.10559(-0.61803)^n+1.8944(1.61803)^n.    
\end{align*}
\end{thm}

\noindent Further, by considering the principal part of the expansion of the generating function around the singularity,   we can gain the growth rate of the recurrence relation in Theorem \ref{recurrence-t} that:

\begin{thm}\label{asymp-t}
If $t(n)$ is the number of maximal independent sets of $T(n)$. Then,
\begin{align*}
   t(n)\approx\frac{1.1708}{0.618^{n+1}}. 
\end{align*}

\end{thm}

\noindent For the last formula of $t(n)$ concerning the average size of a maximal independent set of $T(n)$, we find the first order partial derivative of $T(x, y)$ from Theorem \ref{bivariate-t} (the method can be found in \cite{DM,W2} as well). By letting $y = 1$ in $\frac{\partial T(x, y)}{\partial y}$, we will obtain the recurrence relation of the summation of sizes of all maximal independent sets of $T(n)$. By finding its asymptotic formula and dividing by the approximation of $t(n)$ from Theorem \ref{asymp-t}, we have the following theorem.

\begin{thm}\label{average-t}
If $\overline{t}(n)$ is the average size of a maximal independent set of $T(n)$. Then,
\begin{align*}
   \overline{t}(n)\approx 0.72361n+0.33475. 
\end{align*}
\end{thm}

\subsection{Diamond Cacti}

Let $D(n)$ be the diamond (meta-rectangular) cactus of $n$ diamonds as detailed in Figure \ref{Dn0}.

\vskip -0.5 cm
\setlength{\unitlength}{0.5cm}
\begin{figure}[htb]
\begin{center}
\begin{picture}(13, 6)

%D1
\put(-2, 2){\circle*{0.2}}
\put(0, 4){\circle*{0.2}}
\put(2, 2){\circle*{0.2}}
\put(0, 0){\circle*{0.2}}

\put(-2, 2){\line(1, 1){2}}
\put(-2, 2){\line(1, -1){2}}
\put(0, 4){\line(1, -1){2}}
\put(0, 0){\line(1, 1){2}}

%\put(-2.5, 1.5){\footnotesize$x^{1}_{1}$}
%\put(1.3, 1){\footnotesize$x^{3}_{1} = x^{1}_{2}$}
%\put(-0.2, 4.2){\footnotesize$x^{2}_{1}$}
%\put(-0.2, -0.5){\footnotesize$x^{4}_{1}$}

%D2
\put(2, 2){\circle*{0.2}}
\put(4, 4){\circle*{0.2}}
\put(6, 2){\circle*{0.2}}
\put(4, 0){\circle*{0.2}}

\put(2, 2){\line(1, 1){2}}
\put(2, 2){\line(1, -1){2}}
\put(4, 4){\line(1, -1){2}}
\put(4, 0){\line(1, 1){2}}

%\put(5.5, 1){\footnotesize$x^{3}_{2} = x^{1}_{3}$}
%\put(3.8, 4.2){\footnotesize$x^{2}_{2}$}
%\put(3.8, -0.5){\footnotesize$x^{4}_{2}$}

\put(7.45, 2){\footnotesize$\cdots$}

%Dn
\put(10, 2){\circle*{0.2}}
\put(12, 4){\circle*{0.2}}
\put(14, 2){\circle*{0.2}}
\put(12, 0){\circle*{0.2}}

\put(10, 2){\line(1, 1){2}}
\put(10, 2){\line(1, -1){2}}
\put(12, 4){\line(1, -1){2}}
\put(12, 0){\line(1, 1){2}}

%\put(9.6, 1.5){\footnotesize$x^{1}_{n}$}
%\put(13.8, 1.5){\footnotesize$x^{3}_{n}$}
%\put(11.8, 4.2){\footnotesize$x^{2}_{n}$}
%\put(12, -0.5){\footnotesize$x^{4}_{n}$}

\end{picture}
\end{center}
\vskip 0.5 cm
\caption{\label{Dn0} The diamond cactus $D(n)$ of $n$ diamonds}
\end{figure}
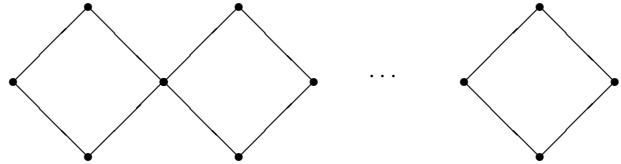
\vskip 5 pt
\vskip 5 pt

\begin{thm}\label{bivariate-d}
Let $d(n, k)$ be the number of maximal independent sets of $D(n)$ containing  $k$ vertices. If $D(x, y)$ is a bi-variate generating function of $d(n, k)$, then 
\begin{align}
D(x, y) = \frac{1 + xy^{2} -xy - x^{2}y^{4} - x^{3}y^{3} + 2x^{2}y^{3}}{1 - xy - xy^{2} + x^{2}y^{3} - x^{3}y^{3}}.\notag
\end{align}
\end{thm}

\begin{thm}\label{recurrence-d}
Let $d(n)$ be the number of maximal independent sets of $D(n)$. Then, $d(1) = 2, d(2) = 4, d(3) = 7$ and
\begin{align*}
    d(n) = 2d(n - 1) - d(n - 2) + d(n - 3)
\end{align*}
\noindent for $n \geq 4$.
\end{thm}

\noindent Next, we give a brief argument to find the exact formula of $d(n)$ from Theorem \ref{recurrence-d}. By replacing $d(n - i)$ with $r^{n - i}$ for all $i \in \{0, 1, 2, 3\}$ and divide by $r^{n - 3}$ throughout the equation, we have
\begin{align*}
    r^{3} - 2r^{2} + r - 1 = 0
\end{align*}
\noindent which we obtain complex roots $0.12256+0.74486i$ and $0.12256-0.74486i$. Then, the polar form is considered together with Euler's formula to cancel the imaginary parts of the solutions. We obtain the exact formula of $d(n)$ as follows.

\begin{thm}\label{exact-d}
If $d(n)$ is the number of maximal independent sets of $D(n)$. Then,
\begin{align*}
d(n) = 1.2672(1.7549)^n - 0.27965\cos(1.4077n - 0.765).  
\end{align*}
\end{thm}

\noindent By similar arguments as Theorem \ref{asymp-t}, we can approximate the recurrence relation in Theorem \ref{recurrence-d} that:

\begin{thm}\label{asymp-d}
If $d(n)$ is the number of maximal independent sets of $D(n)$. Then,
\begin{align*}
   d(n)\approx \frac{0.6213}{0.5698^{n+1}}.
\end{align*}
\end{thm}

\noindent By similar arguments as Theorem \ref{average-t}, we obtain the following theorem.

\begin{thm}\label{average-d}
If $\overline{d}(n)$ is the average size of a maximal independent set of $D(n)$. Then,
\begin{align*}
   \overline{d}(n)\approx 1.2345n+0.89458. 
\end{align*}
\end{thm}

\subsection{Square Cacti}

Let $S(n)$ be the square (ortho-rectangular) cactus of $n$ squares as detailed in Figure \ref{Sn0}. 
\vskip 5 pt

\vskip -0.25 cm
\setlength{\unitlength}{0.5cm}
\begin{figure}[htb]
\begin{center}
\begin{picture}(13, 6.5)

%S1
\put(-2, 3){\circle*{0.2}}
\put(-2, 6){\circle*{0.2}}
\put(1, 3){\circle*{0.2}}
\put(1, 6){\circle*{0.2}}

\put(-2, 3){\line(0, 1){3}}
\put(-2, 6){\line(1, 0){3}}
\put(1, 3){\line(-1, 0){3}}
\put(1, 6){\line(0, -1){3}}

%\put(-2.5, 6.2){\footnotesize$x^{1}_{1}$}
%\put(1, 6.2){\footnotesize$x^{2}_{1}$}
%\put(1.2, 3.3){\footnotesize$x^{3}_{1} = x^{1}_{2}$}
%\put(-2.5, 2.5){\footnotesize$x^{4}_{1}$}

%S2
\put(1, 3){\circle*{0.2}}
\put(4, 3){\circle*{0.2}}
\put(4, 0){\circle*{0.2}}
\put(1, 0){\circle*{0.2}}

\put(1, 3){\line(1, 0){3}}
\put(4, 3){\line(0, -1){3}}
\put(4, 0){\line(-1, 0){3}}
\put(1, 0){\line(0, 1){3}}

%\put(4, 3.2){\footnotesize$x^{2}_{2}$}
%\put(4.3, 0){\footnotesize$x^{3}_{2} = x^{1}_{3}$}
%\put(0.6, -0.6){\footnotesize$x^{4}_{2}$}

\put(7, 3){\footnotesize$\cdots$}

%Sn
\put(10, 3){\circle*{0.2}}
\put(10, 6){\circle*{0.2}}
\put(13, 3){\circle*{0.2}}
\put(13, 6){\circle*{0.2}}

\put(10, 3){\line(0, 1){3}}
\put(10, 6){\line(1, 0){3}}
\put(13, 3){\line(-1, 0){3}}
\put(13, 6){\line(0, -1){3}}

%\put(9.5, 6.2){\footnotesize$x^{1}_{n}$}
%\put(13, 6.2){\footnotesize$x^{2}_{n}$}
%\put(13, 2.6){\footnotesize$x^{3}_{n}$}%
%\put(9.5, 2.6){\footnotesize$x^{4}_{n}$}

\end{picture}
\end{center}
\vskip 0.5 cm
\caption{The square cactus $S(n)$ of $n$ squares}
\label{Sn0} 
\end{figure}
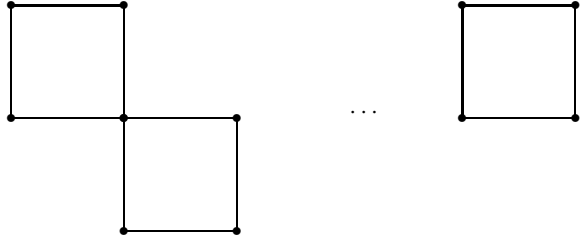
\vskip 5 pt

\begin{thm}\label{bivariate-s}
Let $s(n, k)$ be the number of maximal independent sets of $S(n)$ containing $k$ vertices. If $S(x, y)$ is a bi-variate generating function of $s(n, k)$, then 
\begin{align}
S(x, y) =  \frac{1 - 2xy + 2xy^{2} - 2x^{2}y^{3} + x^{2}y^{2} + x^{2}y^{4}}{1 - 2xy + x^{2}y^{2} - x^{2}y^{3}}\notag
\end{align}
\end{thm}

\begin{thm}\label{recurrence-s}
Let $s(n)$ be the number of maximal independent sets of $S(n)$. Then $s(1) = 2$ and, for $n \geq 2$, we have that
\begin{align*}
    s(n) = 2s(n - 1). 
\end{align*}
Namely, $s(n) = 2^{n}$.
\end{thm}

\noindent By finding the summation of sizes of all maximal independent sets with similar arguments as Theorem \ref{average-t} and divided by $2^{n}$, we obtain the following theorem.

\begin{thm}\label{average-s}
If $\overline{s}(n)$ is the average size of a maximal independent set of $S(n)$. Then,
\begin{align*}
   \overline{s}(n)\approx 1.25n+0.75. 
\end{align*}
\end{thm}

\subsection{Pentagonal Cacti}
\begin{figure}[H]
\centering
\definecolor{ududff}{rgb}{0.30196078431372547,0.30196078431372547,1}
\resizebox{0.8\textwidth}{!}{%

\begin{tikzpicture}[line cap=round,line join=round,>=triangle 45,x=1cm,y=1cm]
\clip(0.07084031875793423,-1.7138790246542182) rectangle (9.988466422567084,3.9142216784276616);
\draw [line width=0.3pt] (1,1)-- (3,1);
\draw [line width=0.3pt] (1,1)-- (0.8111130238479834,1.9874845007960271);
\draw [line width=0.3pt] (0.8111130238479834,1.9874845007960271)-- (1.510822840987893,2.382972658309889);
\draw [line width=0.3pt] (1.510822840987893,2.382972658309889)-- (2.1902512141527324,1.9874845007960271);
\draw [line width=0.3pt] (2.1902512141527324,1.9874845007960271)-- (1.8150445006139406,0);
\draw [line width=0.3pt] (3,1)-- (3.204323412906225,0);
\draw [line width=0.3pt] (1.8150445006139406,0)-- (2.49447287377878,-0.41586661024974864);
\draw [line width=0.3pt] (2.49447287377878,-0.41586661024974864)-- (3.204323412906225,0);
\draw [line width=0.3pt] (5,1)-- (6,1);
\draw [line width=0.3pt] (5,1)-- (4.8065574869367405,1.997625222783562);
\draw [line width=0.3pt] (4.8065574869367405,1.997625222783562)-- (5.4859858601015805,2.4032541022849587);
\draw [line width=0.3pt] (5.4859858601015805,2.4032541022849587)-- (6.195836399229025,2.007765944771097);
\draw [line width=0.3pt] (6.195836399229025,2.007765944771097)-- (6,1);
\begin{scriptsize}
\draw [fill=black] (1,1) circle (1pt);
\draw [fill=black] (2,1) circle (1pt);
\draw [fill=black] (3,1) circle (1pt);
\draw [fill=black] (0.8111130238479834,1.9874845007960271) circle (1pt);
\draw [fill=black] (2.1902512141527324,1.9874845007960271) circle (1pt);
\draw [fill=black] (1.510822840987893,2.382972658309889) circle (1pt);
\draw [fill=black] (1.8150445006139406,0) circle (1pt);
\draw [fill=black] (3.204323412906225,0) circle (1pt);
\draw [fill=black] (2.49447287377878,-0.41586661024974864) circle (1pt);
\draw [fill=black] (5,1) circle (1pt);
\draw [fill=black] (4.8065574869367405,1.997625222783562) circle (1pt);
\draw [fill=black] (5.4859858601015805,2.4032541022849587) circle (1pt);
\draw [fill=black] (6.195836399229025,2.007765944771097) circle (1pt);
\draw [fill=black] (6,1) circle (1pt);

\draw[color=black] (4,1) node {$...$};

\end{scriptsize}
\end{tikzpicture}

 }%
\vskip -1 cm
\caption{The pentagonal cactus $P(n)$ of $n$ pentagons.}
\label{orthopent}
\end{figure}
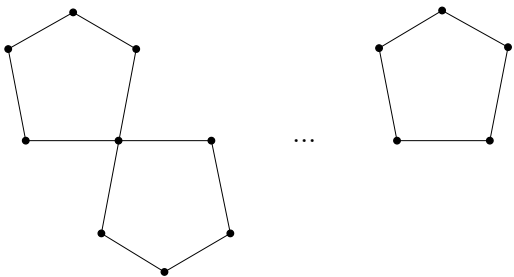

Let $P(n)$ be the (ortho) pentagonal cactus of $n$ pentagons as illustrated in Figure \ref{orthopent}. 

\vskip -1 cm
\vskip 5 pt

\begin{thm}\label{bivariate-p}
Let $p(n, k)$ be the number of maximal independent sets of $P(n)$ containing $k$ vertices. If $P(x, y)$ is a bi-variate generating function of $p(n, k)$, then 
\begin{align}
P(x, y) = \frac{1 + 4x^3y^5 - 4x^3y^4 + 4x^2y^4 - x^2y^3 + 4xy^2}{1 - 5x^{2}y^{3} - xy^{2} - 4x^{3}y^{4}}.\notag
\end{align}
\end{thm}

\begin{thm}\label{recurrence-p}
Let $p(n)$ be the number of maximal independent sets of $P(n)$. Then $p(1) = 5, p(2) = 13, p(3) = 42$ and for $n \geq 4$, we have that
\begin{align*}
    p(n) = p(n - 1) + 5p(n - 2) + 4p(n - 3). 
\end{align*}
\end{thm}

\noindent From Theorem \ref{recurrence-p}, we obtain the exact formula of $p(n)$ by similar arguments of Theorem \ref{exact-d} as follows.

\begin{thm}\label{exact-p}
If $p(n)$ is the number of maximal independent sets of $P(n)$. Then,
\begin{align*}
p(n) = 1.4492(3.0606)^n - 0.56449\cos(3.5897n - 0.428).   
\end{align*}
\end{thm}

\noindent The recurrence relation in
Theorem \ref{recurrence-p} can be approximated as:

\begin{thm}\label{asymp-p}
If $p(n)$ is the number of maximal independent sets of $P(n)$. Then,
\begin{align*}
   p(n)\approx \frac{0.4735}{0.3267^{n+1}}. 
\end{align*}
\end{thm}

\noindent The average size of a maximal independent set is detailed as follows.

\begin{thm}\label{average-p}
If $\overline{p}(n)$ is the average size of a maximal independent set of $P(n)$. Then,
\begin{align*}
   \overline{p}(n)\approx 1.5516n+0.55731. 
\end{align*}
\end{thm}

\subsection{Meta-Pentagonal Cacti}
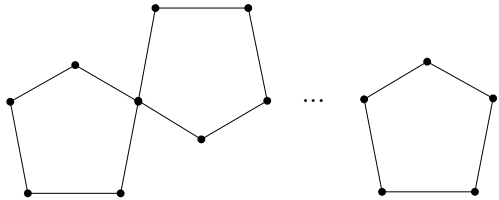
\begin{figure}[H]
\centering
\definecolor{ududff}{rgb}{0.30196078431372547,0.30196078431372547,1}
\resizebox{0.8\textwidth}{!}{%

\begin{tikzpicture}[line cap=round,line join=round,>=triangle 45,x=1cm,y=1cm]
\clip(-0.8215432161451385,-2.798936277320456) rectangle (9.096082887664004,2.829164425761429);
\draw [line width=0.3pt] (0.6247932864612085,-0.997722231544381)-- (0.4359063103091919,-0.010237730748353324);
\draw [line width=0.3pt] (0.4359063103091919,-0.010237730748353324)-- (1.1356161274491015,0.38525042676550875);
\draw [line width=0.3pt] (1.1356161274491015,0.38525042676550875)-- (1.8150445006139408,-0.010237730748353324);
\draw [line width=0.3pt] (3,1)-- (3.204323412906225,0);
\draw [line width=0.3pt] (1.8150445006139406,0)-- (2.49447287377878,-0.41586661024974864);
\draw [line width=0.3pt] (2.49447287377878,-0.41586661024974864)-- (3.204323412906225,0);
\draw [line width=0.3pt] (4.441491495385481,-0.9814692730960854)-- (5.441491495385481,-0.9814692730960854);
\draw [line width=0.3pt] (4.441491495385481,-0.9814692730960854)-- (4.248048982322222,0.016155949687477897);
\draw [line width=0.3pt] (4.248048982322222,0.016155949687477897)-- (4.927477355487062,0.42178482918887417);
\draw [line width=0.3pt] (4.927477355487062,0.42178482918887417)-- (5.637327894614506,0.026296671675012803);
\draw [line width=0.3pt] (5.637327894614506,0.026296671675012803)-- (5.441491495385481,-0.9814692730960854);
\draw [line width=0.3pt] (0.6247932864612085,-0.997722231544381)-- (1.6247932864612085,-0.997722231544381);
\draw [line width=0.3pt] (1.8150445006139408,-0.010237730748353324)-- (1.6247932864612085,-0.997722231544381);
\draw [line width=0.3pt] (1.8150445006139406,0)-- (2,1);
\draw [line width=0.3pt] (2,1)-- (3,1);
\begin{scriptsize}
\draw [fill=black] (0.6247932864612085,-0.997722231544381) circle (1pt);
\draw [fill=black] (1.6247932864612085,-0.997722231544381) circle (1pt);
\draw [fill=black] (3,1) circle (1pt);
\draw [fill=black] (0.4359063103091919,-0.010237730748353324) circle (1pt);
\draw [fill=black] (1.8150445006139408,-0.010237730748353324) circle (1pt);
\draw [fill=black] (1.1356161274491015,0.38525042676550875) circle (1pt);
\draw [fill=black] (1.8150445006139406,0) circle (1pt);
\draw [fill=black] (3.204323412906225,0) circle (1pt);
\draw [fill=black] (2.49447287377878,-0.41586661024974864) circle (1pt);
\draw [fill=black] (4.441491495385481,-0.9814692730960854) circle (1pt);
\draw [fill=black] (4.248048982322222,0.016155949687477897) circle (1pt);
\draw [fill=black] (4.927477355487062,0.42178482918887417) circle (1pt);
\draw [fill=black] (5.637327894614506,0.026296671675012803) circle (1pt);
\draw [fill=black] (5.441491495385481,-0.9814692730960854) circle (1pt);
\draw [fill=black] (2,1) circle (1pt);

\draw[color=black] (3.7,0) node {$...$};

\end{scriptsize}
\end{tikzpicture}

 }%
\vskip -1 cm
\caption{The meta-pentagonal cactus $M(n)$ of $n$ pentagons.}
\label{metapent}
\end{figure}

We let $M(n)$ be the meta-pentagonal cactus of $n$ pentagons as shown in Figure \ref{metapent}.

\vskip -1 cm
\vskip 5 pt

\begin{thm}\label{bivariate-m}
Let $m(n, k)$ be the number of maximal independent sets of $M(n)$ containing $k$ vertices. If $M(x, y)$ is a bi-variate generating function of $m(n, k)$, then 
\begin{align}
M(x, y) = \frac{-11x^{2}y^{4} + 5xy^{2} + 2x^{3}y^{6} + 2x^{3}y^{5} - x^{2}y^{3}}{1 - 4xy^{2} - xy + 3x^{2}y^{3} + 4x^{2}y^{4} - 2x^{3}y^{5}} + 1.\notag
\end{align}
\end{thm}

\begin{thm}\label{recurrence-m}
Let $m(n)$ be the number of maximal independent sets of $M(n)$. Then, $m(1) = 5, m(2) = 13$  for $n \geq 3$, we have that
\begin{align*}
     m(n) = 3m(n - 1) - m(n - 2).  
\end{align*}
\end{thm}

\noindent From Theorem \ref{recurrence-m}, we obtain the exact formula of $m(n)$ as follows.

\begin{thm}\label{exact-m}
If $m(n)$ is the number of maximal independent sets of $M(n)$. Then,
\begin{align*}
m(n) = 0.10573(0.38197)^n+1.89443(2.61803)^n.   
\end{align*}
\end{thm}

\noindent  The recurrence relation in Theorem \ref{recurrence-m} can be approximated as:

\begin{thm}\label{asymp-m}
If $m(n)$ is the number of maximal independent sets of $M(n)$. Then,
\begin{align*}
   m(n)\approx \frac{0.7236}{0.382^{n+1}}. 
\end{align*}
\end{thm}

\noindent The average size of a maximal independent set is detailed as follows.

\begin{thm}\label{average-m}
If $\overline{m}(n)$ is the average size of a maximal independent set of $M(n)$. Then,
\begin{align*}
   \overline{m}(n)\approx 1.7236n+0.24033. 
\end{align*}
\end{thm}

\subsection{Meta-Hexagonal Cacti}
\begin{figure}[H]
\centering
\definecolor{ududff}{rgb}{0.30196078431372547,0.30196078431372547,1}
\resizebox{0.7\textwidth}{!}{%

\begin{tikzpicture}[line cap=round,line join=round,>=triangle 45,x=1cm,y=1cm]
\clip(-0.42006259682951436,-2.0634421903550786) rectangle (6.717157112872825,2.344411290892576);
\draw [line width=0.3pt] (0,1)-- (0,0);
\draw [line width=0.3pt] (0,1)-- (0.5090377499796499,1.3861738649290605);
\draw [line width=0.3pt] (0.5090377499796499,1.3861738649290605)-- (1,1);
\draw [line width=0.3pt] (1,1)-- (1,0);
\draw [line width=0.3pt] (0,0)-- (0.49945650687446763,-0.39225270822641545);
\draw [line width=0.3pt] (0.49945650687446763,-0.39225270822641545)-- (1,0);
\draw [line width=0.3pt] (1,0)-- (1.49925108788594,0.39590652351985406);
\draw [line width=0.3pt] (1.49925108788594,0.39590652351985406)-- (2,0);
\draw [line width=0.3pt] (1,0)-- (1,-1);
\draw [line width=0.3pt] (1,-1)-- (1.5065488585502573,-1.4066428305665222);
\draw [line width=0.3pt] (1.5065488585502573,-1.4066428305665222)-- (2,-1);
\draw [line width=0.3pt] (2,0)-- (2,-1);
\draw [line width=0.3pt] (3.5,1)-- (3.5,0);
\draw [line width=0.3pt] (3.5,1)-- (4.002386425746781,1.3957011045313252);
\draw [line width=0.3pt] (4.002386425746781,1.3957011045313252)-- (4.5,1);
\draw [line width=0.3pt] (4.5,1)-- (4.5,0);
\draw [line width=0.3pt] (4.5,0)-- (4.002386425746781,-0.4141460202193678);
\draw [line width=0.3pt] (3.5,0)-- (4.002386425746781,-0.4141460202193678);
\begin{scriptsize}
\draw [fill=black] (0,0) circle (1pt);
\draw [fill=black] (1,0) circle (1pt);
\draw [fill=black] (1,1) circle (1pt);
\draw [fill=black] (0,1) circle (1pt);
\draw [fill=black] (0.5090377499796499,1.3861738649290605) circle (1pt);
\draw [fill=black] (0.49945650687446763,-0.39225270822641545) circle (1pt);
\draw [fill=black] (1.49925108788594,0.39590652351985406) circle (1pt);
\draw [fill=black] (2,0) circle (1pt);
\draw [fill=black] (1,-1) circle (1pt);
\draw [fill=black] (2,-1) circle (1pt);
\draw [fill=black] (1.5065488585502573,-1.4066428305665222) circle (1pt);
\draw [fill=black] (3.5,0) circle (1pt);
\draw [fill=black] (4.5,0) circle (1pt);
\draw [fill=black] (3.5,1) circle (1pt);
\draw [fill=black] (4.5,1) circle (1pt);
\draw [fill=black] (4.002386425746781,1.3957011045313252) circle (1pt);
\draw [fill=black] (4.002386425746781,-0.4141460202193678) circle (1pt);

\draw[color=black] (2.7,0) node {$...$};

\end{scriptsize}
\end{tikzpicture}

 }%
\vskip -1 cm
\caption{The meta-hexagonal cactus $H(n)$ of $n$ hexagons.}
\label{mhcac}
\end{figure}
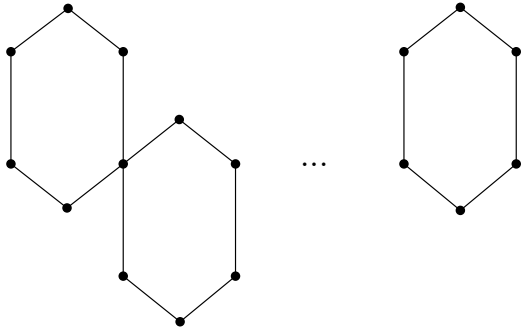

We may let the meta-hexagonal cactus $H(n)$ of $n$ hexagons be the graph in Figure \ref{mhcac}.

\vskip -1 cm
\vskip 5 pt

\begin{thm}\label{bivariate-Hh}
Let $h(n, k)$ be the number of maximal independent sets of $H(n)$ containing $k$ vertices. If $H(x, y)$ is a bi-variate generating function of $h(n, k)$, then 
\begin{align}
H(x, y) &= \frac{1-2x^{2}y^{3}  + xy^{2} - x^{3}y^{5} + x^{2}y^{5} + 5x^{2}y^{4} + xy^{3} - x^{2}y^{6} + x^{3}y^{7}}{1 - 3x^{2}y^{3} - xy^{3} - 2xy^{2} - x^{3}y^{5} + x^{2}y^{5} + x^{2}y^{4} - x^{3}y^{6}}.\notag
\end{align}
\end{thm}

\begin{thm}\label{recurrence-h}
Let $h(n)$ be the number of maximal independent sets of $H(n)$. Then $h(1) = 5, h(2) = 19, h(3) = 64$ and, for $n \geq 4$, we have that
%h(4) = 221, h(5) = 765
\begin{align*}
     h(n) = 3h(n - 1) + h(n - 2) + 2h(n - 3).  
\end{align*}
\end{thm}

\noindent From Theorem \ref{recurrence-h}, we obtain the exact formula of $h(n)$ by similar arguments as Theorem \ref{exact-d} as follows.

\begin{thm}\label{exact-h}
If $h(n)$ is the number of maximal independent sets of $H(n)$. Then,
\begin{align*}
h(n) = 1.5499(3.4567)^n - 0.65668\cos(1.8757n - 0.88031).   
\end{align*}
\end{thm}

\noindent The asymptotic analysis provides that:

\begin{thm}\label{asymp-h}
If $h(n)$ is the number of maximal independent sets of $H(n)$. Then,
\begin{align*}
   h(n)\approx \frac{0.4484}{0.2893^{n+1}}. 
\end{align*}
\end{thm}

\noindent The average size of a maximal independent set is detailed as follows.

\begin{thm}\label{average-h}
If $\overline{h}(n)$ is the average size of a maximal independent set of $H(n)$. Then,
\begin{align*}
   \overline{h}(n)\approx 1.9409n+0.56686. 
\end{align*}
\end{thm}

\subsection{Para-Hexagonal Cacti}
\begin{figure}[H]
\centering
\definecolor{ududff}{rgb}{0.30196078431372547,0.30196078431372547,1}
\resizebox{0.9\textwidth}{!}{%

\begin{tikzpicture}[line cap=round,line join=round,>=triangle 45,x=1cm,y=1cm]
\clip(-0.8083969725917267,-1.5837349250118342) rectangle (8.231945270112451,3.999462124674588);
\draw [line width=0.3pt] (-0.02405899207354767,0.9828713769549792)-- (0.3745474533070522,1.5702914017263894);
\draw [line width=0.3pt] (0.3745474533070522,1.5702914017263894)-- (0.9829467646774412,1.5807810450258788);
\draw [line width=0.3pt] (0.9829467646774412,1.5807810450258788)-- (1.3916045138741162,0.9767710067152169);
\draw [line width=0.3pt] (1.3916045138741162,0.9767710067152169)-- (1.8011389420376196,0.395451352183569);
\draw [line width=0.3pt] (1.3916045138741162,0.9767710067152169)-- (1.7906492987381306,1.5702914017263894);
\draw [line width=0.3pt] (1.3916045138741162,0.9767710067152169)-- (0.9934364079769307,0.37447206558459006);
\draw [line width=0.3pt] (-0.02405899207354767,0.9828713769549792)-- (0.38503709660654173,0.37447206558459006);
\draw [line width=0.3pt] (0.38503709660654173,0.37447206558459006)-- (0.9934364079769307,0.37447206558459006);
\draw [line width=0.3pt] (1.7906492987381306,1.5702914017263894)-- (2.388558966809031,1.5702914017263894);
\draw [line width=0.3pt] (2.388558966809031,1.5702914017263894)-- (2.7976550554891206,0.9723817336554896);
\draw [line width=0.3pt] (2.7976550554891206,0.9723817336554896)-- (2.3990486101085207,0.37447206558459006);
\draw [line width=0.3pt] (1.8011389420376196,0.395451352183569)-- (2.3990486101085207,0.37447206558459006);
\draw [line width=0.3pt] (4.3920808370115205,1.5807810450258788)-- (3.9829847483314316,0.9618920903560004);
\draw [line width=0.3pt] (4.3920808370115205,1.5807810450258788)-- (4.97950086178293,1.5702914017263894);
\draw [line width=0.3pt] (4.97950086178293,1.5702914017263894)-- (5.391604513874117,0.9767710067152169);
\draw [line width=0.3pt] (5.391604513874117,0.9767710067152169)-- (5.00048014838191,0.38496170888407955);
\draw [line width=0.3pt] (3.9829847483314316,0.9618920903560004)-- (4.381591193712031,0.38496170888407955);
\draw [line width=0.3pt] (4.381591193712031,0.38496170888407955)-- (5.00048014838191,0.38496170888407955);
\begin{scriptsize}
\draw [fill=black] (-0.02405899207354767,0.9828713769549792) circle (1pt);
\draw [fill=black] (1.3916045138741162,0.9767710067152169) circle (1pt);
\draw [fill=black] (2.7976550554891206,0.9723817336554896) circle (1pt);
\draw [fill=black] (3.9829847483314316,0.9618920903560004) circle (1pt);
\draw [fill=black] (5.391604513874117,0.9767710067152169) circle (1pt);
\draw [fill=black] (0.3745474533070522,1.5702914017263894) circle (1pt);
\draw [fill=black] (0.9829467646774412,1.5807810450258788) circle (1pt);
\draw [fill=black] (0.38503709660654173,0.37447206558459006) circle (1pt);
\draw [fill=black] (0.9934364079769307,0.37447206558459006) circle (1pt);
\draw [fill=black] (1.7906492987381306,1.5702914017263894) circle (1pt);
\draw [fill=black] (2.388558966809031,1.5702914017263894) circle (1pt);
\draw [fill=black] (1.8011389420376196,0.395451352183569) circle (1pt);
\draw [fill=black] (2.3990486101085207,0.37447206558459006) circle (1pt);
\draw [fill=black] (4.3920808370115205,1.5807810450258788) circle (1pt);
\draw [fill=black] (4.97950086178293,1.5702914017263894) circle (1pt);
\draw [fill=black] (4.381591193712031,0.38496170888407955) circle (1pt);
\draw [fill=black] (5.00048014838191,0.38496170888407955) circle (1pt);

\draw[color=black] (3.3,1) node {$...$};

\end{scriptsize}
\end{tikzpicture}

 }%
\vskip -3 cm
\caption{The para-hexagonal cactus $G(n)$ of $n$ hexagons.}
\label{phn}
\end{figure}
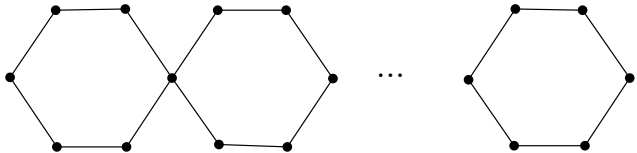

Let $G(n)$ be the para-hexagonal cactus of $n$ hexagons as shown in Figure \ref{phn}.

\vskip -3 cm
\vskip 5 pt

\begin{thm}\label{bivariate-ph}
Let $g(n, k)$ be the number of maximal independent sets of $G(n)$ containing $k$ vertices. If $G(x, y)$ is a bi-variate generating function of $g(n, k)$, then 
\begin{align}
G(x, y) = \frac{-2x^{2}y^{5} - 3x^{2}y^{4} - x^{3}y^{5} + x^{2}y^{6} +x^{3}y^{6} -xy^{2} +4xy^{3} -xy + 1}{1 - xy - xy^{2} +4x^{2}y^{3} -3xy^{2} + x^{2}y^{4} - x^{3}y^{5} - x^{2}y^{5}}.\notag
\end{align}
\end{thm}

\begin{thm}\label{recurrence-ph}
Let $g(n)$ be the number of maximal independent sets of $G(n)$. Then $g(1) = 5, g(2) = 19, g(3) = 76$ and, for $n \geq 4$, we have that
\begin{align*}
     g(n) = 5g(n - 1) - 4g(n - 2) + g(n - 3).  
\end{align*}
\end{thm}

\noindent From Theorem \ref{recurrence-ph}, we obtain the exact formula of $g(n)$ by similar arguments as Theorem \ref{exact-d} as follows.

\begin{thm}\label{exact-g}
If $g(n)$ is the number of maximal independent sets of $G(n)$. Then,
\begin{align*}
g(n) = 1.115(4.0796)^n - 1.3688\cos(0.3777n + 1.5292).   
\end{align*}
\end{thm}

\noindent We can approximate the recurrence relation in Theorem \ref{recurrence-ph} that:

\begin{thm}\label{asymp-ph}
If $g(n)$ is the number of maximal independent sets of $G(n)$. Then,
\begin{align*}
   g(n)\approx \frac{0.2733}{0.2451^{n+1}}. 
\end{align*}
\end{thm}

\noindent The average size of a maximal independent set is detailed as follows.

\begin{thm}\label{average-g}
If $\overline{g}(n)$ is the average size of a maximal independent set of $G(n)$. Then,
\begin{align*}
   \overline{g}(n)\approx 2.9143n+0.88506. 
\end{align*}
\end{thm}

\subsection{Ortho-Hexagonal Cacti}
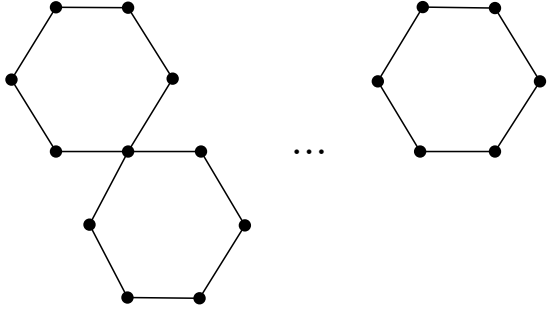
\begin{figure}[H]
\centering
\definecolor{ududff}{rgb}{0.30196078431372547,0.30196078431372547,1}
\resizebox{0.8\textwidth}{!}{%

\begin{tikzpicture}[line cap=round,line join=round,>=triangle 45,x=1cm,y=1cm]
\clip(-0.32537133026435466,-1.5238107651131987) rectangle (5.916913749858403,2.33134280175197);
\draw [line width=0.3pt] (0.7636911477878787,0.983396738469664)-- (1.068196050713675,0.4982589804623496);
\draw [line width=0.3pt] (1.068196050713675,0.4982589804623496)-- (0.7636911477878787,0);
\draw [line width=0.3pt] (2.469888418299665,0.4797553368451276)-- (2.7770094292928365,0.9875825272815013);
\draw [line width=0.3pt] (2.7770094292928365,0.9875825272815013)-- (3.270375092325276,0.9809481903139401);
\draw [line width=0.3pt] (3.270375092325276,0.9809481903139401)-- (3.5778226918823863,0.4797553368451276);
\draw [line width=0.3pt] (3.5778226918823863,0.4797553368451276)-- (3.270375092325276,0);
\draw [line width=0.3pt] (2.469888418299665,0.4797553368451276)-- (2.7589040782508967,0);
\draw [line width=0.3pt] (2.7589040782508967,0)-- (3.270375092325276,0);
\draw [line width=0.3pt] (-0.03310388590148339,0.4916246434947884)-- (0.2709525523276073,0.9861831599960398);
\draw [line width=0.3pt] (-0.03310388590148339,0.4916246434947884)-- (0.2709525523276073,0);
\draw [line width=0.3pt] (0.2709525523276073,0.9861831599960398)-- (0.7636911477878787,0.983396738469664);
\draw [line width=0.3pt] (0.2709525523276073,0)-- (0.7636911477878787,0);
\draw [line width=0.3pt] (0.5,-0.5)-- (0.7636911477878787,0);
\draw [line width=0.3pt] (1.2619236772230777,0)-- (1.561584020328892,-0.5050167892307563);
\draw [line width=0.3pt] (0.5,-0.5)-- (0.7594444433566199,-0.9983600370269141);
\draw [line width=0.3pt] (1.2518740925457486,-1.0033848293655787)-- (1.561584020328892,-0.5050167892307563);
\draw [line width=0.3pt] (0.7636911477878787,0)-- (1.2619236772230777,0);
\draw [line width=0.3pt] (0.7594444433566199,-0.9983600370269141)-- (1.2518740925457486,-1.0033848293655787);
\begin{scriptsize}
\draw [fill=black] (-0.03310388590148339,0.4916246434947884) circle (1pt);
\draw [fill=black] (0.7636911477878787,0.983396738469664) circle (1pt);
\draw [fill=black] (0.7636911477878787,0) circle (1pt);
\draw [fill=black] (1.068196050713675,0.4982589804623496) circle (1pt);
\draw [fill=black] (2.469888418299665,0.4797553368451276) circle (1pt);
\draw [fill=black] (2.7770094292928365,0.9875825272815013) circle (1pt);
\draw [fill=black] (2.7589040782508967,0) circle (1pt);
\draw [fill=black] (3.270375092325276,0) circle (1pt);
\draw [fill=black] (3.5778226918823863,0.4797553368451276) circle (1pt);
\draw [fill=black] (3.270375092325276,0.9809481903139401) circle (1pt);
\draw [fill=black] (0.2709525523276073,0.9861831599960398) circle (1pt);
\draw [fill=black] (0.2709525523276073,0) circle (1pt);
\draw [fill=black] (1.2619236772230777,0) circle (1pt);
\draw [fill=black] (0.7594444433566199,-0.9983600370269141) circle (1pt);
\draw [fill=black] (1.2518740925457486,-1.0033848293655787) circle (1pt);
\draw [fill=black] (1.561584020328892,-0.5050167892307563) circle (1pt);
\draw [fill=black] (0.5,-0.5) circle (1pt);

\draw[color=black] (2,0) node {$...$};

\end{scriptsize}
\end{tikzpicture}

 }%
%\vskip -1 cm
\caption{The ortho-hexagonal cactus $Q(n)$ of $n$ hexagons.}
\label{ohn}
\end{figure}

The ortho-hexagonal cactus $Q(n)$ of $n$ hexagons is shown in Figure \ref{ohn}.

\vskip -2 cm
\vskip 5 pt

\begin{thm}\label{bivariate-qh}
Let $q(n, k)$ be the number of maximal independent sets of $Q(n)$ containing $k$ vertices. If $Q(x, y)$ is a bi-variate generating function of $q(n, k)$, then 
\begin{align*}
Q(x, y) = \frac{x^{2}y^{6}+2xy^{3}+1}{1 - x^{2}y^{5} - x^{2}y^{4} - x^{2}y^{3} - 3xy^2}.
\end{align*}
\end{thm}

\begin{thm}\label{recurrence-qh}
Let $q(n)$ be the number of maximal independent sets of $Q(n)$. Then $q(1) = 5, q(2) = 19$ and, for $n \geq 3$, we have that
%q(3) = 72
\begin{align*}
     q(n) = 3q(n - 1) + 3q(n - 2). 
\end{align*}
\end{thm}

\noindent From Theorem \ref{recurrence-qh}, we obtain the exact formula of $q(n)$ as follows.

\begin{thm}\label{exact-q}
If $q(n)$ is the number of maximal independent sets of $Q(n)$. Then,
\begin{align*}
q(n) = 0.01201(-0.79129)^n+1.32132(3.79129)^n.   
\end{align*}
\end{thm}

\noindent In a similar fashion, we can approximate the recurrence relation in Theorem \ref{recurrence-qh} that:

\begin{thm}\label{asymp-qh}
If $q(n)$ is the number of maximal independent sets of $Q(n)$. Then,
\begin{align*}
   q(n)\approx \frac{0.3485}{0.2638^{n+1}}. 
\end{align*}
\end{thm}

\noindent The average size of a maximal independent set is detailed as follows.

\begin{thm}\label{average-q}
If $\overline{q}(n)$ is the average size of a maximal independent set of $Q(n)$. Then,
\begin{align*}
   \overline{q}(n)\approx 2n+0.41742. 
\end{align*}
\end{thm}

\indent Finally, we summarize all our results in this paper in Tables 1 and 2 as follows. In Table 1, the entries of the colum G.F. are obtained by letting $y = 1$ in each bivariate generating function.

\small{

\[ \begin{array}{c|c|c|l}
\text{Cacti}  & \text {Bivariate Generating Functions}   & \text{G.F.}        & ~~~~~~~~~~~~~~~~\text{Recurrence Relations} \\
\hline
                  &                                              &                                                               &\\
T(n)              & \frac{1 + 2xy + x^{2}y^{2}}{1 - xy - x^{2}y} & \frac{1 + 2x + x^{2}}{1 - x - x^{2}}           &   t(n)=t(n - 1)+t(n - 2)     \\
                  &                        &               &  t(1)=3,t(2)=5,t(3)=8             \\
%\hline
                  &                        &                                   &\\
D(n)              & \frac{1 + xy^{2} -xy - x^{2}y^{4} - x^{3}y^{3} + 2x^{2}y^{3}}{1 - xy - xy^{2} + x^{2}y^{3} - x^{3}y^{3}} & \frac{1 + x^{2} - x^{3}}{1 - 2x + x^{2} - x^{3}}&   d(n)=2d(n - 1)-d(n - 2)+d(n - 3)\\
                  &                        &                                   & d(1)=2,d(2)=4,d(3)=7,d(4)=12      \\
                  &                        &                                   &\\
%\hline
                  &                        &                                   &\\
S(n)              & \frac{1 - 2xy + 2xy^{2} - 2x^{2}y^{3} + x^{2}y^{2} + x^{2}y^{4}}{1 - 2xy + x^{2}y^{2} - x^{2}y^{3}} & \frac{1}{1 - 2x}       &   s(n) = 2s(n-1)           \\
                  &                        &                                   &   s(1)=2                          \\
%\hline
                  &                        &                                   &\\
P(n)              &\frac{1 + 4x^3y^5 - 4x^3y^4 + 4x^2y^4 - x^2y^3 + 4xy^2}{1 - 5x^{2}y^{3} - xy^{2} - 4x^{3}y^{4}} & \frac{3x^{2} + 4x + 1}{ - 4x^{3} - 5x^{2} - x + 1} &   p(n) = p(n - 1) + 5p(n - 2) + 4p(n - 3) \\
                  &                        &                                   &p(1)=5,p(2)=13,p(3)=42,p(4)=127  \\
                  &                        &                                   &\\
%\hline
                  &                        &                                   &\\
M(n)              & \frac{-11x^{2}y^{4} + 5xy^{2} + 2x^{3}y^{6} + 2x^{3}y^{5} - x^{2}y^{3}}{1 - 4xy^{2} - xy + 3x^{2}y^{3} + 4x^{2}y^{4} - 2x^{3}y^{5}} + 1 & \frac{-x^{2} + 2x + 1}{x^{2} - 3x + 1} &  m(n) = 3m(n - 1) - m(n - 2)  \\
                  &                        &                                   & m(1)=5,m(2)=13,m(3)=34,m(4)=89            \\
%\hline
                  &                        &             &\\
H(n)              & \frac{1-2x^{2}y^{3}  + xy^{2} - x^{3}y^{5} + x^{2}y^{5} + 5x^{2}y^{4} + xy^{3} - x^{2}y^{6} + x^{3}y^{7}}{1 - 3x^{2}y^{3} - xy^{3} - 2xy^{2} - x^{3}y^{5} + x^{2}y^{5} + x^{2}y^{4} - x^{3}y^{6}} & \frac{3x^{2} + 2x + 1}{-2x^{3} - x^{2} - 3x + 1}      &   h(n) = 3h(n - 1) + h(n - 2) + 2h(n - 3)    \\
                  &                        &                                  & h(1)=5,h(2)=19,h(3)=64,h(4)=221  \\
%\hline
                  &                        &             &\\
G(n)              & \frac{-2x^{2}y^{5} - 3x^{2}y^{4} - x^{3}y^{5} + x^{2}y^{6} +x^{3}y^{6} -xy^{2} +4xy^{3} -xy + 1}{1 - xy - xy^{2} +4x^{2}y^{3} -3xy^{2} + x^{2}y^{4} - x^{3}y^{5} - x^{2}y^{5}} & \frac{1 -2x^{2}}{1 - 5x + 4x^{2} - x^{3}}    &   g(n) = 5g(n - 1) - 4g(n - 2) + g(n - 3)\\
                  &                        &                                  & g(1)=5,g(2)=19,g(3)=76,g(4)=309  \\
%\hline
                  &                        &             &\\
Q(n)              & \frac{x^{2}y^{6}+2xy^{3}+1}{1 - x^{2}y^{5} - x^{2}y^{4} - x^{2}y^{3} - 3xy^2} & \frac{x^{2} + 2x + 1}{1 - 3x - 3x^{2}}    &   q(n) = 3q(n - 1) + 3q(n - 2)    \\
                  &                        &                                  & q(1)=5,q(2)=19,q(3)=72,q(4)=273 \\
%\hline
\end{array}\]

\vskip 10 pt
\begin{center}
\footnotesize{Table 1: Bivariate generating functions, generating functions and recurrence relations of $n$-gonal cacti.}
\end{center}

}

\small{

\[ \begin{array}{c|c|c|c}
\text{Cacti}  &  \text{Exact Formulae}                      & \text{Asymptotic Formulae}                & \text{Average Size}\\
\hline
                  &                                        & &                                                              \\
T(n)              & 0.10559(-0.61803)^n+1.8944(1.61803)^n           &  \frac{1.17082}{(0.61803)^{n+1}}  & 0.72361n+0.33475\\
                  &&&\\
%\hline
                  &&&\\
D(n)              & 1.2672(1.7549)^n - 0.27965\cos(1.4077n - 0.765) & \frac{0.62126}{(0.56984)^{n+1}}    & 1.2345n+0.89458\\
                  &&&\\
%\hline
                  &&&\\
S(n)              & 2^n                                             & 2^n                                & 1.25n+0.75\\
                  &&&\\
%\hline
                  &&&\\
P(n)              & 1.4492(3.0606)^n - 0.56449\cos(3.5897n - 0.428) &\frac{0.47351}{(0.32673)^{n+1}}     & 1.5516n+0.55731\\
                  &&& \\
%\hline
                  &&&\\
M(n)              & 0.10573(0.38197)^n+1.89443(2.61803)^n           &\frac{0.72361}{(0.38196)^{n+1}}     & 1.7236n+0.24033\\
                  &&&\\
%\hline
                  &&&\\
H(n)              &1.5499(3.4567)^n - 0.65668\cos(1.8757n - 0.88031) & \frac{0.44836}{(0.28930)^{n+1}}   & 1.9409n+0.56686\\
                  &&&\\
%\hline
                  &&&\\
G(n)              & 1.115(4.0796)^n - 1.3688\cos(0.3777n + 1.5292)   & \frac{0.27331}{(0.24512)^{n+1}}   & 2.9143n+0.88506\\
                  &&&\\
%\hline
                  &&&\\
Q(n)              & 0.01201(-0.79129)^n+1.32132(3.79129)^n           & \frac{0.34851}{(0.26376)^{n+1}}   & 2n+0.41742\\
                  &&&\\
%\hline
\end{array}\]

\vskip 10 pt
\begin{center}
\footnotesize{Table 2: Asymptotic formulae, exact formulae and average size of a maximal independent set of the cacti.}
\end{center}
%\label{tab:1}

%\end{table}

}

\section{Proofs}
\subsection{Triangular Cacti}
First, we name all the vertices of $T(n)$ as detailed in Figure \ref{Tn}.

\vskip -0.25 cm
\setlength{\unitlength}{0.8cm}
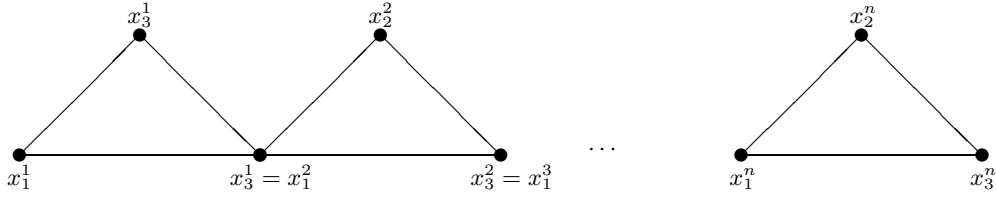
\begin{figure}[htb]
\begin{center}
\begin{picture}(13, 6.5)

%T1
\put(-2, 2){\circle*{0.2}}
\put(0, 4){\circle*{0.2}}
\put(2, 2){\circle*{0.2}}
\put(-2, 2){\line(1, 1){2}}
\put(-2, 2){\line(1, 0){4}}
\put(0, 4){\line(1, -1){2}}

\put(-2.2, 1.5){\footnotesize$x^{1}_{1}$}
\put(1.5, 1.5){\footnotesize$x^{1}_{3} = x^{2}_{1}$}
\put(-0.2, 4.2){\footnotesize$x^{1}_{3}$}

%T2
\put(2, 2){\circle*{0.2}}
\put(4, 4){\circle*{0.2}}
\put(6, 2){\circle*{0.2}}
\put(2, 2){\line(1, 1){2}}
\put(2, 2){\line(1, 0){4}}
\put(4, 4){\line(1, -1){2}}

\put(5.5, 1.5){\footnotesize$x^{2}_{3} = x^{3}_{1}$}
\put(3.8, 4.2){\footnotesize$x^{2}_{2}$}

\put(7.45, 2){\footnotesize$\cdots$}

%Tn
\put(10, 2){\circle*{0.2}}
\put(12, 4){\circle*{0.2}}
\put(14, 2){\circle*{0.2}}
\put(10, 2){\line(1, 1){2}}
\put(10, 2){\line(1, 0){4}}
\put(12, 4){\line(1, -1){2}}

\put(9.8, 1.5){\footnotesize$x^{n}_{1}$}
\put(13.8, 1.5){\footnotesize$x^{n}_{3}$}
\put(11.8, 4.2){\footnotesize$x^{n}_{2}$}

\end{picture}
\end{center}
\vskip -1 cm
\caption{\label{Tn} A labelled triangular cactus of $n$ triangles}
\end{figure}

\noindent Then, we recall that
\vskip 5 pt

\indent $t(n)$ = the number of all maximal independent sets of $T(n)$

\noindent and

\indent $t(n, k)$ = the number of maximal independent sets containing $k$ vertices of $T(n)$.
\vskip 5 pt

\noindent Thus,
\begin{align}\label{tn}
t(n) = \sum_{k \geq 0}t(n, k).    
\end{align}
\vskip 5 pt

\noindent Further, we let 
\begin{align}
T(x) = \sum_{n \geq 0}t(n)x^{n}\notag
\end{align}
\noindent be the generating function of $t(n)$ and we let
\begin{align}
T(x, y) = \sum_{n \geq 0}\sum_{k \geq 0}t(n, k)x^{n}y^{k}\notag
\end{align}
\noindent be the bi-variate generating function of $t(n, k)$. It is worth noting that, when $y = 1$, we have
\begin{align}\label{txytotn}
T(x, 1) = \sum_{n \geq 0}(\sum_{k \geq 0}t(n, k)(1)^{k})x^{n} = \sum_{n \geq 0}t(n)x^{n} = T(x).
\end{align}

\noindent Next, we let $\overline{T}(n)$ be constructed from $T(n)$ by joining one vertex to one vertex of degree two of the $n^{th}$ triangle. 
The graph $\overline{T}(n)$ is illustrated by Figure \ref{figtbarn}.

\vskip -0.25 cm
\setlength{\unitlength}{0.7cm}
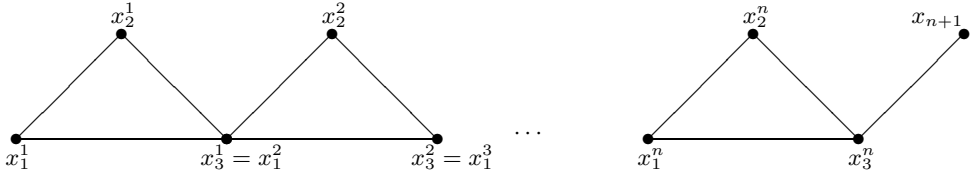
\begin{figure}[htb]
\begin{center}
\begin{picture}(13, 6.5)

%T1
\put(-2, 2){\circle*{0.2}}
\put(0, 4){\circle*{0.2}}
\put(2, 2){\circle*{0.2}}
\put(-2, 2){\line(1, 1){2}}
\put(-2, 2){\line(1, 0){4}}
\put(0, 4){\line(1, -1){2}}

\put(-2.2, 1.5){\footnotesize$x^{1}_{1}$}
\put(1.5, 1.5){\footnotesize$x^{1}_{3} = x^{2}_{1}$}
\put(-0.2, 4.2){\footnotesize$x^{1}_{2}$}

%T2
\put(2, 2){\circle*{0.2}}
\put(4, 4){\circle*{0.2}}
\put(6, 2){\circle*{0.2}}
\put(2, 2){\line(1, 1){2}}
\put(2, 2){\line(1, 0){4}}
\put(4, 4){\line(1, -1){2}}

\put(5.5, 1.5){\footnotesize$x^{2}_{3} = x^{3}_{1}$}
\put(3.8, 4.2){\footnotesize$x^{2}_{2}$}

\put(7.45, 2){\footnotesize$\cdots$}

%Tn
\put(10, 2){\circle*{0.2}}
\put(12, 4){\circle*{0.2}}
\put(14, 2){\circle*{0.2}}
\put(10, 2){\line(1, 1){2}}
\put(10, 2){\line(1, 0){4}}
\put(12, 4){\line(1, -1){2}}

\put(9.8, 1.5){\footnotesize$x^{n}_{1}$}
\put(13.8, 1.5){\footnotesize$x^{n}_{3}$}
\put(11.8, 4.2){\footnotesize$x^{n}_{2}$}

\put(16, 4){\circle*{0.2}}
\put(16, 4){\line(-1, -1){2}}

\put(15, 4.2){\footnotesize$x_{n + 1}$}

\end{picture}
\end{center}
\vskip -1 cm
\caption{\label{figtbarn} The graph $\overline{T}(n)$}
\end{figure}
\vskip 5 pt

\noindent Then, we let 
\vskip 5 pt

\indent $\overline{t}(n, k)$ = the number of maximal independent sets containing $k$ vertices of $\overline{T}(n)$.
\vskip 5 pt

\noindent and 
\begin{align}
\overline{T}(x, y) = \sum_{n \geq 0}\sum_{k \geq 0}\overline{t}(n, k)x^{n}y^{k}\notag
\end{align}
\noindent be the bi-variate generating function of $\overline{t}(n, k)$.
\vskip 5 pt

\indent We are ready to prove Theorem \ref{bivariate-t}.

\noindent \emph{Proof of Theorem \ref{bivariate-t}.} First, we will establish a recurrence relation of $t(n, k)$. Let $x^{n}_{1}$ be the vertex of degree $4$ of the $n^{th}$ triangle of $T_{n}$ and let $D$ be a maximal independent set of $T(n)$ containing $k$ vertices. We distinguish two cases. 
\vskip 5 pt

\noindent \textbf{Case 1:}  $x^{n}_{1} \in D$.\\
\indent Thus, the two vertices of degree two $x^{n}_{2}, x^{n}_{3}$ are not in $D$. Further, $x^{n - 1}_{2}, x^{n - 1}_{1} \notin D$. Removing $x^{n}_{1}, x^{n}_{2}, x^{n}_{3}, x^{n - 1}_{2}, x^{n - 1}_{1}$ from $T(n)$ results in $\overline{T}(n - 3)$. Thus $D = D' \cup \{x^{n}_{1}\}$ where $D'$ is a maximal independent set of $\overline{T}(n - 3)$ containing $k - 1$ vertices. Clearly, there are $\overline{t}(n - 3, k - 1)$ possibilities of $D'$ yielding that there are $\overline{t}(n - 3, k - 1)$ possibilities of $D$.
\vskip 5 pt

\noindent \textbf{Case 2:} $x^{n}_{1} \notin D$.\\
\indent By maximality of $D$, either $x^{n}_{2} \in D$ or $x^{n}_{3} \in D$. In each case, $|D \cap \{x^{n}_{1}, x^{n}_{2}, x^{n}_{3}\}| = 1$. Removing $x^{n}_{1}, x^{n}_{2}, x^{n}_{3}$ from $D$ results in $\overline{T}(n - 2)$. Thus, $D = D' \cup \{w\}$ where $w \in \{x^{n}_{2}, x^{n}_{3}\}$ and $D'$ is a maximal independent set containing $k - 1$ vertices of $\overline{T}(n - 2, k - 1)$. Therefore, for each $w \in \{x^{n}_{2}, x^{n}_{3}\}$, there are $\overline{t}(n - 2, k - 1)$ possibilities of $D'$ yielding that there are totally $2\overline{t}(n - 2, k - 1)$ possibilities of $D$.
\vskip 5 pt

\indent From Cases 1 and 2, we have that
\begin{align}\label{t1}
t(n, k) = \overline{t}(n - 3, k - 1) + 2\overline{t}(n - 2, k - 1)
\end{align}
\noindent For $n \geq 3$ and $k \geq 1$, we multiply $x^{n}y^{k}$ throughout (\ref{t1}) and sum over all $x^{n}y^{k}$. Thus, we have that
\begin{align}\label{t2}
\sum_{n \geq 3}\sum_{k \geq 1}t(n, k)x^{n}y^{k} = \sum_{n \geq 3}\sum_{k \geq 1}\overline{t}(n - 3, k - 1)x^{n}y^{k} + 2\sum_{n \geq 3}\sum_{k \geq 1}\overline{t}(n - 2, k - 1)x^{n}y^{k}.
\end{align}

\noindent We first consider the term $\displaystyle{\sum_{n \geq 3}\sum_{k \geq 1}t(n, k)x^{n}y^{k}}$. It can be checked that 

\begin{center}
    $t(0, 0) = 1, t(0, 1) = 0$ and $ t(0, k) = 0$ for all $k \geq 2$
    \vskip 5 pt
    $t(1, 0) = 0, t(1, 1) = 3$ and $ t(1, k) = 0$ for all $k \geq 2$
    \vskip 5 pt
    $t(2, 0) = 0, t(2, 1) = 1, t(2, 2) = 4$ and $ t(2, k) = 0$ for all $k \geq 3$.
    \vskip 5 pt
\end{center}

Thus,
\begin{align}\label{t3}
\sum_{n \geq 3}\sum_{k \geq 1}t(n, k)x^{n}y^{k} 
                                      &= (\sum_{n \geq 3}\sum_{k \geq 1}t(n, k)x^{n}y^{k}  + t(0, 0) + t(1, 1)xy + t(2, 1)x^{2}y + t(2, 2)x^{2}y^{2})\notag\\
                                      &~~~~ - t(0, 0) - t(1, 1)xy - t(2, 1)x^{2}y - t(2, 2)x^{2}y^{2}\notag\\
                                      &= T(x, y)- t(0, 0) - t(1, 1)xy - t(2, 1)x^{2}y - t(2, 2)x^{2}y^{2}\notag\\
                                      &= T(x, y) - 1 - 3xy - x^{2}y - 4x^{2}y^{2}.
\end{align}

\noindent Now, we consider the term $\displaystyle{\sum_{n \geq 3}\sum_{k \geq 1}\overline{t}(n - 3, k - 1)x^{n}y^{k}}$. Clearly, 
\begin{align}\label{t4}
\sum_{n \geq 3}\sum_{k \geq 1}\overline{t}(n - 3, k - 1)x^{n}y^{k} &= x^{3}y\sum_{n \geq 3}\sum_{k \geq 1}\overline{t}(n - 3, k - 1)x^{n - 3}y^{k - 1}\notag\\
                                      &= x^{3}y\sum_{n \geq 0}\sum_{k \geq 0}\overline{t}(n, k)x^{n}y^{k}\notag\\
                                      &= x^{3}y\overline{T}(x, y).
\end{align}

\noindent Finally, we consider the term $\displaystyle{2\sum_{n \geq 3}\sum_{k \geq 1}\overline{t}(n - 2, k - 1)x^{n}y^{k}}$. It can be check that

\begin{center}
    $\overline{t}(0, 0) = 0, \overline{t}(0, 1) = 2$ and $\overline{t}(0, k) = 0$ for all $k \geq 2$.
\end{center} 

Thus, 
\begin{align}\label{t5}
2\sum_{n \geq 3}\sum_{k \geq 1}\overline{t}(n - 2, k - 1)x^{n}y^{k}
    &= 2x^{2}y\sum_{x \geq 3}\sum_{k \geq 1}\overline{t}(n - 2, k - 1)x^{n - 2}y^{k - 1}\notag\\
    &= 2x^{2}y\sum_{x \geq 1}\sum_{k \geq 0}\overline{t}(n, k)x^{n}y^{k}\notag\\
    &= 2x^{2}y(\sum_{x \geq 1}\sum_{k \geq 0}\overline{t}(n, k)x^{n}y^{k} + \overline{t}(0, 1)y - \overline{t}(0, 1)y)\notag\\
    &= 2x^{2}y(\sum_{x \geq 0}\sum_{k \geq 0}\overline{t}(n, k)x^{n}y^{k} - 2y)\notag\\
    &= 2x^{2}y\overline{T}(x, y) - 4x^{2}y^{2}.
\end{align}

\noindent Plugging (\ref{t3}), (\ref{t4}) and (\ref{t5}) to (\ref{t2}), we have

\begin{align}
T(x, y) - 1 - 3xy - x^{2}y - 4x^{2}y^{2} &= x^{3}y\overline{T}(x, y) + 2x^{2}y\overline{T}(x, y) - 4x^{2}y^{2}\notag
\end{align}

\noindent which can be solved that

\begin{align}\label{semi final t6}
T(x, y) = 1 + 3xy + x^{2}y + (x^{3}y + 2x^{2}y)\overline{T}(x, y).
\end{align}

\indent Next, we will establish the recurrence relation of $\overline{t}(n, k)$. For the graph $\overline{T}(n)$, we let $x_{n + 1}$ be the vertex of degree one who is adjacent to $x^{n}_{3}$ of $T_{n}$. Further, we let $D$ be a maximal independent set of $\overline{T}(n)$ containing $k$ vertices. There are $2$ cases.
\vskip 5 pt

\noindent \textbf{Case 1:} $x_{n + 1} \in D$\\
\indent Thus $x^{n}_{3} \notin D$. Removing $x^{n}_{3}, x_{n + 1}$ from $\overline{T}(n)$ results in $\overline{T}(n - 1)$. Thus, $D = D' \cup \{x_{n + 1}\}$ where $D'$ is a maximal independent set of $\overline{T}(n - 1)$ containing $k - 1$ vertices. Therefore, there are $\overline{t}(n - 1, k - 1)$ possibilities of $D'$ yielding that there are $\overline{t}(n - 1, k - 1)$ possibilities of $D$.
\vskip 5 pt

\noindent \textbf{Case 2:} $x_{n + 1} \notin D$\\
\indent By maximality of $D$, $x^{n}_{3} \in D$. Hence, $x^{n}_{1}, x^{n}_{2} \notin D$. Removing $x^{n}_{1}, x^{n}_{2}, x^{n}_{3}, x_{n + 1}$ from $\overline{T}(n)$ results in $\overline{T}(n - 2)$. Thus, $D = D' \cup \{x^{n}_{3}\}$ where $D'$ is a maximal independent set of $\overline{T}(n - 2)$ containing $k - 1$ vertices. Therefore, there are $\overline{t}(n - 2, k - 1)$ possibilities of $D'$ yielding that there are $\overline{t}(n - 2, k - 1)$ possibilities of $D$.
\vskip 5 pt

\indent From Cases 1 and 2, we have that
\begin{align}\label{tbar1}
    \overline{t}(n, k) = \overline{t}(n - 1, k - 1) + \overline{t}(n - 2, k - 1).
\end{align}
\noindent For $n \geq 2$ and $k \geq 1$, we multiply $x^{n}y^{k}$ throughout (\ref{tbar1})  and sum over all $x^{n}y^{k}$. Thus, 
\begin{align}\label{tbar2}
\sum_{n \geq 2}\sum_{k \geq 1}\overline{t}(n, k)x^{n}y^{k} = \sum_{n \geq 2}\sum_{k \geq 1}\overline{t}(n - 1, k - 1)x^{n}y^{k} + \sum_{n \geq 2}\sum_{k \geq 1}\overline{t}(n - 2, k - 1)x^{n}y^{k}.
\end{align}

\noindent We first consider the term $\displaystyle{\sum_{n \geq 2}\sum_{k \geq 1}\overline{t}(n, k)x^{n}y^{k}}$. It can be checked that 

 \begin{center}
    $\overline{t}(0, 0) = 0, \overline{t}(0, 1) = 2$ and $\overline{t}(0, k) = 0$ for all $k \geq 2$
    \vskip 5 pt
    $\overline{t}(1, 0) = 0, \overline{t}(1, 1) = 1, \overline{t}(1, 2) = 2$ and $\overline{t}(0, k) = 0$ for all $k \geq 3$.
    \vskip 5 pt
\end{center}

\noindent Thus, 
\begin{align}\label{tbar3}
\sum_{n \geq 2}\sum_{k \geq 1}\overline{t}(n, k)x^{n}y^{k}
    &= (\sum_{n \geq 2}\sum_{k \geq 1}\overline{t}(n, k)x^{n}y^{k} + \overline{t}(0, 1)y + \overline{t}(1, 1)xy + \overline{t}(1, 2)xy^{2})\notag\\
    &~~~~ - \overline{t}(0, 1)y - \overline{t}(1, 1)xy - \overline{t}(1, 2)xy^{2}\notag\\
    &= \sum_{n \geq 0}\sum_{k \geq 0}\overline{t}(n, k)x^{n}y^{k} - 2y - xy - 2xy^{2}\notag\\
    &= \overline{T}(x, y)  - 2y - xy - 2xy^{2}
\end{align}
\vskip 5 pt

\noindent We next consider the term $\displaystyle{\sum_{n \geq 2}\sum_{k \geq 1}\overline{t}(n - 1, k - 1)x^{n}y^{k}}$. Recall that 

\begin{center}
    $\overline{t}(0, 0) = 0, \overline{t}(0, 1) = 2$ and $\overline{t}(0, k) = 0$ for all $k \geq 2$.
    \vskip 5 pt
\end{center}

\noindent Thus, 
\begin{align}\label{tbar4}
\sum_{n \geq 2}\sum_{k \geq 1}\overline{t}(n - 1, k - 1)x^{n}y^{k}
    &= xy\sum_{n \geq 2}\sum_{k \geq 1}\overline{t}(n - 1, k - 1)x^{n - 1}y^{k - 1}\notag\\
    &= xy\sum_{n \geq 1}\sum_{k \geq 0}\overline{t}(n, k)x^{n}y^{k}\notag\\
    &= xy(\sum_{n \geq 1}\sum_{k \geq 0}\overline{t}(n, k)x^{n}y^{k} + \overline{t}(0, 1)y - \overline{t}(0, 1)y)\notag\\
    &= xy(\sum_{n \geq 0}\sum_{k \geq 0}\overline{t}(n, k)x^{n}y^{k} - 2y)\notag\\
    &= xy\overline{T}(x, y) - 2xy^{2}
\end{align} 
\vskip 5 pt

\noindent Finally, we consider the term $\displaystyle{\sum_{n \geq 2}\sum_{k \geq 1}\overline{t}(n - 2, k - 1)x^{n}y^{k}}$. Clearly, 

\begin{align}\label{tbar5}
\sum_{n \geq 2}\sum_{k \geq 1}\overline{t}(n - 2, k - 1)x^{n}y^{k}
    &= x^{2}y\sum_{n \geq 2}\sum_{k \geq 1}\overline{t}(n - 2, k - 1)x^{n - 2}y^{k - 1}\notag\\
    &= x^{2}y\sum_{n \geq 0}\sum_{k \geq 0}\overline{t}(n, k)x^{n}y^{k}\notag\\
    &= x^{2}y\overline{T}(x, y).
\end{align}

\noindent Plugging (\ref{tbar3}), (\ref{tbar4}) and (\ref{tbar5}) to (\ref{tbar2}), we have 

\begin{align}
\overline{T}(x, y)  - 2y - xy - 2xy^{2} = xy\overline{T}(x, y) - 2xy^{2} + x^{2}y\overline{T}(x, y)\notag
\end{align}

\noindent which can be solved that

\begin{align}\label{semifinaltbar}
\overline{T}(x, y)   = \frac{2y + xy}{1 - xy - x^{2}y}.
\end{align}

\noindent By plugging (\ref{semifinaltbar}) to (\ref{semi final t6}), we have

\begin{align}
T(x, y) = \frac{1 + 2xy + x^{2}y^{2}}{1 - xy - x^{2}y}
\end{align}
\noindent as required.
\qed
\vskip 5 pt

\indent Next, we will prove Theorem \ref{recurrence-t}.
\vskip 5 pt

\noindent \emph{Proof of Theorem \ref{recurrence-t}} From the equation of Theorem \ref{bivariate-t}, by letting $y = 1$ and (\ref{txytotn}), we have that

\begin{align}
\sum_{n \geq 0}t(n)x^{n} &= T(x, 1)\notag\\
                         &= \frac{1 + 2x + x^{2}}{1 - x - x^{2}}
\end{align}

\noindent which can  be solved that

\begin{align}\label{recursivet-1}
x^{2} + 2x + 1 &= (1 - x - x^{2})\sum_{n \geq 0}t(n)x^{n}\notag\\
               &= \sum_{n \geq 0}t(n)x^{n} - \sum_{n \geq 0}t(n)x^{n + 1} - \sum_{n \geq 0}t(n)x^{n + 2}\notag\\
               &= \sum_{n \geq 0}t(n)x^{n} - \sum_{n \geq 1}t(n - 1)x^{n} - \sum_{n \geq 2}t(n - 2)x^{n}\notag\\
               &= t(0) + t(1)x - t(0)x + \sum_{n \geq 2}(t(n) - t(n - 1) - t(n - 2))x^{n}.
\end{align}
\noindent Because the order of the polynomial on the left hand side of (\ref{recursivet-1}) is two, the coefficients of $x^{n}$ for all $n \geq 3$ must be $0$. Thus, $t(n) - t(n - 1) - t(n - 2) = 0$ implying that
\begin{align}
    t(n) = t(n - 1) + t(n - 2).
\end{align}
\noindent This completes the proof.\qed

\subsection{Diamond Cacti}
First, we may name all the vertices of $D(n)$, the Diamond cacti of $n$ Diamonds, as detailed in Figure \ref{Dn}.

\vskip -0.5 cm
\setlength{\unitlength}{0.8cm}
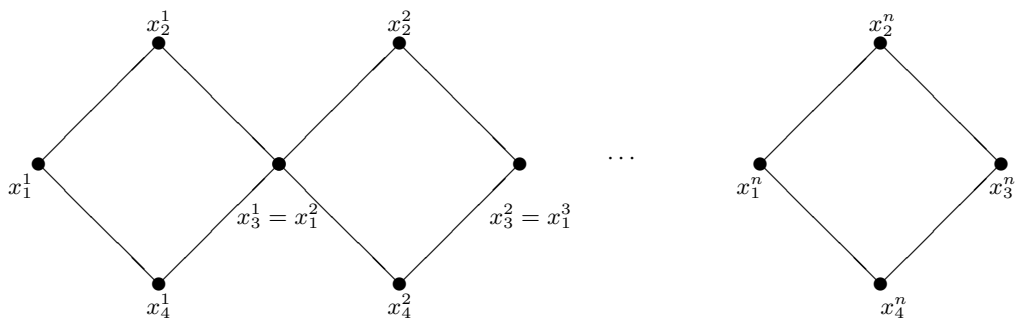
\begin{figure}[htb]
\begin{center}
\begin{picture}(13, 6)

%D1
\put(-2, 2){\circle*{0.2}}
\put(0, 4){\circle*{0.2}}
\put(2, 2){\circle*{0.2}}
\put(0, 0){\circle*{0.2}}

\put(-2, 2){\line(1, 1){2}}
\put(-2, 2){\line(1, -1){2}}
\put(0, 4){\line(1, -1){2}}
\put(0, 0){\line(1, 1){2}}

\put(-2.5, 1.5){\footnotesize$x^{1}_{1}$}
\put(1.3, 1){\footnotesize$x^{1}_{3} = x^{2}_{1}$}
\put(-0.2, 4.2){\footnotesize$x^{1}_{2}$}
\put(-0.2, -0.5){\footnotesize$x^{1}_{4}$}

%D2
\put(2, 2){\circle*{0.2}}
\put(4, 4){\circle*{0.2}}
\put(6, 2){\circle*{0.2}}
\put(4, 0){\circle*{0.2}}

\put(2, 2){\line(1, 1){2}}
\put(2, 2){\line(1, -1){2}}
\put(4, 4){\line(1, -1){2}}
\put(4, 0){\line(1, 1){2}}

\put(5.5, 1){\footnotesize$x^{2}_{3} = x^{3}_{1}$}
\put(3.8, 4.2){\footnotesize$x^{2}_{2}$}
\put(3.8, -0.5){\footnotesize$x^{2}_{4}$}

\put(7.45, 2){\footnotesize$\cdots$}

%Dn
\put(10, 2){\circle*{0.2}}
\put(12, 4){\circle*{0.2}}
\put(14, 2){\circle*{0.2}}
\put(12, 0){\circle*{0.2}}

\put(10, 2){\line(1, 1){2}}
\put(10, 2){\line(1, -1){2}}
\put(12, 4){\line(1, -1){2}}
\put(12, 0){\line(1, 1){2}}

\put(9.6, 1.5){\footnotesize$x^{n}_{1}$}
\put(13.8, 1.5){\footnotesize$x^{n}_{3}$}
\put(11.8, 4.2){\footnotesize$x^{n}_{2}$}
\put(12, -0.5){\footnotesize$x^{n}_{4}$}

\end{picture}
\end{center}
\vskip 0.5 cm
\caption{\label{Dn} A labelled diamond cactus of $n$ diamonds}
\end{figure}
\vskip 5 pt
\vskip 5 pt

\noindent Then, we recall that
\vskip 5 pt

\indent $d(n)$ = the number of all maximal independent sets of $D(n)$

\noindent and

\indent $d(n, k)$ = the number of maximal independent sets containing $k$ vertices of $d(n)$.
\vskip 5 pt

\noindent Thus,
\begin{align}\label{dn}
d(n) = \sum_{k \geq 0}d(n, k).    
\end{align}
\vskip 5 pt

\noindent Further, we let 
\begin{align}
D(x) = \sum_{n \geq 0}d(n)x^{n}\notag
\end{align}
\noindent be the generating function of $d(n)$ and we let
\begin{align}
D(x, y) = \sum_{n \geq 0}\sum_{k \geq 0}d(n, k)x^{n}y^{k}\notag
\end{align}
\noindent be the bi-variate generating function of $d(n, k)$. It is worth noting that, when $y = 1$, we have
\begin{align}\label{txytotnd}
D(x, 1) = \sum_{n \geq 0}(\sum_{k \geq 0}d(n, k)(1)^{k})x^{n} = \sum_{n \geq 0}d(n)x^{n} = D(x).
\end{align}

\noindent Next, we let $\overline{D}(n)$ be constructed from $D(n)$ by joining two vertices to a vertex at distance two from the cut vertex of the $n^{th}$ diamond.
\vskip 5 pt

\vskip -0.25 cm
\setlength{\unitlength}{0.8cm}
\begin{figure}[htb]
\begin{center}
\begin{picture}(13, 6.5)

%D1
\put(-2, 2){\circle*{0.2}}
\put(0, 4){\circle*{0.2}}
\put(2, 2){\circle*{0.2}}
\put(0, 0){\circle*{0.2}}

\put(-2, 2){\line(1, 1){2}}
\put(-2, 2){\line(1, -1){2}}
\put(0, 4){\line(1, -1){2}}
\put(0, 0){\line(1, 1){2}}

%\put(-2.5, 1.5){\footnotesize$x^{1}_{1}$}
%%\put(1.3, 1){\footnotesize$x^{3}_{1} = x^{1}_{2}$}
%\put(-0.2, 4.2){\footnotesize$x^{2}_{1}$}
%\put(-0.2, -0.5){\footnotesize$x^{4}_{1}$}

%D2
\put(2, 2){\circle*{0.2}}
\put(4, 4){\circle*{0.2}}
\put(6, 2){\circle*{0.2}}
\put(4, 0){\circle*{0.2}}

\put(2, 2){\line(1, 1){2}}
\put(2, 2){\line(1, -1){2}}
\put(4, 4){\line(1, -1){2}}
\put(4, 0){\line(1, 1){2}}

%\put(5.5, 1){\footnotesize$x^{3}_{2} = x^{1}_{3}$}
%\put(3.8, 4.2){\footnotesize$x^{2}_{2}$}
%\put(3.8, -0.5){\footnotesize$x^{4}_{2}$}

\put(7.45, 2){\footnotesize$\cdots$}

%Dn
\put(10, 2){\circle*{0.2}}
\put(12, 4){\circle*{0.2}}
\put(14, 2){\circle*{0.2}}
\put(12, 0){\circle*{0.2}}

\put(10, 2){\line(1, 1){2}}
\put(10, 2){\line(1, -1){2}}
\put(12, 4){\line(1, -1){2}}
\put(12, 0){\line(1, 1){2}}

\put(9.6, 1.5){\footnotesize$x^{n}_{1}$}
\put(13.8, 1.5){\footnotesize$x^{n}_{3}$}
\put(11.8, 4.2){\footnotesize$x^{n}_{2}$}
\put(12, -0.5){\footnotesize$x^{n}_{4}$}

\put(16, 4){\circle*{0.2}}
\put(16, 4){\line(-1, -1){2}}
\put(16, 0){\circle*{0.2}}
\put(16, 0){\line(-1, 1){2}}
\put(15.5, 4.3){\footnotesize$x_{n + 1}$}
\put(15.5, -0.5){\footnotesize$x'_{n + 1}$}

\end{picture}
\end{center}
\vskip 0.5 cm
\caption{\label{figdbarn} The graph $\overline{D}(n)$}
\end{figure}
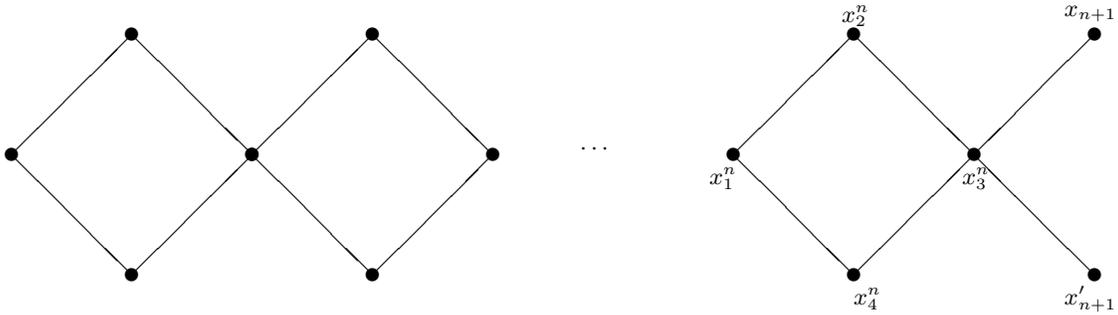
\vskip 5 pt
\vskip 5 pt

\noindent Then, we let
\vskip 5 pt

\indent $\overline{d}(n, k)$ = the number of maximal independent sets containing $k$ vertices of $\overline{D}(n)$.
\vskip 5 pt

\noindent and let
\begin{align}
\overline{D}(x, y) = \sum_{n \geq 0}\sum_{k \geq 0}\overline{d}(n, k)x^{n}y^{k}\notag
\end{align}
\noindent be the bi-variate generating function of $\overline{d}(n, k)$.
\vskip 5 pt

\indent Now, we are ready to prove Theorem \ref{bivariate-d}.
\vskip 5 pt

\noindent \emph{Proof of Theorem \ref{bivariate-d}} First, we will establish a recurrence relation of $d(n, k)$. Let $D$ be a maximal independent set of $D(n)$ containing $k$ vertices. We distinguish two cases. 
\vskip 5 pt

\noindent \textbf{Case 1:}  $x^{n}_{3} \in D$.\\
\indent Because $D$ is independent, $x^{n}_{2}, x^{n}_{4} \notin D$. Removing $x^{n}_{2}, x^{n}_{3}, x^{n}_{4}$ form $D(n)$ results in $D(n-1)$. Thus, $D = D' \cup \{x^{n}_{3}\}$ where  $D'$ is a maximal independent set of $D(n-1)$ containing $k - 1$ vertices.There are $d(n - 1,k - 1)$ possibilities of $D'$ giving that there are $d(n - 1,k - 1)$ possibilities of $D$.
\vskip 5 pt

\noindent \textbf{Case 2:} $x^{n}_{3} \notin D$.\\
\indent By maximality of $D$, $x^{n}_{2} \in D$ or $x^{n}_{4} \in D$. In any case, $x^{n}_{1} \notin D$ where $x^{n}_{2}, x^{n}_{4} \in D$. Removing $x^{n}_{1}, x^{n}_{2}, x^{n}_{3}, x^{n}_{4}$ from $D$ results in $\overline{D}(n - 2)$. Thus, $D = D' \cup \{x^{n}_{2}, x^{n}_{4}\}$ where $D'$ is a maximal independent set of $\overline{D}(n - 2)$. There are $\overline{d}(n - 2, k - 2)$
 possibilities of $D'$ yielding that there are $\overline{d}(n - 2, k - 2)$ possibilities of $D$.
\vskip 5 pt

\indent From Cases 1 and 2, we have that
\begin{align}\label{d11}
d(n, k) = d(n - 1,k - 1) + \overline{d}(n - 2, k - 2)
\end{align}
\noindent For $n \geq 2$ and $k \geq 2$, we multiply $x^{n}y^{k}$ throughout (\ref{d11}) and sum over all $x^{n}y^{k}$. Thus, we have that

\begin{align}\label{d22}
\sum_{n \geq 2}\sum_{k \geq 2}d(n, k)x^{n}y^{k} = \sum_{n \geq 2}\sum_{k \geq 2}d(n-1, k-1)x^{n}y^{k} + \sum_{n \geq 2}\sum_{k \geq 2}\overline{d}(n - 2, k - 2)x^{n}y^{k}
\end{align}

\noindent We first consider the term $\displaystyle{\sum_{n \geq 2}\sum_{k \geq 2}d(n, k)x^{n}y^{k}}$. It can be checked that

\begin{center}
    $d(0, 0) = 1, d(0, 1) = 0, d(0, 2) = 0$, and $ d(0, k) = 0$ for all $k \geq 2$
    \vskip 5 pt
    $d(1, 0) = 0, d(1, 0) = 0, d(1, 2) = 2$, and $ d(0, k) = 0$ for all $k \geq 2$
    \vskip 5 pt
    and $d(n, k) = 0$ for all $n \geq 2, k \leq 1$.
\end{center}

\noindent Clearly, 
\begin{align}\label{d1}
\sum_{n \geq 2}\sum_{k \geq 2}d(n, k)x^{n}y^{k} &= \sum_{n \geq 2}\sum_{k \geq 2}d(n, k)x^{n}y^{k}  + d(0, 0)  + d(1, 2)xy^{2} - d(0, 0) -d(1, 2)xy^{2}\notag\\
                                                &= \sum_{n \geq 0}\sum_{k \geq 0}t(n, k)x^{n}y^{k}-1-2xy^2\notag\\
                                                &= D(x,y)-1-2xy^{2}
\end{align}

\noindent Now, we consider the term $\displaystyle{\sum_{n \geq 2}\sum_{k \geq 2}d(n-1, k-1)x^{n}y^{k}}$. Clearly, 
\begin{align}\label{d2}
\sum_{n \geq 2}\sum_{k \geq 2}d(n-1, k-1)x^{n}y^{k} &= xy(\sum_{n \geq 2}\sum_{k \geq 2}d(n-1, k-1)x^{n-1}y^{k-1})\notag\\
                                                    &= xy(\sum_{n \geq 2}\sum_{k \geq 2}d(n-1, k-1)x^{n-1}y^{k-1}\notag\\ 
                                                    &~~~~ + d(0, 0) - d(0, 0))\notag\\ 
                                                    &= xy(\sum_{n \geq 0}\sum_{k \geq 0}d(n-1, k-1)x^{n-1}y^{k-1} - 1)\notag\\
                                                    &= xyD(x,y)-xy
\end{align}

\noindent Finally, we consider the term $\displaystyle{\sum_{n \geq 2}\sum_{k \geq 2}\overline{d}(n - 2, k - 2)x^{n}y^{k}}$. It can be check that,

\begin{align}\label{dbar1}
\sum_{n \geq 2}\sum_{k \geq 2}\overline{d}(n - 2, k - 2)x^{n}y^{k}
    &= x^{2}y^{2}\sum_{n \geq 2}\sum_{k \geq 2}\overline{d}(n - 2, k - 2)x^{n - 2}y^{k - 2}\notag\\
    &= x^{2}y^{2}\sum_{n \geq 0}\sum_{k \geq 0}\overline{d}(n, k)x^{n}y^{k}\notag\\
    &= x^{2}y^{2}\overline{D}(x, y).
\end{align}

\noindent Plugging (\ref{d1}), (\ref{d2}) and (\ref{dbar1}) to (\ref{d22}), we have 
\begin{align}
D(x,y)-1-2xy^{2} &= xyD(x,y)-xy + x^{2}y^{2}\overline{D}(x, y)\notag
\end{align}

\noindent which can be solved that

\begin{align}\label{semifinald4}
D(x, y) = 1 + 2xy^{2} + xyD(x,y)-xy + x^{2}y^{2}\overline{D}(x, y).
\end{align}
\vskip 5 pt

\indent Next, we will establish the recurrence relation of $\overline{d}(n, k)$. For the graph $\overline{D}(n)$, we let $D$ be a maximal independent set of $\overline{D}(n)$ containing $k$ vertices. There are $2$ cases.
\vskip 5 pt

\noindent \textbf{Case 1:} $x^{n}_{3} \in D$\\
\indent Thus $x^{n}_{2}, x^{n}_{4}, x_{n + 1}, x'_{n + 1} \notin D$. Removing $x^{n}_{2}, x^{n}_{3}, x^{n}_{4}, x_{n + 1}, x'_{n + 1}$ from $\overline{D}(n)$ results in $D(n - 1)$. Thus, $D = D' \cup \{x^{n}_{3}\}$ where $D'$ is a maximal independent set of $D(n - 1)$ containing $k - 1$ vertices. Therefore, there are $d(n - 1, k - 1)$ possibilities of $D'$ yielding that there are $d(n - 1, k - 1)$ possibilities of $D$.
\vskip 5 pt

\noindent \textbf{Case 2:} $x^{n}_{3} \notin D$\\
\indent By maximality of $D$, $x_{n + 1}, x'_{n + 1} \in D$. Removing $x^{n}_{3}, x_{n + 1}, x'_{n + 1}$ from $\overline{D}(n)$ results in $\overline{D}(n - 1)$. Thus, $D = D' \cup \{x_{n + 1}, x'_{n + 1}\}$ where $D'$ is a maximal independent set of $\overline{D}(n - 1)$ containing $k - 2$ vertices. Therefore, there are $\overline{d}(n - 1, k - 2)$ possibilities of $D'$ yielding that there are $\overline{d}(n - 1, k - 2)$ possibilities of $D$.
\vskip 5 pt

\indent From Cases 1 and 2, we have that
\begin{align}\label{dbar4}
    \overline{d}(n, k) = d(n - 1, k - 1) + \overline{d}(n - 1, k - 2).
\end{align}
\vskip 5 pt

\noindent For $n \geq 1$ and $k \geq 2$, we multiply $x^{n}y^{k}$ throughout (\ref{dbar4}) and sum over all $x^{n}y^{k}$. Thus, we have that

\begin{align}\label{dbar444}
\sum_{n \geq 1}\sum_{k \geq 2}\overline{d}(n , k )x^{n}y^{k} = \sum_{n \geq 1}\sum_{k \geq 2}d(n-1 , k-1 )x^{n}y^{k} + \sum_{n \geq 1}\sum_{k \geq 2}\overline{d}(n - 1, k - 2)x^{n}y^{k}
\end{align}

\noindent We first consider the term $\displaystyle{\sum_{n \geq 1}\sum_{k \geq 2}\overline{d}(n , k )x^{n}y^{k}}$. It can be checked that

 \begin{center}
    $\overline{d}(0, 0) = 0, \overline{d}(0, 1) = 1 ,\overline{d}(0, 2) = 1$ and $\overline{d}(0, k) = 0$ for all $k \geq 3$
    \vskip 5 pt
\end{center}

\begin{align}\label{dbar2}
\sum_{n \geq 1}\sum_{k \geq 2}\overline{d}(n , k )x^{n}y^{k}
    &= \sum_{n \geq 1}\sum_{k \geq 2}\overline{d}(n , k)x^{n}y^{k}\notag\\
    & + \overline{d}(0, 1)y  + \overline{d}(0, 2)y^{2} - \overline{d}(0, 1)y -\overline{d}(0, 2)y^{2}\notag\\
    &= \sum_{n \geq 0}\sum_{k \geq 0}\overline{d}(n, k)x^{n}y^{k} -y-y^{2}\notag\\
    &= \overline{D}(x, y)-y-y^{2}
\end{align}    
\vskip 5 pt

\noindent We next consider the term $\displaystyle{\sum_{n \geq 1}\sum_{k \geq 2}d(n-1 , k-1 )x^{n}y^{k}}$. Recall that

\begin{center}
    $d(0, 0) = 1$, and $d(n, 0) = 0$ for all $n \geq 1$.
    \vskip 5 pt
\end{center}

\begin{align}\label{d3}
\sum_{n \geq 1}\sum_{k \geq 2}d(n-1 , k-1 )x^{n}y^{k} &= xy(\sum_{n \geq 1}\sum_{k \geq 2}d(n-1, k-1)x^{n-1}y^{k-1})\notag\\
                                                &= xy(\sum_{n \geq 0}\sum_{k \geq 1}d(n, k)x^{n}y^{k})\notag\\
                                                &= xy(\sum_{n \geq 0}\sum_{k \geq 1}d(n, k)x^{n}y^{k}+ d(0, 0) - d(0, 0))\notag\\
                                                &= xy(\sum_{n \geq 0}\sum_{k \geq 0}d(n, k)x^{n}y^{k}-1)\notag\\
                                                &= xyD(x,y)-xy\notag\\
\end{align}
\vskip 5 pt

\noindent Finally, we consider the term $\displaystyle{\sum_{n \geq 1}\sum_{k \geq 2}\overline{d}(n - 1, k - 2)x^{n}y^{k}}$. Clearly, 

\begin{align}\label{dbar3}
\sum_{n \geq 1}\sum_{k \geq 2}\overline{d}(n - 1, k - 2)x^{n}y^{k}
    &= xy^{2}(\sum_{n \geq 1}\sum_{k \geq 2}\overline{d}(n - 1, k - 2)x^{n - 1}y^{k - 2})\notag\\
    &=  xy^{2}(\sum_{n \geq 0}\sum_{k \geq 0}\overline{d}(n, k)x^{n}y^{k} )\notag\\
    &=  xy^{2}\overline{D}(x, y)
\end{align}    
\vskip 5 pt

\noindent Plugging (\ref{dbar2}), (\ref{d3}) and (\ref{dbar3}) to (\ref{dbar444}), we have 

\begin{align}
\overline{D}(x, y)-y-y^{2} &= xyD(x,y)-xy + xy^{2}\overline{D}(x, y).\notag
\end{align}

\noindent which can be solved that

\begin{align}\label{semifinaldbar4}
\overline{D}(x, y) &= y + y^{2} + xyD(x,y)-xy + xy^{2}\overline{D}(x, y).
\end{align}
\vskip 5 pt

\noindent By plugging (\ref{semifinald4}) to (\ref{semifinaldbar4}), we have

\begin{align}
D(x, y) = \frac{1 + xy^{2} -xy - x^{2}y^{4} - x^{3}y^{3} + 2x^{2}y^{3}}{1 - xy - xy^{2} + x^{2}y^{3} - x^{3}y^{3}}.\notag
\end{align}

\noindent as required. This proves Theorem \ref{bivariate-d}.
\qed

\noindent \emph{Proof of Theorem \ref{recurrence-d}} By Theorem \ref{bivariate-d} with $y = 1$, we have that
\begin{align}
\sum_{n \geq 0}d(n)x^{n} &= D(x, 1)\notag\\
                         &= \frac{1 + x^{2} - x^{3}}{1 - 2x + x^{2} - x^{3}}
\end{align}

\noindent which can  be solved that

\begin{align}\label{recursived-1}
1 + x^{2} - x^{3} &= (1 - 2x + x^{2} - x^{3})\sum_{n \geq 0}d(n)x^{n}\notag\\
                  &= \sum_{n \geq 0}d(n)x^{n}  - \sum_{n \geq 0}2d(n)x^{n + 1} + \sum_{n \geq 0}d(n)x^{n + 2} - \sum_{n \geq 0}d(n)x^{n + 3}\notag\\
                  &= \sum_{n \geq 0}d(n)x^{n} - \sum_{n \geq 1}2d(n - 1)x^{n} + \sum_{n \geq 2}d(n - 2)x^{n} - \sum_{n \geq 3}2d(n - 3)x^{n}\notag\\
                  &= d(0) + d(1)x + d(2)x^{2} + d(3)x^{3} - 2d(0)x - 2d(1)x^{2} - 2d(2)x^{3} + d(0)x^{2}\notag\\
                  & +d(1)x^{3} - 2d(0)x^{3} + \sum_{n \geq 4}(d(n) - 2d(n) - d(n - 2) - d(n - 3))x^{n}
\end{align}
\noindent Because the order of the polynomial on the left hand side of (\ref{recursived-1}) is three, the coefficients of $x^{n}$ for all $n \geq 4$ must be $0$. Thus, $d(n) - 2d(n - 1) + d(n - 2)- d(n - 3) = 0$ implying that
\begin{align*}
    d(n) = 2d(n - 1) - d(n - 2) + d(n - 3).
\end{align*}
This proves Theorem \ref{recurrence-d}. \qed

\subsection{Square Cacti}
First we may name all the vertices of $S(n)$, a square cacti of $n$ squares, as detailed in Figure \ref{Sn}. 
\vskip 5 pt

\vskip -0.25 cm
\setlength{\unitlength}{0.8cm}
\begin{figure}[htb]
\begin{center}
\begin{picture}(13, 6.5)

%S1
\put(-2, 3){\circle*{0.2}}
\put(-2, 6){\circle*{0.2}}
\put(1, 3){\circle*{0.2}}
\put(1, 6){\circle*{0.2}}

\put(-2, 3){\line(0, 1){3}}
\put(-2, 6){\line(1, 0){3}}
\put(1, 3){\line(-1, 0){3}}
\put(1, 6){\line(0, -1){3}}

\put(-2.5, 6.2){\footnotesize$x^{1}_{1}$}
\put(1, 6.2){\footnotesize$x^{1}_{2}$}
\put(1.2, 3.3){\footnotesize$x^{1}_{3} = x^{2}_{1}$}
\put(-2.5, 2.5){\footnotesize$x^{1}_{4}$}

%S2
\put(1, 3){\circle*{0.2}}
\put(4, 3){\circle*{0.2}}
\put(4, 0){\circle*{0.2}}
\put(1, 0){\circle*{0.2}}

\put(1, 3){\line(1, 0){3}}
\put(4, 3){\line(0, -1){3}}
\put(4, 0){\line(-1, 0){3}}
\put(1, 0){\line(0, 1){3}}

\put(4, 3.2){\footnotesize$x^{2}_{2}  = x^{3}_{4}$}
\put(4.3, 0){\footnotesize$x^{2}_{3}$}
\put(0.6, -0.6){\footnotesize$x^{2}_{4}$}

\put(7, 3){\footnotesize$\cdots$}

%Sn
\put(10, 3){\circle*{0.2}}
\put(10, 6){\circle*{0.2}}
\put(13, 3){\circle*{0.2}}
\put(13, 6){\circle*{0.2}}

\put(10, 3){\line(0, 1){3}}
\put(10, 6){\line(1, 0){3}}
\put(13, 3){\line(-1, 0){3}}
\put(13, 6){\line(0, -1){3}}

\put(9.5, 6.2){\footnotesize$x^{n}_{1}$}
\put(13, 6.2){\footnotesize$x^{n}_{2}$}
\put(13, 2.6){\footnotesize$x^{n}_{3}$}
\put(9.5, 2.6){\footnotesize$x^{n}_{4}$}

\end{picture}
\end{center}
\vskip 0.5 cm
\caption{A labelled square cacti of $n$ Squares}
\label{Sn}
\end{figure}
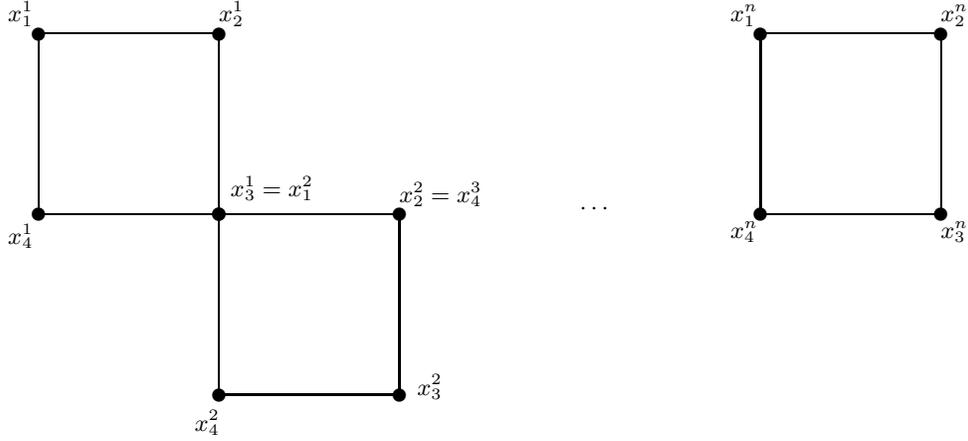

\indent Then, we let
\vskip 5 pt

\indent $s(n)$ = the number of all maximal independent sets of $S(n)$

\noindent and

\indent $s(n, k)$ = the number of maximal independent sets containing $k$ vertices of $S(n)$.
\vskip 5 pt

\noindent Thus,
\begin{align*}
s(n) = \sum_{k \geq 0}s(n, k)    
\end{align*}
\vskip 5 pt

\noindent Further, we let 
\begin{align}
S(x) = \sum_{n \geq 0}s(n)x^{n}\notag
\end{align}
\noindent be the generating function of $s(n)$ and we let
\begin{align}
S(x, y) = \sum_{n \geq 0}\sum_{k \geq 0}s(n, k)x^{n}y^{k}\notag
\end{align}
\noindent be the bi-variate generating function of $s(n, k)$. It is worth noting that, when $y = 1$, we have
\begin{align}\label{txytotns}
S(x, 1) = \sum_{n \geq 0}(\sum_{k \geq 0}s(n, k)(1)^{k})x^{n} = \sum_{n \geq 0}s(n)x^{n} = S(x).
\end{align}

\vskip -0.25 cm
\setlength{\unitlength}{0.7cm}
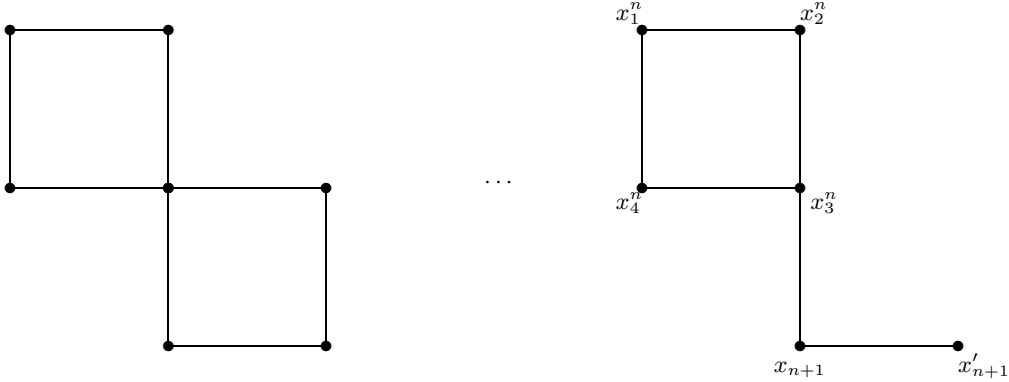
\begin{figure}[htb]
\begin{center}
\begin{picture}(13, 6.5)

%S1
\put(-2, 3){\circle*{0.2}}
\put(-2, 6){\circle*{0.2}}
\put(1, 3){\circle*{0.2}}
\put(1, 6){\circle*{0.2}}

\put(-2, 3){\line(0, 1){3}}
\put(-2, 6){\line(1, 0){3}}
\put(1, 3){\line(-1, 0){3}}
\put(1, 6){\line(0, -1){3}}

%\put(-2.5, 6.2){\footnotesize$x^{1}_{1}$}
%\put(1, 6.2){\footnotesize$x^{2}_{1}$}
%\put(1.2, 3.3){\footnotesize$x^{3}_{1} = x^{1}_{2}$}
%\put(-2.5, 2.5){\footnotesize$x^{4}_{1}$}

%S2
\put(1, 3){\circle*{0.2}}
\put(4, 3){\circle*{0.2}}
\put(4, 0){\circle*{0.2}}
\put(1, 0){\circle*{0.2}}

\put(1, 3){\line(1, 0){3}}
\put(4, 3){\line(0, -1){3}}
\put(4, 0){\line(-1, 0){3}}
\put(1, 0){\line(0, 1){3}}

%\put(4, 3.2){\footnotesize$x^{2}_{2}$}
%\put(4.3, 0){\footnotesize$x^{3}_{2} = x^{1}_{3}$}
%\put(0.6, -0.6){\footnotesize$x^{4}_{2}$}

\put(7, 3){\footnotesize$\cdots$}

%Sn
\put(10, 3){\circle*{0.2}}
\put(10, 6){\circle*{0.2}}
\put(13, 3){\circle*{0.2}}
\put(13, 6){\circle*{0.2}}

\put(10, 3){\line(0, 1){3}}
\put(10, 6){\line(1, 0){3}}
\put(13, 3){\line(-1, 0){3}}
\put(13, 6){\line(0, -1){3}}

\put(9.5, 6.2){\footnotesize$x^{n}_{1}$}
\put(13, 6.2){\footnotesize$x^{n}_{2}$}
\put(13.2, 2.6){\footnotesize$x^{n}_{3}$}
\put(9.5, 2.6){\footnotesize$x^{n}_{4}$}

\put(13, 0){\circle*{0.2}}
\put(16, 0){\circle*{0.2}}
\put(13, 0){\line(0, 1){3}}
\put(16, 0){\line(-1, 0){3}}
\put(12.5, -0.5){\footnotesize$x_{n+1}$}
\put(16, -0.5){\footnotesize$x'_{n+1}$}

\end{picture}
\end{center}
\vskip 1 cm
\caption{The graph $\overline{S}(n)$.}
\label{Snbar}
\end{figure}
\vskip 5 pt

\noindent Next, we let $\overline{S}(n)$ be constructed from $S(n)$ by joining an end vertex of a path of length one to a vertex at distance one from the cut vertex of the $n^{th}$ square. The graph $\overline{S}(n)$ is shown in Figure \ref{Snbar}.
\vskip 5 pt

\noindent Then, we let 
\vskip 5 pt

\indent $\overline{s}(n, k)$ = the number of maximal independent sets containing $k$ vertices of $\overline{S}(n)$.
\vskip 5 pt

\noindent and let
\begin{align}
\overline{S}(x, y) = \sum_{n \geq 0}\sum_{k \geq 0}\overline{s}(n, k)x^{n}y^{k}\notag
\end{align}
\noindent be the bi-variate generating function of $\overline{t}(n, k)$.
\vskip 5 pt

\indent Now, we are ready to prove Theorem \ref{bivariate-s}.
\vskip 5 pt

\noindent \emph{Proof of Theorem \ref{bivariate-s}} First, we will establish a recurrence relation of $s(n, k)$. Let $D$ be a maximal independent set of $S(n)$ containing $k$ vertices. We distinguish two cases. 
\vskip 5 pt

\noindent \textbf{Case 1:}  $x^{n}_{2} \in D$.\\
\indent Because $D$ is independent, $x^{n}_{1}, x^{n}_{3} \notin D$. Removing $x^{n}_{1}, x^{n}_{2}, x^{n}_{3}$ form $D(n)$ results in $S(n-1)$. Thus, $D = D' \cup \{x^{n}_{2}\}$ where  $D'$ is a maximal independent set of $S(n-1)$ containing $k - 1$ vertices. There are $s(n - 1,k - 1)$ possibilities of $D'$ giving that there are $s(n - 1,k - 1)$ possibilities of $D$.
\vskip 5 pt

\noindent \textbf{Case 2:} $x^{n}_{2} \notin D$.\\
\indent By maximality of $D$, $x^{n}_{1} \in D$ or $x^{n}_{3} \in D$. We may assume that $|D \cap \{x^{n}_{1}, x^{n}_{3}\}| = 1$. Suppose without loss of generality that $x^{n}_{1} \in D$ but $x^{n}_{3} \notin D$. Because $D$ is independent, $x^{n}_{4} \notin D$. Thus, $D \cup \{x^{n}_{3}\}$ is an independent set contradicting the maximality of $D$. Thus, $|D \cap \{x^{n}_{1}, x^{n}_{3}\}| = 2$. Because $D$ is independent, $x^{n}_{4} \notin D$. Removing $x^{n}_{1}, x^{n}_{2}, x^{n}_{3}, x^{n}_{4}$ from $D$ results in $\overline{S}(n - 2)$. Thus, $D = D' \cup \{x^{n}_{1}, x^{n}_{3}\}$ where $D'$ is a maximal independent set of $\overline{S}(n - 2)$ containing $k - 2$ vertices. There are $\overline{s}(n - 2, k - 2)$ possibilities of $D'$ yielding that there are $\overline{s}(n - 2, k - 2)$ possibilities of $D$.
\vskip 5 pt

\indent From Cases 1 and 2, we have that
\begin{align}\label{s1}
s(n, k) = s(n - 1,k - 1) + \overline{s}(n - 2, k - 2)
\end{align}
\noindent For $n \geq 2$ and $k \geq 2$, we multiply $x^{n}y^{k}$ throughout (\ref{s1}) and sum over all $x^{n}y^{k}$. Thus, we have that
\begin{align}\label{s2}
\sum_{n \geq 2}\sum_{k \geq 2}s(n, k)x^{n}y^{k} = \sum_{n \geq 2}\sum_{k \geq 2}s(n-1, k-1)x^{n}y^{k} + \sum_{n \geq 2}\sum_{k \geq 2}\overline{s}(n - 2, k - 2)x^{n}y^{k}
\end{align}

\noindent We first consider the term $\displaystyle{\sum_{n \geq 2}\sum_{k \geq 2}s(n, k)x^{n}y^{k}}$. It can be checked that 

\begin{center}
    $s(0, 0) = 1, s(0, 1) = 0 $ and $s(0, k) = 0 $ for all $k \geq 2$
    \vskip 5 pt
    $s(1, 0) = 0 , s(1, 1) = 0 , s(1, 2) = 2$ and $s(0,k) = 0$ for all $k \geq 3$
    \vskip 5 pt
    $s(n, k) = 0$ for all $n \geq 2, k \leq 1$.
\end{center}

\noindent Thus,
\begin{align}\label{s3}
\sum_{n \geq 2}\sum_{k \geq 2}s(n, k)x^{n}y^{k} &= \sum_{n \geq 2}\sum_{k \geq 2}s(n, k)x^{n}y^{k}  + s(0, 0)  + s(1, 2)xy^{2} - s(0, 0) - s(1, 2)xy^{2}\notag\\
                                                &= \sum_{n \geq 0}\sum_{k \geq 0}s(n, k)x^{n}y^{k}-1-2xy^2\notag\\
                                                &= S(x,y)-1-2xy^{2}
\end{align}
\vskip 5 pt

\noindent We next consider the term $\displaystyle{\sum_{n \geq 2}\sum_{k \geq 2}s(n-1, k-1)x^{n}y^{k}}$. Recall that
\vskip 5 pt

\begin{align}\label{s4}
\sum_{n \geq 2}\sum_{k \geq 2}s(n-1, k-1)x^{n}y^{k} &= xy(\sum_{n \geq 2}\sum_{k \geq 2}s(n-1, k-1)x^{n-1}y^{k-1})\notag\\
                                                    &= xy(\sum_{n \geq 1}\sum_{k \geq 1}s(n, k)x^{n}y^{k} + s(0, 0) - s(0, 0))\notag\\  
                                                     &= xy(\sum_{n \geq 0}\sum_{k \geq 0}s(n, k)x^{n}y^{k}-1)\notag\\ 
                                                    &= xyS(x,y)-xy
\end{align}

\noindent Finally, we consider the term $\displaystyle{\sum_{n \geq 2}\sum_{k \geq 2}\overline{s}(n - 2, k - 2)x^{n}y^{k}}$. Clearly,
\vskip 5 pt

\begin{align}\label{sbar1}
\sum_{n \geq 2}\sum_{k \geq 2}\overline{s}(n - 2, k - 2)x^{n}y^{k}
    &= x^{2}y^{2}(\sum_{n \geq 2}\sum_{k \geq 2}\overline{s}(n - 2, k - 2)x^{n - 2}y^{k - 2})\notag\\
    &= x^{2}y^{2}(\sum_{n \geq 0}\sum_{k \geq 0}\overline{s}(n , k )x^{n}y^{k})\notag\\
    &=  x^{2}y^{2}\overline{S}(x, y)
\end{align}

\noindent Plugging (\ref{s3}), (\ref{s4}) and (\ref{sbar1}) to (\ref{s2}), we have 

\begin{align}
S(x,y)-1-2xy^{2} &= xyS(x,y)-xy + x^{2}y^{2}\overline{S}(x, y)\notag
\end{align}

\noindent which can be solved that

\begin{align}\label{semifinals4}
S(x,y) &= 1+ 2xy^{2} + xyS(x,y)-xy + x^{2}y^{2}\overline{S}(x, y).
\end{align}
\vskip 5 pt

\indent Next, we will establish the recurrence relation of $\overline{s}(n, k)$. For the graph $\overline{S}(n)$, we let $D$ be a maximal independent set of $\overline{S}(n)$ containing $k$ vertices. There are $2$ cases.
\vskip 5 pt

\noindent \textbf{Case 1:} $x'_{n + 1} \in D$\\
\indent Thus $ x_{n + 1} \notin D$. Removing $x_{n + 1}, x'_{n + 1}$ from $\overline{S}(n)$ results in $S(n)$. Thus, $D = D' \cup \{x'_{n+1}\}$ where $D'$ is a maximal independent set of $S(n)$ containing $k - 1$ vertices. Therefore, there are $s(n , k - 1)$ possibilities of $D'$ yielding that there are $s(n , k - 1)$ possibilities of $D$.
\vskip 5 pt

\noindent \textbf{Case 2:} $x'_{n + 1} \notin D$\\
\indent By maximality of $D$, $x_{n + 1} \in D$. Because $D$ is independent, $x^{n}_{3} \notin D$. Removing $x^{n}_{3}, x_{n + 1}, x'_{n + 1}$ from $\overline{S}(n)$ results in $\overline{S}(n - 1)$. Thus, $D = D' \cup \{x_{n + 1}\}$ where $D'$ is a maximal independent set of $\overline{s}(n - 1)$ containing $k - 1$ vertices. Therefore, there are $\overline{s}(n - 1, k - 1)$ possibilities of $D'$ yielding that there are $\overline{s}(n - 1, k - 1)$ possibilities of $D$.
\vskip 5 pt

\indent From Cases 1 and 2, we have that
\begin{align}\label{sbar4}
    \overline{s}(n, k) = s(n,k-1) + \overline{s}(n - 1, k - 1).
\end{align}
\vskip 5 pt

\noindent For $n \geq 1$ and $k \geq 1$, we multiply $x^{n}y^{k}$ throughout (\ref{sbar4}) and sum over all $x^{n}y^{k}$. Thus, we have that

\begin{align}\label{sbar444}
\sum_{n \geq 1}\sum_{k \geq 1}\overline{s}(n, k) = \sum_{n \geq 1}\sum_{k \geq 1}s(n,k-1) + \sum_{n \geq 1}\sum_{k \geq 1}\overline{s}(n - 1, k - 1).
\end{align}

\noindent We first consider the term $\displaystyle{\sum_{n \geq 1}\sum_{k \geq 1}\overline{s}(n , k)x^{n}y^{k}}$. 
\vskip 5 pt

\noindent Next, it can be checked that
 \begin{center}
    $\overline{s}(0, 0) = 0, \overline{s}(0, 1) = 1 ,\overline{s}(0, 2) = 1$ and $\overline{s}(0, k) = 0$ for all $k \geq 3$
    \vskip 5 pt
    $\overline{s}(1, 0) = 0, \overline{s}(1, 1) = 0, \overline{s}(1, 2) = 1$ and $\overline{s}(1, k) = 0$ for all $k \geq 3$.
    \vskip 5 pt
\end{center}

\noindent Thus, 

\begin{align}\label{sbar2}
\sum_{n \geq 1}\sum_{k \geq 1}\overline{s}(n , k)x^{n}y^{k}
    &= \sum_{n \geq 1}\sum_{k \geq 1}\overline{s}(n , k )x^{n}y^{k})\notag\\
    & + s(0,1)y + s(0,2)y^{2} - s(0,1)y - s(0,2)y^{2}\notag\\
    &=  \overline{S}(x, y)-y-y^{2}
\end{align}
\vskip 5 pt

\noindent Now, we consider the term $\displaystyle{\sum_{n \geq 1}\sum_{k \geq 1}s(n,k-1)x^{n}y^{k}}$. Clearly,

\begin{center}
    $s(0, 0) = 1, s(0, 1) = 0 $ and $s(0, k) = 0 $ for all $k \geq 2$
    \vskip 5 pt
\end{center}

\begin{align}\label{s5}
\sum_{n \geq 1}\sum_{k \geq 1}s(n, k-1)x^{n}y^{k} &= y(\sum_{n \geq 1}\sum_{k \geq 1}s(n, k-1)x^{n}y^{k-1} + s(0, 0) - s(0, 0))\notag\\
                                                &= y\sum_{n \geq 0}\sum_{k \geq 0}s(n, k)x^{n}y^{k}-y\notag\\
                                                &= yS(x,y)-y\notag\\
\end{align}
\vskip 5 pt

\noindent Finally, we consider the term $\displaystyle{\sum_{n \geq 1}\sum_{k \geq 1}\overline{s}(n - 1, k - 1)x^{n}y^{k}}$. It can be check that,
\vskip 5 pt

\begin{align}\label{sbar3}
\sum_{n \geq 1}\sum_{k \geq 1}\overline{s}(n - 1, k - 1)x^{n}y^{k}
    &= xy(\sum_{n \geq 1}\sum_{k \geq 1}\overline{s}(n - 1, k - 1)x^{n - 1}y^{k - 1})\notag\\
    &= xy(\sum_{n \geq 0}\sum_{k \geq 0}\overline{s}(n , k )x^{n}y^{k})\notag\\
    &=  xy\overline{S}(x, y)
\end{align} 

\noindent Plugging (\ref{sbar2}), (\ref{s5}) and (\ref{sbar3}) to (\ref{sbar444}), we have 

\begin{align}
\overline{S}(x, y)-y-y^{2} &= yS(x,y)-y + xy\overline{S}(x, y).\notag
\end{align}

\noindent This together with (\ref{semifinals4}) yield

\begin{align}
S(x, y) =  \frac{1 - 2xy + 2xy^{2} - 2x^{2}y^{3} + x^{2}y^{2} + x^{2}y^{4}}{1 - 2xy + x^{2}y^{2} - x^{2}y^{3}}\notag
\end{align}

\noindent as required. This proves Theorem \ref{bivariate-s}.\qed

\noindent \emph{Proof of Theorem \ref{recurrence-s}}
By Theorem \ref{bivariate-s} with $y = 1$, we have that 
\begin{align}
\sum_{n \geq 0}s(n)x^{n} &= (x, 1)\notag\\
                         &= \frac{1}{1 - 2x}
\end{align}

\noindent which can  be solved that

\begin{align}\label{recursives-1}
1  &= (1 - 2x )\sum_{n \geq 0}s(n)x^{n}\notag\\
   &= \sum_{n \geq 0}s(n)x^{n}  - \sum_{n \geq 0}2s(n)x^{n + 1}\notag\\
   &= \sum_{n \geq 0}s(n)x^{n} - \sum_{n \geq 1}2s(n - 1)x^{n}\notag\\
   &= s(0) +  \sum_{n \geq 1}(s(n) - 2s(n - 1))x^{n}
\end{align}
\noindent Because the order of the polynomial on the left hand side of (\ref{recursives-1}) is zero, the coefficients of $x^{n}$ for all $n \geq 1$ must be $0$. Thus, $s(n) - 2s(n - 1) = 0$ implying that

\begin{align}\label{recursives-2}
    s(n) = 2s(n - 1). 
\end{align}
\noindent This proves Theorem \ref{recurrence-s}.\qed

\subsection{Pentagonal Cacti}
\begin{figure}[H]
\centering
\definecolor{ududff}{rgb}{0.30196078431372547,0.30196078431372547,1}
\resizebox{1\textwidth}{!}{%

\begin{tikzpicture}[line cap=round,line join=round,>=triangle 45,x=1cm,y=1cm]
\clip(0.07084031875793423,-1.7138790246542182) rectangle (9.988466422567084,3.9142216784276616);
\draw [line width=0.3pt] (1,1)-- (3,1);
\draw [line width=0.3pt] (1,1)-- (0.8111130238479834,1.9874845007960271);
\draw [line width=0.3pt] (0.8111130238479834,1.9874845007960271)-- (1.510822840987893,2.382972658309889);
\draw [line width=0.3pt] (1.510822840987893,2.382972658309889)-- (2.1902512141527324,1.9874845007960271);
\draw [line width=0.3pt] (2.1902512141527324,1.9874845007960271)-- (1.8150445006139406,0);
\draw [line width=0.3pt] (3,1)-- (3.204323412906225,0);
\draw [line width=0.3pt] (1.8150445006139406,0)-- (2.49447287377878,-0.41586661024974864);
\draw [line width=0.3pt] (2.49447287377878,-0.41586661024974864)-- (3.204323412906225,0);
\draw [line width=0.3pt] (5,1)-- (6,1);
\draw [line width=0.3pt] (5,1)-- (4.8065574869367405,1.997625222783562);
\draw [line width=0.3pt] (4.8065574869367405,1.997625222783562)-- (5.4859858601015805,2.4032541022849587);
\draw [line width=0.3pt] (5.4859858601015805,2.4032541022849587)-- (6.195836399229025,2.007765944771097);
\draw [line width=0.3pt] (6.195836399229025,2.007765944771097)-- (6,1);
\begin{scriptsize}
\draw [fill=black] (1,1) circle (1pt);
\draw [fill=black] (2,1) circle (1pt);
\draw [fill=black] (3,1) circle (1pt);
\draw [fill=black] (0.8111130238479834,1.9874845007960271) circle (1pt);
\draw [fill=black] (2.1902512141527324,1.9874845007960271) circle (1pt);
\draw [fill=black] (1.510822840987893,2.382972658309889) circle (1pt);
\draw [fill=black] (1.8150445006139406,0) circle (1pt);
\draw [fill=black] (3.204323412906225,0) circle (1pt);
\draw [fill=black] (2.49447287377878,-0.41586661024974864) circle (1pt);
\draw [fill=black] (5,1) circle (1pt);
\draw [fill=black] (4.8065574869367405,1.997625222783562) circle (1pt);
\draw [fill=black] (5.4859858601015805,2.4032541022849587) circle (1pt);
\draw [fill=black] (6.195836399229025,2.007765944771097) circle (1pt);
\draw [fill=black] (6,1) circle (1pt);

\draw[color=black] (4,1) node {$...$};

%%%%Pent1
\draw[color=black] (0.7,1) node {\tiny$x^{1}_{1}$};

\draw[color=black] (0.5,2) node {\tiny$x^{1}_{2}$};

\draw[color=black] (1.6,2.7) node {\tiny$x^{1}_{3}$};

\draw[color=black] (2.5,2) node {\tiny$x^{1}_{4}$};

\draw[color=black] (1.8,1.2) node {\tiny$x^{1}_{5}$};

%%%%Pent2
\draw[color=black] (3.3,1) node {\tiny$x^{2}_{5}$};

\draw[color=black] (1.6,0) node {\tiny$x^{2}_{2}$};

\draw[color=black] (2.5,-0.7) node {\tiny$x^{2}_{3}$};

\draw[color=black] (3.5,0) node {\tiny$x^{2}_{4}$};

\draw[color=black] (2.2,0.7) node {\tiny$x^{2}_{1}$};

%%%pentn
\draw[color=black] (4.7,1) node {\tiny$x^{n}_{1}$};

\draw[color=black] (4.5,2) node {\tiny$x^{n}_{2}$};

\draw[color=black] (5.6,2.7) node {\tiny$x^{n}_{3}$};

\draw[color=black] (6.5,2) node {\tiny$x^{n}_{4}$};

\draw[color=black] (6.3,1) node {\tiny$x^{n}_{5}$};

\end{scriptsize}
\end{tikzpicture}

 }%
\vskip -1 cm
\caption{A labelled pentagonal cactus of $n$ pentagons.}
\label{pentwithlabel}
\end{figure}
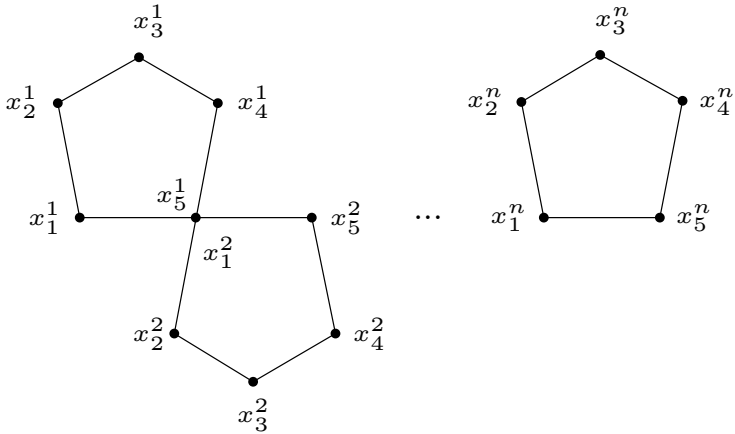

First, we name all the vertices of $P(n)$ as shown in Figure \ref{pentwithlabel}

\vskip -1 cm
\vskip 5 pt

\noindent Then, recall that
\vskip 5 pt

\indent $p(n)$ = the number of all maximal independent sets of $P(n)$

\noindent and

\indent $p(n, k)$ = the number of maximal independent sets containing $k$ vertices of $p(n)$.
\vskip 5 pt

\noindent Thus,
\begin{align}\label{pn}
p(n) = \sum_{k \geq 0}p(n, k).    
\end{align}
\vskip 5 pt

\noindent Further, we let 
\begin{align}
P(x) = \sum_{n \geq 0}p(n)x^{n}\notag
\end{align}
\noindent be the generating function of $p(n)$ and we let
\begin{align}
P(x, y) = \sum_{n \geq 0}\sum_{k \geq 0}p(n, k)x^{n}y^{k}\notag
\end{align}
\noindent be the bi-variate generating function of $p(n, k)$. It is worth noting that, when $y = 1$, we have
\begin{align}\label{txytotnp}
P(x, 1) = \sum_{n \geq 0}(\sum_{k \geq 0}p(n, k)(1)^{k})x^{n} = \sum_{n \geq 0}p(n)x^{n} = P(x).
\end{align}

%\vskip 1 cm
\begin{figure}[H]
\centering
\definecolor{ududff}{rgb}{0.30196078431372547,0.30196078431372547,1}
\resizebox{1\textwidth}{!}{%

\begin{tikzpicture}[line cap=round,line join=round,>=triangle 45,x=1cm,y=1cm]
\clip(-0.10155195503015937,-2.1600707921057554) rectangle (9.816074148778982,3.4680299109761292);
\draw [line width=0.3pt] (1,1)-- (3,1);
\draw [line width=0.3pt] (1,1)-- (0.8111130238479834,1.9874845007960271);
\draw [line width=0.3pt] (0.8111130238479834,1.9874845007960271)-- (1.510822840987893,2.382972658309889);
\draw [line width=0.3pt] (1.510822840987893,2.382972658309889)-- (2.1902512141527324,1.9874845007960271);
\draw [line width=0.3pt] (2.1902512141527324,1.9874845007960271)-- (1.8150445006139406,0);
\draw [line width=0.3pt] (3,1)-- (3.204323412906225,0);
\draw [line width=0.3pt] (1.8150445006139406,0)-- (2.49447287377878,-0.41586661024974864);
\draw [line width=0.3pt] (2.49447287377878,-0.41586661024974864)-- (3.204323412906225,0);
\draw [line width=0.3pt] (4.148179353047069,1.0608443319252094)-- (5.148179353047069,1.0608443319252094);
\draw [line width=0.3pt] (4.148179353047069,1.0608443319252094)-- (3.954736839983809,2.0584695547087715);
\draw [line width=0.3pt] (3.954736839983809,2.0584695547087715)-- (4.634165213148649,2.464098434210168);
\draw [line width=0.3pt] (4.634165213148649,2.464098434210168)-- (5.344015752276094,2.0686102766963064);
\draw [line width=0.3pt] (5.344015752276094,2.0686102766963064)-- (5.148179353047069,1.0608443319252094);
\draw [line width=0.3pt] (5.148179353047069,1.0608443319252094)-- (4.928249353047069,0.07088433192520943);
\draw [line width=0.3pt] (4.928249353047069,0.07088433192520943)-- (5.658378133889677,-0.3448815563370043);
\draw [line width=0.3pt] (5.658378133889677,-0.3448815563370043)-- (6.358087951029586,0.0810287671394623);
\begin{scriptsize}
\draw [fill=black] (1,1) circle (1pt);
\draw [fill=black] (2,1) circle (1pt);
\draw [fill=black] (3,1) circle (1pt);
\draw [fill=black] (0.8111130238479834,1.9874845007960271) circle (1pt);
\draw [fill=black] (2.1902512141527324,1.9874845007960271) circle (1pt);
\draw [fill=black] (1.510822840987893,2.382972658309889) circle (1pt);
\draw [fill=black] (1.8150445006139406,0) circle (1pt);
\draw [fill=black] (3.204323412906225,0) circle (1pt);
\draw [fill=black] (2.49447287377878,-0.41586661024974864) circle (1pt);
\draw [fill=black] (4.148179353047069,1.0608443319252094) circle (1pt);
\draw [fill=black] (3.954736839983809,2.0584695547087715) circle (1pt);
\draw [fill=black] (4.634165213148649,2.464098434210168) circle (1pt);
\draw [fill=black] (5.344015752276094,2.0686102766963064) circle (1pt);
\draw [fill=black] (5.148179353047069,1.0608443319252094) circle (1pt);
\draw [fill=black] (4.928249353047069,0.07088433192520943) circle (1pt);
\draw [fill=black] (5.658378133889677,-0.3448815563370043) circle (1pt);
\draw [fill=black] (6.358087951029586,0.0810287671394623) circle (1pt);

%%%pentn
\draw[color=black] (3.8,1) node {\tiny$x^{n}_{1}$};

\draw[color=black] (3.6,2) node {\tiny$x^{n}_{2}$};

\draw[color=black] (4.7,2.7) node {\tiny$x^{n}_{3}$};

\draw[color=black] (5.6,2) node {\tiny$x^{n}_{4}$};

\draw[color=black] (6,1) node {\tiny$x^{n}_{5}=x^{n+1}_{1}$};

%%pentn+1

\draw[color=black] (4.5,0) node {\tiny$x^{n+1}_{2}$};

\draw[color=black] (5.8,-0.7) node {\tiny$x^{n+1}_{3}$};

\draw[color=black] (6.8,0) node {\tiny$x^{n+1}_{4}$};

\end{scriptsize}
\end{tikzpicture}

 }%
\vskip -1.5 cm
\caption{The graphs $\overline{P}(n)$}
\label{pentbar}
\end{figure}
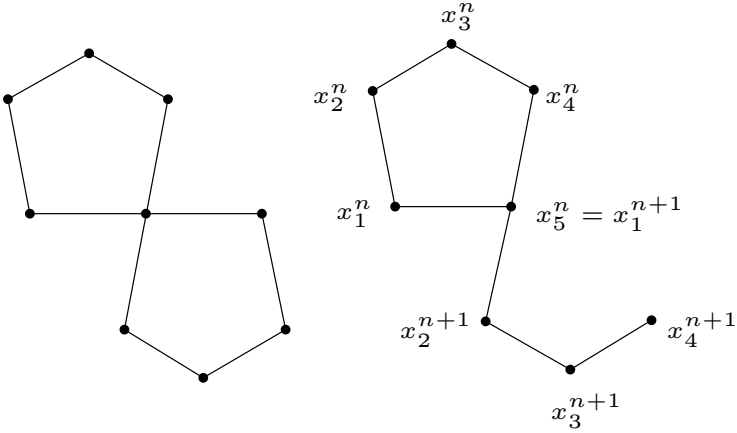
\vskip 5 pt

\noindent Next, we let $\overline{P}(n)$ be constructed from $p(n)$ by joining an end vertex of a path of length two to a vertex at distance one from the cut vertex of the $n^{th}$ pentagon. The graph $\overline{P}(n)$ is shown in Figure \ref{pentbar}.
\vskip 5 pt

\noindent Then, we let
\vskip 5 pt

\indent $\overline{p}(n, k)$ = the number of maximal independent sets containing $k$ vertices of $\overline{P}(n)$.
\vskip 5 pt

\noindent and let
\begin{align}
\overline{P}(x, y) = \sum_{n \geq 0}\sum_{k \geq 0}\overline{p}(n, k)x^{n}y^{k}\notag
\end{align}
\noindent be the bi-variate generating function of $\overline{p}(n, k)$. Now, we are ready to prove Theorem \ref{bivariate-p}
\vskip 5 pt

\noindent \emph{Proof of Theorem \ref{bivariate-p}} Let $D$ be a maximal independent set of $P(n)$ containing $k$ vertices. We distinguish $2$ cases. 
\vskip 5 pt

\noindent \textbf{Case 1:} $x^{n}_{4} \in D$\\
\indent Thus, $x^{n}_{3}, x^{n}_{5} \notin D$. We further consider $2$ subcases.
\vskip 5 pt

\noindent \textbf{Subcase 1.1:} $x^{n}_{2} \in D$\\
\indent Thus, $x^{n}_{1} \notin D$. Removing $x^{n}_{1}, ...,  x^{n}_{5}$ from $P(n)$ results in $\overline{P}(n - 2)$. Thus, $D = D' \cup \{x^{n}_{2}, x^{n}_{4}\}$ where $D'$ is a maximal independent set of $\overline{P}(n - 2)$ containing $k - 2$ vertices. Therefore, there are $\overline{p}(n - 2, k - 2)$ possibilities of $D'$ yielding that there are $\overline{p}(n - 2, k - 2)$ possibilities of $D$.
\vskip 5 pt

\noindent \textbf{Subcase 1.2:} $x^{n}_{2} \notin D$\\
\indent By maximality of $D$, $x^{n}_{1} \in D$. Thus, $x^{n-1}_{1}, x^{n-1}_{4} \notin D$. By maximality of $D$, either $x^{n-1}_{2} \in D$ or $x^{n-1}_{3} \in D$. Removing all vertices of the $n^{th}$ and $(n-1)^{th}$ pentagons results in $\overline{P}(n - 3)$. Thus, for any $u \in \{x^{n-1}_{2}, x^{n-1}_{3}\}$, we have that $D = D' \cup \{u, x^{n}_{1}, x^{n}_{4}\}$ $D'$ is a maximal independent set of $\overline{P}(n - 3)$ containing $k - 3$ vertices. Therefore, there are $\overline{p}(n - 3, k - 3)$ possibilities of $D'$ yielding that there are $2\overline{p}(n - 3, k - 3)$ possibilities of $D$.
\vskip 5 pt

\noindent \textbf{Case 2:} $x^{n}_{4} \notin D$\\
\indent We further have the following $3$ subcases.
\vskip 5 pt

\noindent \textbf{Subcase 2.1:} $n^{n}_{3} \notin D$ and $x^{n}_{5} \in D$.\\
\indent Thus, $x^{n}_{2} \in D$. Removing all the vertices of the $n^{th}$ pentagon results in $\overline{P}(n - 2)$. We have that there are $\overline{p}(n - 2, k - 2)$ possibilities of $D$.
\vskip 5 pt

\noindent \textbf{Subcase 2.2:} $n^{n}_{3} \in D$ and $x^{n}_{5} \notin D$.\\
\indent Thus, $x^{n}_{2} \notin D$. By the maximality of $D$, $x^{n}_{1} \in D$ and this implies that either $x^{n - 1}_{2} \in D$ or $x^{n - 1}_{3} \in D$. Removing all the vertices of the $(n - 1)^{th}$ and $n^{th}$ pentagons results in $\overline{P}(n - 3)$. We have that there are $2\overline{p}(n - 3, k - 3)$ possibilities of $D$.
\vskip 5 pt

\noindent \textbf{Subcase 2.3:} $n^{n}_{3} \in D$ and $x^{n}_{5} \in D$.\\
\indent Thus, $x^{n}_{1}, x^{n}_{2} \notin D$. Removing all the vertices of the $n^{th}$ pentagon results in $\overline{P}(n - 2)$. We have that there are $\overline{p}(n - 2, k - 2)$ possibilities of $D$.
\vskip 5 pt

\noindent From, all the cases, we have that
\begin{align}\label{p11}
p(n, k) = 3\overline{p}(n - 2, k - 2) + 4\overline{p}(n - 3, k - 3)
\end{align}
\noindent For $n \geq 3 $ and $k \geq 3 $, we multiply $x^{n}y^{k}$ throughout (\ref{p11}) and sum over all $x^{n}y^{k}$. Thus, we have that
\begin{align}\label{p22}
\sum_{n \geq 3}\sum_{k \geq 3}p(n, k)x^{n}y^{k} = \sum_{n \geq 3}\sum_{k \geq 3}3\overline{p}(n-2, k-2)x^{n}y^{k} + \sum_{n \geq 3}\sum_{k \geq 3}4\overline{p}(n - 3, k - 3)x^{n}y^{k}
\end{align}

\noindent We first consider the term $\displaystyle{\sum_{n \geq 3}\sum_{k \geq 3}p(n, k)x^{n}y^{k}}$. It can be checked that

\begin{center}
    $p(0, 0) = 1$, and $ p(0, k) = 0$ for all $k \geq 1$
    \vskip 5 pt
    $p(1, 0) = 0, p(1, 1) = 0, p(1, 2) = 5$, and $ p(1, k) = 0$ for all $k \geq 3$
    \vskip 5 pt
    $p(2, 0) = p(2, 1) = p(2, 2) = 0, p(2, 3) = 4, p(2, 4) = 9$, and $p(2, k) = 0$ for all $k \geq 5$.
\end{center}
\noindent Thus, 
\begin{align}\label{p1}
\sum_{n \geq 3}\sum_{k \geq 3}p(n, k)x^{n}y^{k} &= \sum_{n \geq 3}\sum_{k \geq 3}p(n, k)x^{n}y^{k}  + p(0, 0)  + p(1, 2)xy^{2} + p(2,3)x^{2}y^{3} + p(2, 4)x^{2}y^{4}\notag\\
                                                &-p(0, 0) -p(1, 2)xy^{2} - p(2,3)x^{2}y^{3} - p(2, 4)x^{2}y^{4}\notag\\
                                                &= \sum_{n \geq 0}\sum_{k \geq 0}p(n, k)x^{n}y^{k}-1-5xy^2 - 4x^2y^3 - 9x^{2}y^{4}\notag\\
                                                &= P(x,y)-1-5xy^2 - 4x^2y^3 - 9x^{2}y^{4.}
\end{align}

\noindent Now, we consider the term $\displaystyle{\sum_{n \geq 3}\sum_{k \geq 3}3\overline{p}(n - 2, k - 2)x^{n}y^{k}}$. clearly,

\begin{center}
    $\overline{p}(0, 0) = 0, \overline{p}(0, 1) = 0 ,\overline{p}(0, 2) = 3$ and $\overline{p}(0, k) = 0$ for all $k \geq 3$
    \vskip 5 pt
    $\overline{p}(1, 0) = 0, \overline{p}(1, 1) = 0 ,\overline{p}(1, 2) = 0, \overline{p}(1, 3) = 7$, $\overline{p}(1, 4) = 3$, and $\overline{p}(1, k) = 0$ for all $k \geq 5$
    \vskip 5 pt
    $\overline{p}(n, 0) = 0$ for all $n \geq 0$.
\end{center}
\noindent Thus,
\begin{align}\label{pbar1}
\sum_{n \geq 3}\sum_{k \geq 3}3\overline{p}(n - 2, k - 2)x^{n}y^{k}
    &= 3x^{2}y^{2}\sum_{n \geq 3}\sum_{k \geq 3}\overline{p}(n - 2, k - 2)x^{n - 2}y^{k - 2}\notag\\
    &= 3x^{2}y^{2}(\sum_{n \geq 1}\sum_{k \geq 1}\overline{p}(n, k)x^{n}y^{k} + \overline{p}(0,2)y^{2} - \overline{p}(0,2)y^{2})\notag\\
    &= 3x^{2}y^{2}(\sum_{n \geq 0}\sum_{k \geq 0}\overline{p}(n, k)x^{n}y^{k} - 3y^{2})\notag\\
    &= 3x^{2}y^{2}\overline{P}(x, y) - 9x^{2}y^{4}.
\end{align}

\noindent Finally, we consider the term $\displaystyle{\sum_{n \geq 3}\sum_{k \geq 3}4\overline{p}(n - 3, k - 3)x^{n}y^{k}}$. It can be check that,

\begin{align}\label{pbar2}
\sum_{n \geq 3}\sum_{k \geq 3}4\overline{p}(n - 3, k - 3)x^{n}y^{k}
    &= 4x^{3}y^{3}\sum_{n \geq 3}\sum_{k \geq 3}\overline{p}(n - 3, k - 3)x^{n - 3}y^{k - 3}\notag\\
    &= 4x^{3}y^{3}\sum_{n \geq 0}\sum_{k \geq 0}\overline{p}(n, k)x^{n}y^{k}\notag\\
    &= 4x^{3}y^{3}\overline{P}(x, y).
\end{align}

\noindent Plugging (\ref{p1}), (\ref{pbar1}) and (\ref{pbar2}) to (\ref{p22}), we have 
\begin{align}
P(x,y) - 1- 5xy^2 - 4x^{2}y^{3} -9x^{2}y^{4} = 3x^{2}y^{2}\overline{P}(x, y) - 9x^{2}y^{4} + 4x^{3}y^{3}\overline{P}(x, y)\notag
\end{align}

\noindent which can be solved that

\begin{align}\label{semifinalp4}
P(x, y) = 1 + 5xy^2 + 4x^{2}y^{3} + 3x^{2}y^{2}\overline{P}(x, y) + 4x^{3}y^{3}\overline{P}(x, y).
\end{align}
\vskip 5 pt

\indent Next, we will find a recurrence relation of $\overline{p}(n, k)$. Let $D$ be a maximal independent set of $\overline{P}(n)$ containing $k$ vertices. We distinguish $2$ cases.
\vskip 5 pt

\noindent \textbf{Case 1:} $x^{n + 1}_{3} \in D$\\
\indent Thus, $x^{n + 1}_{2}, x^{n + 1}_{4} \notin D$. Removing $x^{n + 1}_{2}, x^{n + 1}_{3}, x^{n + 1}_{4}$ from $\overline{P}(n)$ results in $P(n)$. Thus, there are $p(n, k - 1)$ psiibilities of $D$.
\vskip 5 pt

\noindent \textbf{Case 2:} $x^{n + 1}_{3} \notin D$\\
\indent By maximality of $D$,  $x^{n+1}_{4} \in D$. We further distinguish $2$ subcases.
\vskip 5 pt

\noindent \textbf{Case 2.1:} $x^{n + 1}_{2} \notin D$\\
\indent By maximality of $D$, $x^{n}_{5} \in D$. Thus, $x^{n}_{4}, x^{n}_{1} \notin D$. By maximality of $D$, either $x^{n}_{2} \in D$ or $x^{n}_{3} \in D$. Removing all vertices of the $n^{th}$ pentagon and $x^{n + 1}_{2}, x^{n + 1}_{3}, x^{n + 1}_{4}$ from $\overline{P}(n)$ results in $\overline{P}(n - 2)$. Thus, for any $u \in \{x^{n}_{2}, x^{n}_{3}\}$, we have that $D = D' \cup \{u, x^{n + 1}_{4}, x^{n}_{5}\}$ where $D'$ is a maximal independent set of $\overline{P}(n - 2)$ containing $k - 3$ vertices. Therefore, there are $\overline{p}(n - 2, k - 3)$ possibilities of $D'$ yielding that there are $2\overline{p}(n - 2, k - 3)$ possibilities of $D$.  
\vskip 5 pt

\noindent \textbf{Case 2.2:} $x^{n + 1}_{2} \in D$\\
\indent Thus, $x^{n}_{5} \notin D$. Removing $x^{n}_{5}, x^{n + 1}_{2}, x^{n + 1}_{3}, x^{n + 1}_{4}$ from $\overline{P}(n)$ results in $\overline{P}(n - 1)$. There are $\overline{p}(n - 1, k - 2)$ possibilities of $D$.
\vskip 5 pt

\indent From all the cases, we have that
\begin{align}\label{pbar3}
    \overline{p}(n, k) = p(n, k - 1) + 2\overline{p}(n - 2, k - 3) + \overline{p}(n - 1, k - 2).
\end{align}
\vskip 5 pt

\noindent For $n \geq 2, k \geq 3$, we multiply $x^{n}y^{k}$ throughout (\ref{pbar3}) and sum over all $n, k$. We have that

\begin{align}\label{pbar33}
    \displaystyle{
    \sum_{n \geq 2}\sum_{k \geq 3}\overline{p}(n, k) = \sum_{n \geq 2}\sum_{k \geq 3}p(n, k - 1) + \sum_{n \geq 2}\sum_{k \geq 3}2\overline{p}(n - 2, k - 3) + \sum_{n \geq 2}\sum_{k \geq 3}\overline{p}(n - 1, k - 2).}
\end{align}

\noindent We first consider the term $\displaystyle{\sum_{n \geq 2}\sum_{k \geq 3}\overline{p}(n , k )x^{n}y^{k}}$. It can be checked that,

  \begin{center}
    $\overline{p}(0, 0) = 0, \overline{p}(0, 1) = 0 ,\overline{p}(0, 2) = 3$ and $\overline{p}(0, k) = 0$ for all $k \geq 3$
    \vskip 5 pt
    $\overline{p}(1, 0) = 0, \overline{p}(1, 1) = 0 ,\overline{p}(1, 2) = 0, \overline{p}(1, 3) = 7$, $\overline{p}(1, 4) = 3$, and $\overline{p}(1, k) = 0$ for all $k \geq 5$
    \vskip 5 pt
    $\overline{p}(n, k) = 0$ for all $n \geq 2, k \leq 2$.
\end{center}

\begin{align}\label{pbar4}
\sum_{n \geq 2}\sum_{k \geq 3}\overline{p}(n , k )x^{n}y^{k}
    &= \sum_{n \geq 2}\sum_{k \geq 3}\overline{p}(n , k )x^{n}y^{k} + \overline{p}(0,2)y^{2} + \overline{p}(1,3)xy^{3} + \overline{p}(1, 4)xy^4\notag\\
    & - \overline{p}(0,2)y^{2} - \overline{p}(1,3)xy^{3} - \overline{p}(1, 4)xy^4 \notag\\
    &= \sum_{n \geq 0}\sum_{k \geq 0}\overline{p}(n , k )x^{n}y^{k} - 3y^{2} - 7xy^{3} - 3xy^4\notag\\
    &= \overline{P}(x, y) - 3y^{2} - 7xy^{3} - 3xy^4
\end{align}    
\vskip 5 pt

\noindent We next consider the term $\displaystyle{\sum_{n \geq 2}\sum_{k \geq 3}p(n , k-1 )x^{n}y^{k}}$. Recall that

\begin{align}\label{p2}
\sum_{n \geq 2}\sum_{k \geq 3}p(n , k-1 )x^{n}y^{k} &= y(\sum_{n \geq 2}\sum_{k \geq 3}p(n , k-1)x^{n}y^{k-1})\notag\\
                                                               &= y(\sum_{n \geq 2}\sum_{k \geq 2}p(n, k)x^{n}y^{k})\notag\\
                                                               &= y(\sum_{n \geq 2}\sum_{k \geq 2}p(n, k)x^{n}y^{k} + p(0,0) + p(1,2)xy^{2} - p(0,0) - p(1,2)xy^{2})\notag\\ 
                                                               &= y(\sum_{n \geq 0}\sum_{k \geq 0}p(n, k)x^{n}y^{k} - 1 - 5xy^{2})\notag\\
                                                               &= yP(x,y) - y - 5xy^{3}
\end{align}
\vskip 5 pt

\noindent We next consider the term $\displaystyle{\sum_{n \geq 2}\sum_{k \geq 3}2\overline{p}(n - 2 , k - 3)x^{n}y^{k}}$. It can be checked that,

\begin{align}\label{pbar5}
\sum_{n \geq 2}\sum_{k \geq 3}2\overline{p}(n - 2 , k - 3 )x^{n}y^{k}
    &= 2x^{2}y^{3}(\sum_{n \geq 2}\sum_{k \geq 3}\overline{p}(n - 2, k - 3)x^{n-2}y^{k-3}) \notag\\
    &= 2x^{2}y^{3}(\sum_{n \geq 0}\sum_{k \geq 0}\overline{p}(n , k )x^{n}y^{k})\notag\\
    &= 2x^{2}y^{3}\overline{P}(x, y)
\end{align}    
\vskip 5 pt

\noindent Finally, we consider the term $\displaystyle{\sum_{n \geq 2}\sum_{k \geq 3}\overline{p}(n - 1, k - 2)x^{n}y^{k}}$. Clearly, 

\begin{align}\label{pbar6}
\sum_{n \geq 2}\sum_{k \geq 3}\overline{p}(n - 1, k - 2)x^{n}y^{k}
    &= xy^{2}(\sum_{n \geq 2}\sum_{k \geq 3}\overline{p}(n - 1, k - 2)x^{n - 1}y^{k - 2})\notag\\
    &=  xy^{2}(\sum_{n \geq 1}\sum_{k \geq 1}\overline{p}(n, k)x^{n}y^{k} )\notag\\
    &=  xy^{2}(\sum_{n \geq 1}\sum_{k \geq 1}\overline{p}(n, k)x^{n}y^{k} + \overline{p}(0,2)y^{2} - \overline{p}(0,2)y^{2})\notag\\
    &=  xy^{2}(\sum_{n \geq 0}\sum_{k \geq 0}\overline{p}(n, k)x^{n}y^{k} - 3y^{2})\notag\\
    &=  xy^{2}\overline{P}(x, y) - 3xy^{4}
\end{align}    
\vskip 5 pt

\noindent Plugging (\ref{pbar4}), (\ref{p2}) , (\ref{pbar5}) and (\ref{pbar6}) to (\ref{pbar33}), we have 

\begin{align}
\overline{P}(x, y) - 3y^{2} - 7xy^{3} - 3xy^4 = yP(x,y) - y - 5xy^{3} + 2x^{2}y^{3}\overline{P}(x, y) + xy^{2}\overline{P}(x, y) - 3xy^{4}.\notag
\end{align}

\noindent which can be solved that

\begin{align}\label{semifinalpbar7}
\overline{P}(x, y) &= 3y^{2} + 7xy^{3} + yP(x,y) - y - 5xy^{3} + 2x^{2}y^{3}\overline{P}(x, y) + xy^{2}\overline{P}(x, y).
\end{align}
\vskip 5 pt

\noindent By (\ref{semifinalp4}) and (\ref{semifinalpbar7}), we have
\begin{align}
P(x, y) = \frac{1 + 4x^3y^5 - 4x^3y^4 + 4x^2y^4 - x^2y^3 + 4xy^2}{1 - 5x^{2}y^{3} - xy^{2} - 4x^{3}y^{4}}\notag
\end{align}
\noindent as required.
\qed

\noindent We will prove Theorem \ref{recurrence-p}.
\vskip 5 pt

\noindent \emph{Proof of Theorem \ref{recurrence-p}} By Theorem \ref{bivariate-p} with $y = 1$, we have that
\vskip 5 pt

\begin{align}
\sum_{n \geq 0}p(n)x^{n} &= P(x, 1)\notag\\
                         &= \frac{3x^{2} + 4x + 1}{ - 4x^{3} - 5x^{2} - x + 1}
\end{align}

\noindent which can  be solved that

\begin{align}\label{recursived-p1}
3x^{2} + 4x + 1 &= (- 4x^{3} - 5x^{2} - x + 1)\sum_{n \geq 0}p(n)x^{n}\notag\\
                  &= \sum_{n \geq 0}p(n)x^{n}  - \sum_{n \geq 0}5p(n)x^{n + 2} - \sum_{n \geq 0}p(n)x^{n + 1} - \sum_{n \geq 0}4p(n)x^{n + 3}\notag\\
                  &= \sum_{n \geq 0}p(n)x^{n} - \sum_{n \geq 2}5p(n - 2)x^{n} - \sum_{n \geq 1}p(n - 1)x^{n} - \sum_{n \geq 3}4p(n - 3)x^{n}\notag\\
                  &= p(0) + p(1)x + p(2)x^{2} - 5p(0)x^{2} - p(0)x - p(1)x^{2}\notag\\
                  & + \sum_{n \geq 3}(p(n) - 5p(n-2) - p(n - 1) - 4p(n - 3))x^{n}
\end{align}
\noindent Because the order of the polynomial on the left hand side of (\ref{recursived-p1}) is two, the coefficients of $x^{n}$ for all $n \geq 3$ must be $0$. Thus, $p(n) - 5p(n-2) - p(n - 1) - 4p(n - 3) = 0$ implying that

\begin{align*}
    p(n) = p(n - 1) + 5p(n - 2) + 4p(n - 3)
\end{align*}

\subsection{Meta-Pentagonal Cacti}
\begin{figure}[H]
\centering
\definecolor{ududff}{rgb}{0.30196078431372547,0.30196078431372547,1}
\resizebox{1\textwidth}{!}{%

\begin{tikzpicture}[line cap=round,line join=round,>=triangle 45,x=1cm,y=1cm]
\clip(-0.8215432161451385,-2.798936277320456) rectangle (9.096082887664004,2.829164425761429);
\draw [line width=0.3pt] (0.6247932864612085,-0.997722231544381)-- (0.4359063103091919,-0.010237730748353324);
\draw [line width=0.3pt] (0.4359063103091919,-0.010237730748353324)-- (1.1356161274491015,0.38525042676550875);
\draw [line width=0.3pt] (1.1356161274491015,0.38525042676550875)-- (1.8150445006139408,-0.010237730748353324);
\draw [line width=0.3pt] (3,1)-- (3.204323412906225,0);
\draw [line width=0.3pt] (1.8150445006139406,0)-- (2.49447287377878,-0.41586661024974864);
\draw [line width=0.3pt] (2.49447287377878,-0.41586661024974864)-- (3.204323412906225,0);
\draw [line width=0.3pt] (4.441491495385481,-0.9814692730960854)-- (5.441491495385481,-0.9814692730960854);
\draw [line width=0.3pt] (4.441491495385481,-0.9814692730960854)-- (4.248048982322222,0.016155949687477897);
\draw [line width=0.3pt] (4.248048982322222,0.016155949687477897)-- (4.927477355487062,0.42178482918887417);
\draw [line width=0.3pt] (4.927477355487062,0.42178482918887417)-- (5.637327894614506,0.026296671675012803);
\draw [line width=0.3pt] (5.637327894614506,0.026296671675012803)-- (5.441491495385481,-0.9814692730960854);
\draw [line width=0.3pt] (0.6247932864612085,-0.997722231544381)-- (1.6247932864612085,-0.997722231544381);
\draw [line width=0.3pt] (1.8150445006139408,-0.010237730748353324)-- (1.6247932864612085,-0.997722231544381);
\draw [line width=0.3pt] (1.8150445006139406,0)-- (2,1);
\draw [line width=0.3pt] (2,1)-- (3,1);
\begin{scriptsize}
\draw [fill=black] (0.6247932864612085,-0.997722231544381) circle (1pt);
\draw [fill=black] (1.6247932864612085,-0.997722231544381) circle (1pt);
\draw [fill=black] (3,1) circle (1pt);
\draw [fill=black] (0.4359063103091919,-0.010237730748353324) circle (1pt);
\draw [fill=black] (1.8150445006139408,-0.010237730748353324) circle (1pt);
\draw [fill=black] (1.1356161274491015,0.38525042676550875) circle (1pt);
\draw [fill=black] (1.8150445006139406,0) circle (1pt);
\draw [fill=black] (3.204323412906225,0) circle (1pt);
\draw [fill=black] (2.49447287377878,-0.41586661024974864) circle (1pt);
\draw [fill=black] (4.441491495385481,-0.9814692730960854) circle (1pt);
\draw [fill=black] (4.248048982322222,0.016155949687477897) circle (1pt);
\draw [fill=black] (4.927477355487062,0.42178482918887417) circle (1pt);
\draw [fill=black] (5.637327894614506,0.026296671675012803) circle (1pt);
\draw [fill=black] (5.441491495385481,-0.9814692730960854) circle (1pt);
\draw [fill=black] (2,1) circle (1pt);

\draw[color=black] (3.9,0) node {$...$};

%%%%Pent1
\draw[color=black] (1,0.7) node {\tiny$x^{1}_{2}$};

\draw[color=black] (0.2,0.1) node {\tiny$x^{1}_{1}$};

\draw[color=black] (1.8,-1.2) node {\tiny$x^{1}_{4}$};

\draw[color=black] (0.3,-1) node {\tiny$x^{1}_{5}$};

\draw[color=black] (1.4,-0.2) node {\tiny$x^{1}_{3}$};

%%%%Pent2
\draw[color=black] (1.8,1.2) node {\tiny$x^{2}_{5}$};

\draw[color=black] (3.3,1) node {\tiny$x^{2}_{4}$};

\draw[color=black] (2.5,-0.7) node {\tiny$x^{2}_{2}$};

\draw[color=black] (3.5,0) node {\tiny$x^{2}_{3}$};

\draw[color=black] (2.1,0.2) node {\tiny$x^{2}_{1}$};

%%%%Pentn
\draw[color=black] (4.9,0.7) node {\tiny$x^{n}_{2}$};

\draw[color=black] (4.6,-0.2) node {\tiny$x^{n}_{1}$};

\draw[color=black] (5.7,-1.2) node {\tiny$x^{n}_{4}$};

\draw[color=black] (4.3,-1.2) node {\tiny$x^{n}_{5}$};

\draw[color=black] (6,0.1) node {\tiny$x^{n}_{3}$};

\end{scriptsize}
\end{tikzpicture}

 }%
\vskip -1.5 cm
\caption{A labelled meta-pentagonal cactus of $n$ pentagons.}
\label{metapentwithlabel}
\end{figure}
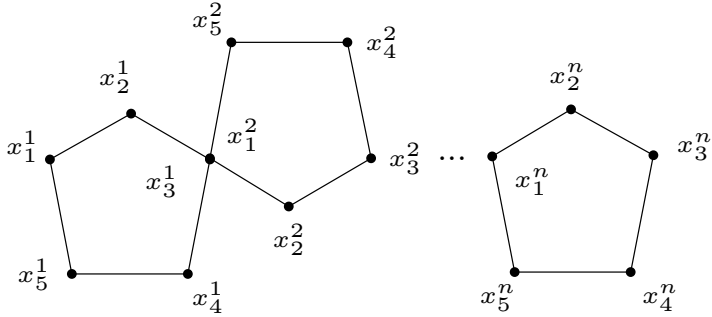

First, we many name all the vertices of $M(n)$ as shown in Figure \ref{metapentwithlabel}. 

\vskip -2 cm
\vskip 5 pt

\noindent Then, we recall that:
\vskip 5 pt

\indent $m(n)$ = the number of all maximal independent sets of $M(n)$

\noindent and

\indent $M(n, k)$ = the number of maximal independent sets containing $k$ vertices of $M(n)$.
\vskip 5 pt

\noindent Thus,
\begin{align}\label{Mn}
m(n) = \sum_{k \geq 0}m(n, k).    
\end{align}
\vskip 5 pt

\noindent Further, we let 
\begin{align}
M(x) = \sum_{n \geq 0}m(n)x^{n}\notag
\end{align}
\noindent be the generating function of $m(n)$ and we let
\begin{align}
M(x, y) = \sum_{n \geq 0}\sum_{k \geq 0}m(n, k)x^{n}y^{k}\notag
\end{align}
\noindent be the bi-variate generating function of $m(n, k)$. It is worth noting that, when $y = 1$, we have
\begin{align}\label{txytotnm}
M(x, 1) = \sum_{n \geq 0}(\sum_{k \geq 0}m(n, k)(1)^{k})x^{n} = \sum_{n \geq 0}m(n)x^{n} = M(x).
\end{align}

\begin{figure}[H]
\centering
\definecolor{ududff}{rgb}{0.30196078431372547,0.30196078431372547,1}
\resizebox{1.2\textwidth}{!}{%

\begin{tikzpicture}[line cap=round,line join=round,>=triangle 45,x=1cm,y=1cm]
\clip(-1.600981787405596,-2.7897522087931304) rectangle (10.42567829554255,3.986462987196297);
\draw [line width=0.3pt] (-0.19861937101705318,0.9910168657903746)-- (0.48573348818055523,1.4061161410413825);
\draw [line width=0.3pt] (0.48573348818055523,1.4061161410413825)-- (1.181305246709273,0.9910168657903746);
\draw [line width=0.3pt] (1.181305246709273,0.9910168657903746)-- (1.9105337032313154,0.5759175905393669);
\draw [line width=0.3pt] (1.9105337032313154,0.5759175905393669)-- (2.5836676630978146,1.002235765121483);
\draw [line width=0.3pt] (2.5836676630978146,1.002235765121483)-- (2.3817274751378643,2.0007178055901234);
\draw [line width=0.3pt] (2.3817274751378643,2.0007178055901234)-- (1.394464334000332,1.9894989062590158);
\draw [line width=0.3pt] (1.394464334000332,1.9894989062590158)-- (1.181305246709273,0.9910168657903746);
\draw [line width=0.3pt] (-0.19861937101705318,0.9910168657903746)-- (-0.00970093979975073,0);
\draw [line width=0.3pt] (-0.00970093979975073,0)-- (0.9902990602002507,0);
\draw [line width=0.3pt] (0.9902990602002507,0)-- (1.181305246709273,0.9910168657903746);
\draw [line width=0.3pt] (-0.19861937101705318,0.9910168657903746)-- (0.48573348818055523,1.4061161410413825);
\draw [line width=0.3pt] (0.48573348818055523,1.4061161410413825)-- (1.181305246709273,0.9910168657903746);
\draw [line width=0.3pt] (1.181305246709273,0.9910168657903746)-- (1.9105337032313154,0.5759175905393669);
\draw [line width=0.3pt] (1.9105337032313154,0.5759175905393669)-- (2.5836676630978146,1.002235765121483);
\draw [line width=0.3pt] (2.5836676630978146,1.002235765121483)-- (2.3817274751378643,2.0007178055901234);
\draw [line width=0.3pt] (2.3817274751378643,2.0007178055901234)-- (1.394464334000332,1.9894989062590158);
\draw [line width=0.3pt] (1.394464334000332,1.9894989062590158)-- (1.181305246709273,0.9910168657903746);
\draw [line width=0.3pt] (-0.19861937101705318,0.9910168657903746)-- (0.48573348818055523,1.4061161410413825);
\draw [line width=0.3pt] (0.48573348818055523,1.4061161410413825)-- (1.181305246709273,0.9910168657903746);
\draw [line width=0.3pt] (1.181305246709273,0.9910168657903746)-- (1.9105337032313154,0.5759175905393669);
\draw [line width=0.3pt] (1.9105337032313154,0.5759175905393669)-- (2.5836676630978146,1.002235765121483);
\draw [line width=0.3pt] (2.5836676630978146,1.002235765121483)-- (2.3817274751378643,2.0007178055901234);
\draw [line width=0.3pt] (2.3817274751378643,2.0007178055901234)-- (1.394464334000332,1.9894989062590158);
\draw [line width=0.3pt] (1.394464334000332,1.9894989062590158)-- (1.181305246709273,0.9910168657903746);
\draw [line width=0.3pt] (-0.19861937101705318,0.9910168657903746)-- (0.48573348818055523,1.4061161410413825);
\draw [line width=0.3pt] (0.48573348818055523,1.4061161410413825)-- (1.181305246709273,0.9910168657903746);
\draw [line width=0.3pt] (1.181305246709273,0.9910168657903746)-- (1.9105337032313154,0.5759175905393669);
\draw [line width=0.3pt] (1.9105337032313154,0.5759175905393669)-- (2.5836676630978146,1.002235765121483);
\draw [line width=0.3pt] (2.5836676630978146,1.002235765121483)-- (2.3817274751378643,2.0007178055901234);
\draw [line width=0.3pt] (2.3817274751378643,2.0007178055901234)-- (1.394464334000332,1.9894989062590158);
\draw [line width=0.3pt] (1.394464334000332,1.9894989062590158)-- (1.181305246709273,0.9910168657903746);
\draw [line width=0.3pt] (3.795308790857518,0.9910168657903748)-- (4.513318348048452,1.383678342379166);
\draw [line width=0.3pt] (4.513318348048452,1.383678342379166)-- (5.2088901065771696,0.9797979664592665);
\draw [line width=0.3pt] (5.2088901065771696,0.9797979664592665)-- (5,0);
\draw [line width=0.3pt] (3.795308790857518,0.9910168657903748)-- (4,0);
\draw [line width=0.3pt] (4,0)-- (5,0);
\draw [line width=0.3pt] (5.2088901065771696,0.9797979664592665)-- (5.399611395206011,1.9782800069279074);
\begin{scriptsize}
\draw [fill=black] (-0.00970093979975073,0) circle (1pt);
\draw [fill=black] (0.9902990602002507,0) circle (1pt);
\draw [fill=black] (-0.19861937101705318,0.9910168657903746) circle (1pt);
\draw [fill=black] (1.181305246709273,0.9910168657903746) circle (1pt);
\draw [fill=black] (0.48573348818055523,1.4061161410413825) circle (1pt);
\draw [fill=black] (1.394464334000332,1.9894989062590158) circle (1pt);
\draw [fill=black] (2.3817274751378643,2.0007178055901234) circle (1pt);
\draw [fill=black] (2.5836676630978146,1.002235765121483) circle (1pt);
\draw [fill=black] (1.9105337032313154,0.5759175905393669) circle (1pt);
\draw [fill=black] (-0.19861937101705318,0.9910168657903746) circle (1pt);
\draw [fill=black] (1.181305246709273,0.9910168657903746) circle (1pt);
\draw [fill=black] (0.48573348818055523,1.4061161410413825) circle (1pt);
\draw [fill=black] (1.394464334000332,1.9894989062590158) circle (1pt);
\draw [fill=black] (2.3817274751378643,2.0007178055901234) circle (1pt);
\draw [fill=black] (2.5836676630978146,1.002235765121483) circle (1pt);
\draw [fill=black] (1.9105337032313154,0.5759175905393669) circle (1pt);
\draw [fill=black] (-0.19861937101705318,0.9910168657903746) circle (1pt);
\draw [fill=black] (1.181305246709273,0.9910168657903746) circle (1pt);
\draw [fill=black] (0.48573348818055523,1.4061161410413825) circle (1pt);
\draw [fill=black] (1.394464334000332,1.9894989062590158) circle (1pt);
\draw [fill=black] (2.3817274751378643,2.0007178055901234) circle (1pt);
\draw [fill=black] (2.5836676630978146,1.002235765121483) circle (1pt);
\draw [fill=black] (1.9105337032313154,0.5759175905393669) circle (1pt);
\draw [fill=black] (-0.19861937101705318,0.9910168657903746) circle (1pt);
\draw [fill=black] (1.181305246709273,0.9910168657903746) circle (1pt);
\draw [fill=black] (0.48573348818055523,1.4061161410413825) circle (1pt);
\draw [fill=black] (1.394464334000332,1.9894989062590158) circle (1pt);
\draw [fill=black] (2.3817274751378643,2.0007178055901234) circle (1pt);
\draw [fill=black] (2.5836676630978146,1.002235765121483) circle (1pt);
\draw [fill=black] (1.9105337032313154,0.5759175905393669) circle (1pt);
\draw [fill=black] (3.795308790857518,0.9910168657903748) circle (1pt);
\draw [fill=black] (5.2088901065771696,0.9797979664592665) circle (1pt);
\draw [fill=black] (4,0) circle (1pt);
\draw [fill=black] (5,0) circle (1pt);
\draw [fill=black] (4.513318348048452,1.383678342379166) circle (1pt);
\draw [fill=black] (5.399611395206011,1.9782800069279074) circle (1pt);

\draw[color=black] (3.2,1) node {$...$};

%%%%Pentn
\draw[color=black] (4.5,1.7) node {\tiny$x^{n}_{2}$};

\draw[color=black] (4.1,0.9) node {\tiny$x^{n}_{1}$};

\draw[color=black] (5.4,0) node {\tiny$x^{n}_{4}$};

\draw[color=black] (3.7,0) node {\tiny$x^{n}_{5}$};

\draw[color=black] (5.6,1.2) node {\tiny$x^{n}_{3}$};

\draw[color=black] (5.8,2.4) node {\tiny$x^{n+1}_{5}$};

\end{scriptsize}
\end{tikzpicture}

 }%
\vskip -3 cm
\caption{The graph $\overline{M}(n)$.}
\label{mpbar}
\end{figure}
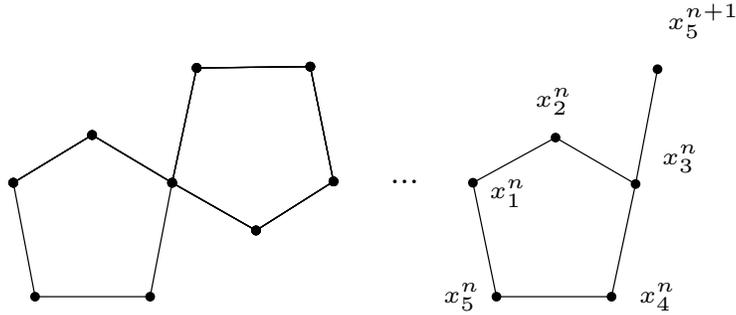

\begin{figure}[H]
\centering
\definecolor{ududff}{rgb}{0.30196078431372547,0.30196078431372547,1}
\resizebox{1.2\textwidth}{!}{%

\begin{tikzpicture}[line cap=round,line join=round,>=triangle 45,x=1cm,y=1cm]
\clip(-1.600981787405596,-2.7897522087931304) rectangle (10.42567829554255,3.986462987196297);
\draw [line width=0.3pt] (-0.19861937101705318,0.9910168657903746)-- (0.48573348818055523,1.4061161410413825);
\draw [line width=0.3pt] (0.48573348818055523,1.4061161410413825)-- (1.181305246709273,0.9910168657903746);
\draw [line width=0.3pt] (1.181305246709273,0.9910168657903746)-- (1.9105337032313154,0.5759175905393669);
\draw [line width=0.3pt] (1.9105337032313154,0.5759175905393669)-- (2.5836676630978146,1.002235765121483);
\draw [line width=0.3pt] (2.5836676630978146,1.002235765121483)-- (2.3817274751378643,2.0007178055901234);
\draw [line width=0.3pt] (2.3817274751378643,2.0007178055901234)-- (1.394464334000332,1.9894989062590158);
\draw [line width=0.3pt] (1.394464334000332,1.9894989062590158)-- (1.181305246709273,0.9910168657903746);
\draw [line width=0.3pt] (-0.19861937101705318,0.9910168657903746)-- (-0.00970093979975073,0);
\draw [line width=0.3pt] (-0.00970093979975073,0)-- (0.9902990602002507,0);
\draw [line width=0.3pt] (0.9902990602002507,0)-- (1.181305246709273,0.9910168657903746);
\draw [line width=0.3pt] (-0.19861937101705318,0.9910168657903746)-- (0.48573348818055523,1.4061161410413825);
\draw [line width=0.3pt] (0.48573348818055523,1.4061161410413825)-- (1.181305246709273,0.9910168657903746);
\draw [line width=0.3pt] (1.181305246709273,0.9910168657903746)-- (1.9105337032313154,0.5759175905393669);
\draw [line width=0.3pt] (1.9105337032313154,0.5759175905393669)-- (2.5836676630978146,1.002235765121483);
\draw [line width=0.3pt] (2.5836676630978146,1.002235765121483)-- (2.3817274751378643,2.0007178055901234);
\draw [line width=0.3pt] (2.3817274751378643,2.0007178055901234)-- (1.394464334000332,1.9894989062590158);
\draw [line width=0.3pt] (1.394464334000332,1.9894989062590158)-- (1.181305246709273,0.9910168657903746);
\draw [line width=0.3pt] (-0.19861937101705318,0.9910168657903746)-- (0.48573348818055523,1.4061161410413825);
\draw [line width=0.3pt] (0.48573348818055523,1.4061161410413825)-- (1.181305246709273,0.9910168657903746);
\draw [line width=0.3pt] (1.181305246709273,0.9910168657903746)-- (1.9105337032313154,0.5759175905393669);
\draw [line width=0.3pt] (1.9105337032313154,0.5759175905393669)-- (2.5836676630978146,1.002235765121483);
\draw [line width=0.3pt] (2.5836676630978146,1.002235765121483)-- (2.3817274751378643,2.0007178055901234);
\draw [line width=0.3pt] (2.3817274751378643,2.0007178055901234)-- (1.394464334000332,1.9894989062590158);
\draw [line width=0.3pt] (1.394464334000332,1.9894989062590158)-- (1.181305246709273,0.9910168657903746);
\draw [line width=0.3pt] (-0.19861937101705318,0.9910168657903746)-- (0.48573348818055523,1.4061161410413825);
\draw [line width=0.3pt] (0.48573348818055523,1.4061161410413825)-- (1.181305246709273,0.9910168657903746);
\draw [line width=0.3pt] (1.181305246709273,0.9910168657903746)-- (1.9105337032313154,0.5759175905393669);
\draw [line width=0.3pt] (1.9105337032313154,0.5759175905393669)-- (2.5836676630978146,1.002235765121483);
\draw [line width=0.3pt] (2.5836676630978146,1.002235765121483)-- (2.3817274751378643,2.0007178055901234);
\draw [line width=0.3pt] (2.3817274751378643,2.0007178055901234)-- (1.394464334000332,1.9894989062590158);
\draw [line width=0.3pt] (1.394464334000332,1.9894989062590158)-- (1.181305246709273,0.9910168657903746);
\draw [line width=0.3pt] (3.795308790857518,0.9910168657903748)-- (4.513318348048452,1.383678342379166);
\draw [line width=0.3pt] (4.513318348048452,1.383678342379166)-- (5.2088901065771696,0.9797979664592665);
\draw [line width=0.3pt] (5.2088901065771696,0.9797979664592665)-- (5,0);
\draw [line width=0.3pt] (3.795308790857518,0.9910168657903748)-- (4,0);
\draw [line width=0.3pt] (4,0)-- (5,0);
\draw [line width=0.3pt] (5.2088901065771696,0.9797979664592665)-- (5.399611395206011,1.9782800069279074);
\draw [line width=0.3pt] (5.399611395206011,1.9782800069279074)-- (6.409312335005763,1.9894989062590156);
\draw [line width=0.3pt] (5.2088901065771696,0.9797979664592665)-- (5.926899663768104,0.5983553892015835);
\begin{scriptsize}
\draw [fill=black] (-0.00970093979975073,0) circle (1pt);
\draw [fill=black] (0.9902990602002507,0) circle (1pt);
\draw [fill=black] (-0.19861937101705318,0.9910168657903746) circle (1pt);
\draw [fill=black] (1.181305246709273,0.9910168657903746) circle (1pt);
\draw [fill=black] (0.48573348818055523,1.4061161410413825) circle (1pt);
\draw [fill=black] (1.394464334000332,1.9894989062590158) circle (1pt);
\draw [fill=black] (2.3817274751378643,2.0007178055901234) circle (1pt);
\draw [fill=black] (2.5836676630978146,1.002235765121483) circle (1pt);
\draw [fill=black] (1.9105337032313154,0.5759175905393669) circle (1pt);
\draw [fill=black] (-0.19861937101705318,0.9910168657903746) circle (1pt);
\draw [fill=black] (1.181305246709273,0.9910168657903746) circle (1pt);
\draw [fill=black] (0.48573348818055523,1.4061161410413825) circle (1pt);
\draw [fill=black] (1.394464334000332,1.9894989062590158) circle (1pt);
\draw [fill=black] (2.3817274751378643,2.0007178055901234) circle (1pt);
\draw [fill=black] (2.5836676630978146,1.002235765121483) circle (1pt);
\draw [fill=black] (1.9105337032313154,0.5759175905393669) circle (1pt);
\draw [fill=black] (-0.19861937101705318,0.9910168657903746) circle (1pt);
\draw [fill=black] (1.181305246709273,0.9910168657903746) circle (1pt);
\draw [fill=black] (0.48573348818055523,1.4061161410413825) circle (1pt);
\draw [fill=black] (1.394464334000332,1.9894989062590158) circle (1pt);
\draw [fill=black] (2.3817274751378643,2.0007178055901234) circle (1pt);
\draw [fill=black] (2.5836676630978146,1.002235765121483) circle (1pt);
\draw [fill=black] (1.9105337032313154,0.5759175905393669) circle (1pt);
\draw [fill=black] (-0.19861937101705318,0.9910168657903746) circle (1pt);
\draw [fill=black] (1.181305246709273,0.9910168657903746) circle (1pt);
\draw [fill=black] (0.48573348818055523,1.4061161410413825) circle (1pt);
\draw [fill=black] (1.394464334000332,1.9894989062590158) circle (1pt);
\draw [fill=black] (2.3817274751378643,2.0007178055901234) circle (1pt);
\draw [fill=black] (2.5836676630978146,1.002235765121483) circle (1pt);
\draw [fill=black] (1.9105337032313154,0.5759175905393669) circle (1pt);
\draw [fill=black] (3.795308790857518,0.9910168657903748) circle (1pt);
\draw [fill=black] (5.2088901065771696,0.9797979664592665) circle (1pt);
\draw [fill=black] (4,0) circle (1pt);
\draw [fill=black] (5,0) circle (1pt);
\draw [fill=black] (4.513318348048452,1.383678342379166) circle (1pt);
\draw [fill=black] (5.399611395206011,1.9782800069279074) circle (1pt);
\draw [fill=black] (6.409312335005763,1.9894989062590156) circle (1pt);
\draw [fill=black] (5.926899663768104,0.5983553892015835) circle (1pt);

\draw[color=black] (3.2,1) node {$...$};

%%%%Pentn
\draw[color=black] (4.5,1.7) node {\tiny$x^{n}_{2}$};

\draw[color=black] (4.2,0.8) node {\tiny$x^{n}_{1}$};

\draw[color=black] (5.4,0) node {\tiny$x^{n}_{4}$};

\draw[color=black] (3.7,0) node {\tiny$x^{n}_{5}$};

\draw[color=black] (6,1.2) node {\tiny$x^{n}_{3}= x^{n+1}_{1}$};

%%%pentn+1
\draw[color=black] (5.5,2.4) node {\tiny$x^{n+1}_{5}$};

\draw[color=black] (7,2.4) node {\tiny$x^{n+1}_{4}$};

\draw[color=black] (6.5,0.55) node {\tiny$x^{n+1}_{2}$};

\end{scriptsize}
\end{tikzpicture}

 }%
\vskip -3 cm
\caption{The graph $\tilde{M}(n)$.}
\label{mptilde}
\end{figure}
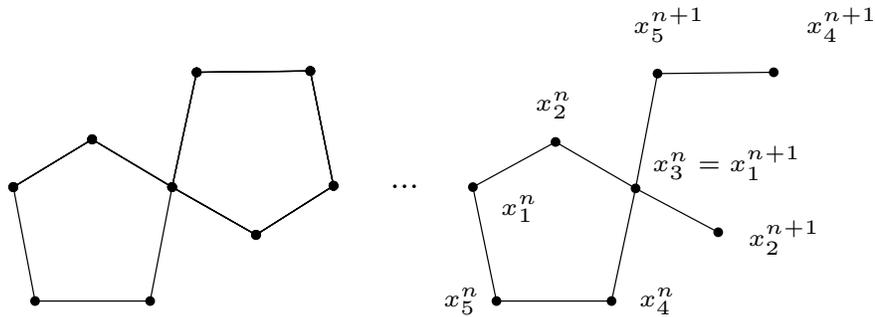

\noindent Next, we let $\overline{M}(n)$ be constructed from $M(n)$ by joining one vertex to a vertex at distance two from the cut vertex of the $n^{th}$ pentagon. The graph $\overline{M}(n)$ is illustrated by Figure \ref{mpbar}.\indent Further, we define the graph $\tilde{M}$ from $M(n)$ by identifying a vertex of degree two of a path of length three to a vertex at distance two from the cut vertex of the $n^{th}$ pentagon. The graph $\tilde{M}(n)$ is shown in Figure \ref{mptilde}.

\vskip -1 cm
\vskip 5 pt

\vskip -2 cm
\vskip 5 pt

\noindent Then, we let
\vskip 5 pt

\indent $\overline{m}(n, k)$ = the number of maximal independent sets containing $k$ vertices of $\overline{M}(n)$,
\vskip 5 pt

\indent $\tilde{m}(n, k)$ = the number of maximal independent sets containing $k$ vertices of $\tilde{M}(n)$,
\vskip 5 pt

\noindent and let
\begin{align}
\overline{M}(x, y) = \sum_{n \geq 0}\sum_{k \geq 0}\overline{m}(n, k)x^{n}y^{k}\notag
\end{align}
\begin{align}
\tilde{M}(x, y) = \sum_{n \geq 0}\sum_{k \geq 0}\tilde{m}(n, k)x^{n}y^{k}.\notag
\end{align}
\noindent Now, we are ready to prove Theorem \ref{bivariate-m}
\vskip 5 pt

\noindent \emph{Proof of Theorem \ref{bivariate-m}} Let $D$ be a maximal independent set of $M(n)$ containing $k$ vertices. We distinguish $2$ cases. 
\vskip 5 pt

\noindent \textbf{Case 1:} $x^{n}_{4} \in D$\\
\indent Thus, $x^{n}_{3}, x^{n}_{5} \notin D$. Removing $x^{n}_{3}, x^{n}_{4} , x^{n}_{5}$ from $M(n)$ result in $\overline{M}(n - 1)$. There are $\overline{m}(n - 1, k - 1)$ possibilities of $D$.
\vskip 5 pt

\noindent \textbf{Case 2:} $x^{n}_{4} \notin D$\\
\indent We further have the following $3$ subcases.
\vskip 5 pt

\noindent \textbf{Subcase 2.1:} $n^{n}_{3} \notin D$ and $x^{n}_{5} \in D$.\\
\indent Thus, $x^{n}_{1} \notin D$ and $x^{n}_{2} \in D$. Removing all the vertices of the $n^{th}$ pentagon results in $\tilde{M}(n - 2)$. We have that there are $\tilde{m}(n - 2, k - 2)$ possibilities of $D$.
\vskip 5 pt

\noindent \textbf{Subcase 2.2:} $n^{n}_{3} \in D$ and $x^{n}_{5} \notin D$.\\
\indent Thus, $x^{n}_{2} \notin D$. By the maximality of $D$, $x^{n}_{1} \in D$ and this implies that $x^{n - 1}_{2}, x^{n - 1}_{4} \notin D$. Removing all the vertices of the $n^{th}$ pentagon and $x^{n - 1}_{2}, x^{n - 1}_{4}$ results in $\overline{M}(n - 2)$. We have that there are $\overline{m}(n - 2, k - 2)$ possibilities of $D$.
\vskip 5 pt

\noindent \textbf{Subcase 2.3:} $n^{n}_{3} \in D$ and $x^{n}_{5} \in D$.\\
\indent Thus, $x^{n}_{2} \in D$. Removing all the vertices of the $n^{th}$ pentagon results in $\tilde{M}(n - 2)$. We have that there are $\tilde{m}(n - 2, k - 2)$ possibilities of $D$.
\vskip 5 pt

\indent From, all the cases, we have that
\begin{align}\label{m11}
m(n, k) = \overline{m}(n - 1, k - 1) + 2\Tilde{m}(n - 2, k - 2) +  \overline{m}(n - 2, k - 2)
\end{align}
\noindent For $n \geq 2 $ and $k \geq 2 $, we multiply $x^{n}y^{k}$ throughout (\ref{m11}) and sum over all $x^{n}y^{k}$. Thus, we have that

\begin{align}\label{m1111}
\sum_{n \geq 2}\sum_{k \geq 2}m(n, k)x^{n}y^{k} &= \sum_{n \geq 2}\sum_{k \geq 2}\overline{m}(n - 1, k - 1)x^{n}y^{k} + 2\sum_{n \geq 2}\sum_{k \geq 2}\Tilde{m}(n - 2, k - 2)x^{n}y^{k} \notag\\
& \text{ }  + \sum_{n \geq 2}\sum_{k \geq 2}\overline{m}(n - 2, k - 2)x^{n}y^{k}.
\end{align}

\noindent We first consider the term $\displaystyle{\sum_{n \geq 2}\sum_{k \geq 2}m(n, k)x^{n}y^{k}}$. It can be checked that

\begin{center}
    $m(0, 0) = 1, m(0, 1) = 0, m(0, 2) = 0$, and $ m(0, k) = 0$ for all $k \geq 3$
    \vskip 5 pt
    $m(1, 0) = 0, m(1, 0) = 0, m(1, 2) = 5$, and $ m(0, k) = 0$ for all $k \geq 3$
    \vskip 5 pt
\end{center}

\begin{align}\label{m1}
\sum_{n \geq 2}\sum_{k \geq 2}m(n, k)x^{n}y^{k} &= \sum_{n \geq 2}\sum_{k \geq 2}m(n, k)x^{n}y^{k}  + m(0, 0)  + m(1, 2)xy^{2}  -m(0, 0) -m(1, 2)xy^{2}\notag\\
                  &= \sum_{n \geq 0}\sum_{k \geq 0}m(n, k)x^{n}y^{k} - 1 - 5xy^2 \notag\\
                  &= M(x,y) - 1 - 5xy^2
\end{align}

\noindent Now, we consider the term $\displaystyle{\sum_{n \geq 2}\sum_{k \geq 2}\overline{m}(n - 1, k - 1)x^{n}y^{k}}$. clearly,

\begin{center}
    $\overline{m}(0, 0) = 0, \overline{m}(0, 1) = 2$ and $\overline{m}(0, k) = 0$ for all $k \geq 2$
    \vskip 5 pt
    $\overline{m}(1, 0) = 0, \overline{m}(1, 1) = 0$ and $\overline{m}(1, k) = 0$ for all $k \geq 2$
    \vskip 5 pt
\end{center}

\begin{align}\label{mbar1}
\sum_{n \geq 2}\sum_{k \geq 2}\overline{m}(n - 1, k - 1)x^{n}y^{k}
    &= xy(\sum_{n \geq 2}\sum_{k \geq 2}\overline{m}(n - 1, k - 1)x^{n - 1}y^{k - 1})\notag\\
    &= xy(\sum_{n \geq 2}\sum_{k \geq 2}\overline{m}(n - 1, k - 1)x^{n - 1}y^{k - 1} + \overline{m}(0,1)y - \overline{m}(0,1)y)\notag\\
    &= xy(\sum_{n \geq 0}\sum_{k \geq 0}\overline{m}(n, k)x^{n}y^{k} - 2y)\notag\\
    &= xy\overline{M}(x, y) - 2xy^{2}.
\end{align}

\noindent Next, we consider the term $\displaystyle{\sum_{n \geq 2}\sum_{k \geq 2}\Tilde{2m}(n - 2, k - 2)x^{n}y^{k}}$. It can be check that,

\begin{align}\label{mtil1}
\sum_{n \geq 2}\sum_{k \geq 2}\Tilde{2m}(n - 2, k - 2)x^{n}y^{k}
    &= 2x^{2}y^{2}(\sum_{n \geq 2}\sum_{k \geq 2}\Tilde{m}(n - 2, k - 2)x^{n - 2}y^{k - 2})\notag\\
    &= 2x^{2}y^{2}(\sum_{n \geq 0}\sum_{k \geq 0}\Tilde{m}(n , k)x^{n}y^{k}\notag\\
    &= 2x^{2}y^{2}\Tilde{M}(x, y).
\end{align}

\noindent Finally, we consider the term $\displaystyle{\sum_{n \geq 2}\sum_{k \geq 2}\overline{m}(n - 2, k - 2)x^{n}y^{k}}$. It can be check that,

\begin{align}\label{mbar2}
\sum_{n \geq 2}\sum_{k \geq 2}\overline{m}(n - 2, k - 2)x^{n}y^{k}
    &= x^{2}y^{2}(\sum_{n \geq 2}\sum_{k \geq 2}\overline{m}(n - 2, k - 2)x^{n - 2}y^{k - 2})\notag\\
    &= x^{2}y^{2}(\sum_{n \geq 0}\sum_{k \geq 0}\overline{m}(n, k)x^{n}y^{k})\notag\\
    &= x^{2}y^{2}\overline{M}(x, y).
\end{align}

\noindent Plugging (\ref{m1}), (\ref{mbar1}) ,(\ref{mtil1}) and (\ref{mbar2}) to (\ref{m1111}), we have 
\begin{align}
M(x,y) - 1- 5xy^2 &= xy\overline{M}(x, y) - 2xy^{2} + 2x^{2}y^{2}\Tilde{M}(x, y) + x^{2}y^{2}\overline{M}(x, y)\notag
\end{align}

\noindent which can be solved that

\begin{align}\label{semifinalm4}
M(x,y) &= 1 + 5xy^2 + xy\overline{M}(x, y) - 2xy^{2} + 2x^{2}y^{2}\Tilde{M}(x, y) + x^{2}y^{2}\overline{M}(x, y).
\end{align}
\vskip 5 pt

\indent Next, we will find the recurrence relation of $\overline{m}(n, k)$. Let $D$ be a maximal independent set of $\overline{M}(n)$ containing $k$ vertices. We distinguish $2$ cases. 
\vskip 5 pt

\noindent \textbf{Case 1:} $x^{n + 1}_{2} \in D$\\
\indent Thus, $x^{n}_{3} \notin D$. Removing $x^{n}_{3}, x^{n + 1}_{2}$ from $\overline{M}(n)$ results in $\tilde{M}(n - 1)$. There are $\tilde{m}(n - 1, k - 1)$ possibilities of $D$.
\vskip 5 pt

\noindent \textbf{Case 2:} $x^{n + 1}_{2} \notin D$\\
\indent By the maximality of $D$, $x^{n}_{3} \in D$. Thus, $x^{n}_{2}, x^{n}_{4} \notin D$. Removing the vertices $x^{n}_{2}, x^{n}_{3}, x^{n}_{4}, x^{n + 1}_{2}$ from $\overline{M}(n)$ results in $\overline{M}(n - 1)$. There are $\overline{m}(n - 1, k - 1)$ possibilities of $D$.
\vskip 5 pt

\indent From the two cases, we have that
\begin{align}\label{mbar3}
\overline{m}(n, k) = \Tilde{m}(n - 1, k - 1) + \overline{m}(n - 1, k - 1)
\end{align}
\vskip 5 pt

\noindent For $n \geq 1$ and $k \geq 1$, we multiply $x^{n}y^{k}$ throughout (\ref{mbar3}) and sum over all $x^{n}y^{k}$. Thus, we have that

\begin{align}\label{mbar3333}
\sum_{n \geq 1}\sum_{k \geq 1}\overline{m}(n, k)x^{n}y^{k} = \sum_{n \geq 1}\sum_{k \geq 1}\Tilde{m}(n - 1, k - 1)x^{n}y^{k} + \sum_{n \geq 1}\sum_{k \geq 1}\overline{m}(n - 1, k - 1)x^{n}y^{k}
\end{align}

\noindent We first consider the term $\displaystyle{\sum_{n \geq 1}\sum_{k \geq 1}\overline{m}(n , k )x^{n}y^{k}}$. It can be checked that, 

\begin{center}
    $\overline{m}(0, 0) = 0, \overline{m}(0, 1) = 2$ and $\overline{m}(0, k) = 0$ for all $k \geq 2$
    \vskip 5 pt
    $\overline{m}(1, 0) = 0, \overline{m}(1, 1) = 0$ and $\overline{m}(1, k) = 0$ for all $k \geq 2$
    \vskip 5 pt
\end{center}

\begin{align}\label{mbar4}
\sum_{n \geq 1}\sum_{k \geq 1}\overline{m}(n, k)x^{n}y^{k}
    &= \sum_{n \geq 1}\sum_{k \geq 1}\overline{m}(n - 1 , k - 1)x^{n }y^{k} + \overline{m}(0, 1)y - \overline{m}(0, 1)y \notag\\
    &= \sum_{n \geq 0}\sum_{k \geq 0}\overline{m}(n , k )x^{n}y^{k} - 2y\notag\\
    &= \overline{M}(x, y) - 2y
\end{align}

\noindent We next consider the term $\displaystyle{\sum_{n \geq 1}\sum_{k \geq 1}\Tilde{m}(n - 1 , k - 1 )x^{n}y^{k}}$. It can be checked that,

\begin{align}\label{mtil2}
\sum_{n \geq 1}\sum_{k \geq 1}\Tilde{m}(n - 1 , k - 1)x^{n}y^{k}
    &= xy(\sum_{n \geq 1}\sum_{k \geq 1}\Tilde{m}(n - 1 , k - 1 )x^{n - 1}y^{k - 1})\notag\\
    &= xy(\sum_{n \geq 0}\sum_{k \geq 0}\Tilde{m}(n , k )x^{n}y^{k} )\notag\\
    &= xy\Tilde{M}(x, y)
\end{align}    
\vskip 5 pt

\noindent Finally, We next consider the term $\displaystyle{\sum_{n \geq 1}\sum_{k \geq 1}\overline{m}(n - 1 , k - 1)x^{n}y^{k}}$. Cleary 

\begin{align}\label{mbar5}
\sum_{n \geq 1}\sum_{k \geq 1}\overline{m}(n - 1 , k - 1)x^{n}y^{k}
    &= xy(\sum_{n \geq 1}\sum_{k \geq 1}\overline{m}(n - 1 , k -1 )x^{n - 1}y^{k - 1})\notag\\
    &= xy(\sum_{n \geq 0}\sum_{k \geq 0}\overline{m}(n , k )x^{n}y^{k})\notag\\
    &= xy\overline{M}(x, y)
\end{align}    
\vskip 5 pt

\noindent Plugging (\ref{mbar4}) , (\ref{mtil2}) and (\ref{mbar5}) to (\ref{mbar3333}) we have 

\begin{align}
\overline{M}(x, y) - 2y &= xy\Tilde{M}(x, y) + xy\overline{M}(x, y).\notag
\end{align}

\noindent which can be solved that

\begin{align}\label{semifinalmbar6}
\overline{M}(x, y) &= 2y + xy\Tilde{M}(x, y) + xy\overline{M}(x, y).
\end{align}
\vskip 5 pt

\indent Finally, we will find the recurrence relation of $\tilde{m}(n, k)$. Let $D$ be a maximal independent set of $\tilde{M}(n)$ containing $k$ vertices. We distinguish $2$ cases. 
\vskip 5 pt

\noindent \textbf{Case 1:} $x^{n+1}_{2} \in D$\\
\indent Thus, $x^{n}_{2}, x^{n}_{4}, x^{n+1}_{1}, x^{n+1}_{3} \notin D$. By the maximality of $D$, $x^{n + 1}_{4} \in D$. Removing $x^{n}_{2}, x^{n}_{4}, x^{n+1}_{1},$ $... , x^{n+1}_{4}$ from $\tilde{M}(n)$ results in $\overline{M}(n - 1)$. There are $\overline{m}(n - 1, k - 2)$ possibilities of $D$.
\vskip 5 pt

\noindent \textbf{Case 2:} $x^{n + 1}_{2} \notin D$\\
\indent By the maximality of $D$, $x^{n + 1}_{1} \in D$. Further, $x^{n+1}_{3} \in D$ or $x^{n+1}_{4} \in D$. Removing all the vertices $x^{n+1}_{1}, ... , x^{n+1}_{4}$ results in $\tilde{M}(n - 1)$. We have that there are $2\tilde{m}(n - 1, k - 2)$ possibilities of $D$.
\vskip 5 pt

\indent From all the cases, we have that
\begin{align}\label{mtil3}
    \Tilde{m}(n, k) = \overline{m}(n - 1, k - 2) + \Tilde{2m}(n - 1, k - 2)
\end{align}
\vskip 5 pt

\noindent For $n \geq 1$ and $k \geq 2$, we multiply $x^{n}y^{k}$ throughout (\ref{mtil3}) and sum over all $x^{n}y^{k}$. Thus, we have that

\begin{align}\label{mtil3333}
\sum_{n \geq 1}\sum_{k \geq 2}\Tilde{m}(n, k)x^{n}y^{k} = \sum_{n \geq 1}\sum_{k \geq 2}\overline{m}(n - 1, k - 2)x^{n}y^{k} + \sum_{n \geq 1}\sum_{k \geq 2}\Tilde{2m}(n - 1, k - 2)x^{n}y^{k}
\end{align}

\noindent We first consider the term $\displaystyle{\sum_{n \geq 1}\sum_{k \geq 2}\Tilde{m}(n , k )x^{n}y^{k}}$. It can be checked that, 

\begin{center}
    $\Tilde{m}(0, 0) = 0, \Tilde{m}(0, 1) = 0 , \Tilde{m}(0, 2) = 3$ and $\Tilde{m}(0, k) = 0$ for all $k \geq 3$
    \vskip 5 pt
    $\Tilde{m}(1, 0) = 0, \Tilde{m}(1, 1) = 0, \Tilde{m}(0,2) = 0 $ and $\Tilde{m}(0, k) = 0$ for all $k \geq 3$
    \vskip 5 pt
\end{center}

\begin{align}\label{mtil4}
\sum_{n \geq 1}\sum_{k \geq 2}\Tilde{m}(n , k)x^{n}y^{k}
    &= \sum_{n \geq 1}\sum_{k \geq 2}\Tilde{m}(n , k)x^{n}y^{k}
    + \Tilde{m}(0,2)y^{2} -\Tilde{m}(0,2)y^{2}  \notag\\
    &= \sum_{n \geq 0}\sum_{k \geq 0}\Tilde{m}(n , k )x^{n}y^{k} - 3y^{2}\notag\\
    &= \Tilde{M}(x, y) - 3y^{2}
\end{align}    
\vskip 5 pt

\noindent Now, we consider the term $\displaystyle{\sum_{n \geq 1}\sum_{k \geq 2}\overline{m}(n - 1 , k - 2)x^{n}y^{k}}$. Cleary,

\begin{align}\label{mbar7}
\sum_{n \geq 1}\sum_{k \geq 2}\overline{m}(n - 1 , k - 2)x^{n}y^{k}
    &= xy^{2}(\sum_{n \geq 1}\sum_{k \geq 2}\overline{m}(n - 1 , k - 2 )x^{n - 1}y^{k - 2})\notag\\
    &= xy^{2}(\sum_{n \geq 0}\sum_{k \geq 0}\overline{m}(n, k)x^{n} y^{k} )\notag\\
    &= xy^{2}\overline{M}(x, y) 
\end{align}    
\vskip 5 pt

\noindent Finally, we consider the term
$\displaystyle{\sum_{n \geq 1}\sum_{k \geq 2}2\Tilde{m}(n - 1 , k - 2)x^{n}y^{k}}$. It can be checked that, 

\begin{align}\label{mtil5}
\sum_{n \geq 1}\sum_{k \geq 2}2\Tilde{m}(n - 1 , k - 2)x^{n}y^{k}
    &= 2xy^{2}(\sum_{n \geq 1}\sum_{k \geq 2}\Tilde{m}(n - 1 , k - 2 )x^{n - 1}y^{k - 2})\notag\\
    &= 2xy^{2}(\sum_{n \geq 0}\sum_{k \geq 0}\Tilde{m}(n , k )x^{n}y^{k} )\notag\\
    &= 2xy^{2}\Tilde{M}(x, y)
\end{align}

\noindent Plugging (\ref{mtil4}), (\ref{mbar7}) , and (\ref{mtil5}) to (\ref{mtil3333}), we have 

\begin{align}
\Tilde{M}(x , y) - 3y^{2} = xy^{2}\overline{M}(x , y) + 2xy^{2}\Tilde{M}(x , y)\notag
\end{align}

\noindent which can be solved that

\begin{align}\label{semifinalmbar7}
\Tilde{M}(x, y) = 3y^{2} + xy^{2}\overline{M}(x , y) + 2xy^{2}\Tilde{M}(x , y)
\end{align}
\vskip 5 pt

\noindent By plugging (\ref{semifinalm4}) to (\ref{semifinalmbar6}) to (\ref{semifinalmbar7}) , we have
\begin{align}
M(x, y) = \frac{-11x^{2}y^{4} + 5xy^{2} + 2x^{3}y^{6} + 2x^{3}y^{5} - x^{2}y^{3}}{1 - 4xy^{2} - xy + 3x^{2}y^{3} + 4x^{2}y^{4} - 2x^{3}y^{5}} + 1\notag
\end{align}
\noindent as required.
\qed

\indent Now, we are ready to prove Theorem \ref{recurrence-m}
\vskip 5 pt

\noindent \emph{Proof of Theorem \ref{recurrence-m}} By Theorem \ref{bivariate-m} with $y = 1$, we have that
\vskip 5 pt

\begin{align}
\sum_{n \geq 0}m(n)x^{n} = M(x, 1) = \frac{1 - 5x^{2} + 2x^{3}}{1 - 5x + 7x^{2} - 2x^{3}} = \frac{-x^{2} + 2x + 1}{x^{2} - 3x + 1}
\end{align}
 
 \noindent which can  be solved that

\begin{align}\label{recursived-m2}
-x^{2} + 2x + 1   &= (x^{2} - 3x + 1)\sum_{n \geq 0}m(n)x^{n}\notag\\
                  &= \sum_{n \geq 0}m(n)x^{n + 2}  - 3\sum_{n \geq 0}m(n)x^{n + 1} + \sum_{n \geq 0}n(n)x^{n} \notag\\
                  &= \sum_{n \geq 2}m(n-2)x^{n} - 3\sum_{n \geq 1}m(n - 1)x^{n} + \sum_{n \geq 0}m(n)x^{n} \notag\\
                  &= m(0) - 3m(0)x + m(1)x + \sum_{n \geq 2}(m(n-2) - 3m(n-1) + m(n))x^{n}
\end{align}
\noindent Because the order of the polynomial on the left hand side of (\ref{recursived-m2}) is two, the coefficients of $x^{n}$ for all $n \geq 3$ must be $0$. Thus, $m(n) - 3m(n-1) + m(n - 2) = 0$ implying that

\begin{align*}
    m(n) = 3m(n - 1) - m(n - 2). 
\end{align*}

\subsection{Meta-Hexagonal Cacti}
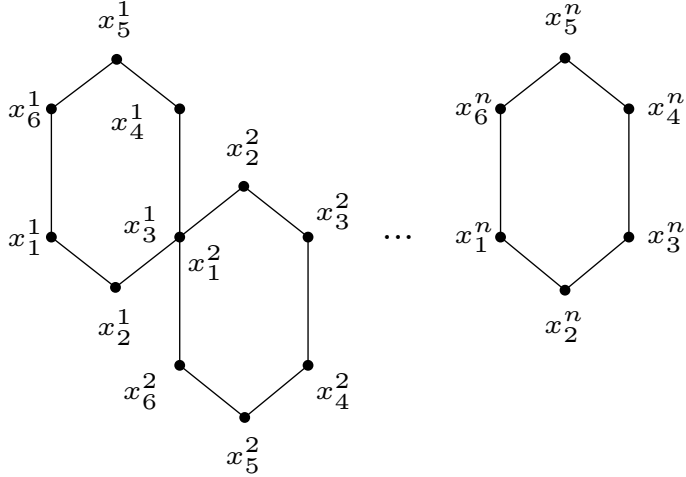
\begin{figure}[H]
\centering
\definecolor{ududff}{rgb}{0.30196078431372547,0.30196078431372547,1}
\resizebox{0.8\textwidth}{!}{%

\begin{tikzpicture}[line cap=round,line join=round,>=triangle 45,x=1cm,y=1cm]
\clip(-0.42006259682951436,-2.0634421903550786) rectangle (6.717157112872825,2.344411290892576);
\draw [line width=0.3pt] (0,1)-- (0,0);
\draw [line width=0.3pt] (0,1)-- (0.5090377499796499,1.3861738649290605);
\draw [line width=0.3pt] (0.5090377499796499,1.3861738649290605)-- (1,1);
\draw [line width=0.3pt] (1,1)-- (1,0);
\draw [line width=0.3pt] (0,0)-- (0.49945650687446763,-0.39225270822641545);
\draw [line width=0.3pt] (0.49945650687446763,-0.39225270822641545)-- (1,0);
\draw [line width=0.3pt] (1,0)-- (1.49925108788594,0.39590652351985406);
\draw [line width=0.3pt] (1.49925108788594,0.39590652351985406)-- (2,0);
\draw [line width=0.3pt] (1,0)-- (1,-1);
\draw [line width=0.3pt] (1,-1)-- (1.5065488585502573,-1.4066428305665222);
\draw [line width=0.3pt] (1.5065488585502573,-1.4066428305665222)-- (2,-1);
\draw [line width=0.3pt] (2,0)-- (2,-1);
\draw [line width=0.3pt] (3.5,1)-- (3.5,0);
\draw [line width=0.3pt] (3.5,1)-- (4.002386425746781,1.3957011045313252);
\draw [line width=0.3pt] (4.002386425746781,1.3957011045313252)-- (4.5,1);
\draw [line width=0.3pt] (4.5,1)-- (4.5,0);
\draw [line width=0.3pt] (4.5,0)-- (4.002386425746781,-0.4141460202193678);
\draw [line width=0.3pt] (3.5,0)-- (4.002386425746781,-0.4141460202193678);
\begin{scriptsize}
\draw [fill=black] (0,0) circle (1pt);
\draw [fill=black] (1,0) circle (1pt);
\draw [fill=black] (1,1) circle (1pt);
\draw [fill=black] (0,1) circle (1pt);
\draw [fill=black] (0.5090377499796499,1.3861738649290605) circle (1pt);
\draw [fill=black] (0.49945650687446763,-0.39225270822641545) circle (1pt);
\draw [fill=black] (1.49925108788594,0.39590652351985406) circle (1pt);
\draw [fill=black] (2,0) circle (1pt);
\draw [fill=black] (1,-1) circle (1pt);
\draw [fill=black] (2,-1) circle (1pt);
\draw [fill=black] (1.5065488585502573,-1.4066428305665222) circle (1pt);
\draw [fill=black] (3.5,0) circle (1pt);
\draw [fill=black] (4.5,0) circle (1pt);
\draw [fill=black] (3.5,1) circle (1pt);
\draw [fill=black] (4.5,1) circle (1pt);
\draw [fill=black] (4.002386425746781,1.3957011045313252) circle (1pt);
\draw [fill=black] (4.002386425746781,-0.4141460202193678) circle (1pt);

\draw[color=black] (2.7,0) node {$...$};

%%%%hex1
\draw[color=black] (0.5,1.7) node {\tiny$x^{1}_{5}$};

\draw[color=black] (-0.2,1) node {\tiny$x^{1}_{6}$};

\draw[color=black] (0.5,-0.7) node {\tiny$x^{1}_{2}$};

\draw[color=black] (-0.2,0) node {\tiny$x^{1}_{1}$};

\draw[color=black] (0.7,0.1) node {\tiny$x^{1}_{3}$};

\draw[color=black] (0.6,0.9) node {\tiny$x^{1}_{4}$};

%%hex2
\draw[color=black] (1.2,-0.2) node {\tiny$x^{2}_{1}$};

\draw[color=black] (1.5,0.7) node {\tiny$x^{2}_{2}$};

\draw[color=black] (2.2,0.2) node {\tiny$x^{2}_{3}$};

\draw[color=black] (2.2,-1.2) node {\tiny$x^{2}_{4}$};

\draw[color=black] (1.5,-1.7) node {\tiny$x^{2}_{5}$};

\draw[color=black] (0.7,-1.2) node {\tiny$x^{2}_{6}$};

%%%%hexn
\draw[color=black] (4,1.7) node {\tiny$x^{n}_{5}$};

\draw[color=black] (3.3,1) node {\tiny$x^{n}_{6}$};

\draw[color=black] (4,-0.7) node {\tiny$x^{n}_{2}$};

\draw[color=black] (3.3,0) node {\tiny$x^{n}_{1}$};

\draw[color=black] (4.8,1) node {\tiny$x^{n}_{4}$};

\draw[color=black] (4.8,0) node {\tiny$x^{n}_{3}$};

\end{scriptsize}
\end{tikzpicture}

 }%
%\vskip 1 cm
\caption{The graphs $H(n)$ whose all vertices are labelled.}
\label{mh}
\end{figure}

We first name all the vertices of $H(n)$ as shown in Figure \ref{mh}.

%\vskip 1 cm
\vskip 5 pt

\noindent Then, we recall that
\vskip 5 pt

\indent $h(n)$ = the number of all maximal independent sets of $H(n)$

\noindent and

\indent $H(n, k)$ = the number of maximal independent sets containing $k$ vertices of $H(n)$.
\vskip 5 pt

\noindent Thus,
\begin{align}\label{Hn}
h(n) = \sum_{k \geq 0}h(n, k).    
\end{align}
\vskip 5 pt

\noindent Further, we let 
\begin{align}
H(x) = \sum_{n \geq 0}h(n)x^{n}\notag
\end{align}
\noindent be the generating function of $h(n)$ and we let
\begin{align}
H(x, y) = \sum_{n \geq 0}\sum_{k \geq 0}h(n, k)x^{n}y^{k}\notag
\end{align}
\noindent be the bi-variate generating function of $h(n, k)$. It is worth noting that, when $y = 1$, we have
\begin{align}\label{txytotnh}
H(x, 1) = \sum_{n \geq 0}(\sum_{k \geq 0}h(n, k)(1)^{k})x^{n} = \sum_{n \geq 0}h(n)x^{n} = H(x).
\end{align}

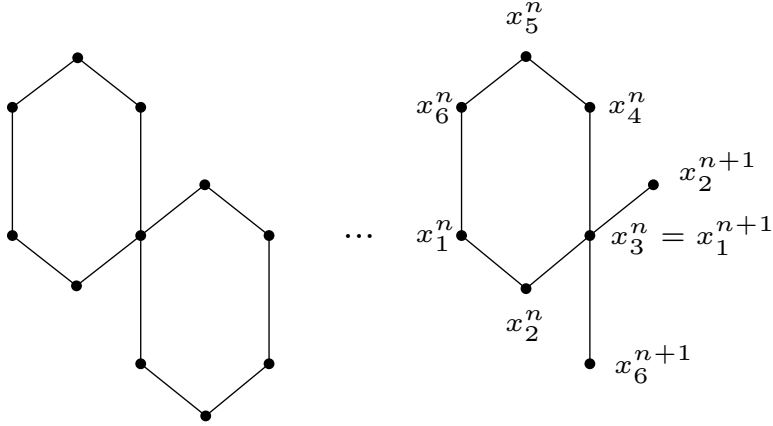
\begin{figure}[H]
\centering
\definecolor{ududff}{rgb}{0.30196078431372547,0.30196078431372547,1}
\resizebox{0.8\textwidth}{!}{%

\begin{tikzpicture}[line cap=round,line join=round,>=triangle 45,x=1cm,y=1cm]
\clip(-0.42006259682951436,-2.06344219035508) rectangle (6.717157112872825,2.344411290892576);
\draw [line width=0.3pt] (0,1)-- (0,0);
\draw [line width=0.3pt] (0,1)-- (0.5090377499796499,1.3861738649290605);
\draw [line width=0.3pt] (0.5090377499796499,1.3861738649290605)-- (1,1);
\draw [line width=0.3pt] (1,1)-- (1,0);
\draw [line width=0.3pt] (0,0)-- (0.49945650687446763,-0.39225270822641545);
\draw [line width=0.3pt] (0.49945650687446763,-0.39225270822641545)-- (1,0);
\draw [line width=0.3pt] (1,0)-- (1.49925108788594,0.39590652351985406);
\draw [line width=0.3pt] (1.49925108788594,0.39590652351985406)-- (2,0);
\draw [line width=0.3pt] (1,0)-- (1,-1);
\draw [line width=0.3pt] (1,-1)-- (1.5065488585502573,-1.4066428305665222);
\draw [line width=0.3pt] (1.5065488585502573,-1.4066428305665222)-- (2,-1);
\draw [line width=0.3pt] (2,0)-- (2,-1);
\draw [line width=0.3pt] (3.5,1)-- (3.5,0);
\draw [line width=0.3pt] (3.5,1)-- (4.002386425746781,1.3957011045313252);
\draw [line width=0.3pt] (4.002386425746781,1.3957011045313252)-- (4.5,1);
\draw [line width=0.3pt] (4.5,1)-- (4.5,0);
\draw [line width=0.3pt] (4.5,0)-- (4.002386425746781,-0.4141460202193678);
\draw [line width=0.3pt] (3.5,0)-- (4.002386425746781,-0.4141460202193678);
\draw [line width=0.3pt] (4.5,0)-- (4.994883236093936,0.39590652351985367);
\draw [line width=0.3pt] (4.5,0)-- (4.5,-1);
\begin{scriptsize}
\draw [fill=black] (0,0) circle (1pt);
\draw [fill=black] (1,0) circle (1pt);
\draw [fill=black] (1,1) circle (1pt);
\draw [fill=black] (0,1) circle (1pt);
\draw [fill=black] (0.5090377499796499,1.3861738649290605) circle (1pt);
\draw [fill=black] (0.49945650687446763,-0.39225270822641545) circle (1pt);
\draw [fill=black] (1.49925108788594,0.39590652351985406) circle (1pt);
\draw [fill=black] (2,0) circle (1pt);
\draw [fill=black] (1,-1) circle (1pt);
\draw [fill=black] (2,-1) circle (1pt);
\draw [fill=black] (1.5065488585502573,-1.4066428305665222) circle (1pt);
\draw [fill=black] (3.5,0) circle (1pt);
\draw [fill=black] (4.5,0) circle (1pt);
\draw [fill=black] (3.5,1) circle (1pt);
\draw [fill=black] (4.5,1) circle (1pt);
\draw [fill=black] (4.002386425746781,1.3957011045313252) circle (1pt);
\draw [fill=black] (4.002386425746781,-0.4141460202193678) circle (1pt);
\draw [fill=black] (4.994883236093936,0.39590652351985367) circle (1pt);
\draw [fill=black] (4.5,-1) circle (1pt);

\draw[color=black] (2.7,0) node {$...$};

%%%%hexn
\draw[color=black] (4,1.7) node {\tiny$x^{n}_{5}$};

\draw[color=black] (3.3,1) node {\tiny$x^{n}_{6}$};

\draw[color=black] (4,-0.7) node {\tiny$x^{n}_{2}$};

\draw[color=black] (3.3,0) node {\tiny$x^{n}_{1}$};

\draw[color=black] (4.8,1) node {\tiny$x^{n}_{4}$};

\draw[color=black] (5.3,0) node {\tiny$x^{n}_{3}=x^{n + 1}_{1}$};

\draw[color=black] (5.5,0.5) node {\tiny$x^{n+1}_{2}$};

\draw[color=black] (5,-1) node {\tiny$x^{n+1}_{6}$};

\end{scriptsize}
\end{tikzpicture}

 }%
%\vskip -1 cm
\caption{The graph $\overline{H}(n)$.}
\label{hbar}
\end{figure}

\begin{figure}[H]
\centering
\definecolor{ududff}{rgb}{0.30196078431372547,0.30196078431372547,1}
\resizebox{0.8\textwidth}{!}{%

\begin{tikzpicture}[line cap=round,line join=round,>=triangle 45,x=1cm,y=1cm]
\clip(-0.42006259682951436,-2.0634421903550786) rectangle (6.717157112872825,2.344411290892576);
\draw [line width=0.3pt] (0,1)-- (0,0);
\draw [line width=0.3pt] (0,1)-- (0.5090377499796499,1.3861738649290605);
\draw [line width=0.3pt] (0.5090377499796499,1.3861738649290605)-- (1,1);
\draw [line width=0.3pt] (1,1)-- (1,0);
\draw [line width=0.3pt] (0,0)-- (0.49945650687446763,-0.39225270822641545);
\draw [line width=0.3pt] (0.49945650687446763,-0.39225270822641545)-- (1,0);
\draw [line width=0.3pt] (1,0)-- (1.49925108788594,0.39590652351985406);
\draw [line width=0.3pt] (1.49925108788594,0.39590652351985406)-- (2,0);
\draw [line width=0.3pt] (1,0)-- (1,-1);
\draw [line width=0.3pt] (1,-1)-- (1.5065488585502573,-1.4066428305665222);
\draw [line width=0.3pt] (1.5065488585502573,-1.4066428305665222)-- (2,-1);
\draw [line width=0.3pt] (2,0)-- (2,-1);
\draw [line width=0.3pt] (3.5,1)-- (3.5,0);
\draw [line width=0.3pt] (3.5,1)-- (4.002386425746781,1.3957011045313252);
\draw [line width=0.3pt] (4.002386425746781,1.3957011045313252)-- (4.5,1);
\draw [line width=0.3pt] (4.5,1)-- (4.5,0);
\draw [line width=0.3pt] (4.5,0)-- (4.002386425746781,-0.4141460202193678);
\draw [line width=0.3pt] (3.5,0)-- (4.002386425746781,-0.4141460202193678);
\draw [line width=0.3pt] (4.5,0)-- (4.994883236093936,0.39590652351985367);
\draw [line width=0.3pt] (4.5,0)-- (4.5,-1);
\draw [line width=0.3pt] (4.994883236093936,0.39590652351985367)-- (5.5,0);
\draw [line width=0.3pt] (5.5,0)-- (5.5,-1);
\begin{scriptsize}
\draw [fill=black] (0,0) circle (1pt);
\draw [fill=black] (1,0) circle (1pt);
\draw [fill=black] (1,1) circle (1pt);
\draw [fill=black] (0,1) circle (1pt);
\draw [fill=black] (0.5090377499796499,1.3861738649290605) circle (1pt);
\draw [fill=black] (0.49945650687446763,-0.39225270822641545) circle (1pt);
\draw [fill=black] (1.49925108788594,0.39590652351985406) circle (1pt);
\draw [fill=black] (2,0) circle (1pt);
\draw [fill=black] (1,-1) circle (1pt);
\draw [fill=black] (2,-1) circle (1pt);
\draw [fill=black] (1.5065488585502573,-1.4066428305665222) circle (1pt);
\draw [fill=black] (3.5,0) circle (1pt);
\draw [fill=black] (4.5,0) circle (1pt);
\draw [fill=black] (3.5,1) circle (1pt);
\draw [fill=black] (4.5,1) circle (1pt);
\draw [fill=black] (4.002386425746781,1.3957011045313252) circle (1pt);
\draw [fill=black] (4.002386425746781,-0.4141460202193678) circle (1pt);
\draw [fill=black] (4.994883236093936,0.39590652351985367) circle (1pt);
\draw [fill=black] (4.5,-1) circle (1pt);
\draw [fill=black] (5.5,0) circle (1pt);
\draw [fill=black] (5.5,-1) circle (1pt);

\draw[color=black] (2.7,0) node {$...$};

%%%%hexn
\draw[color=black] (4,1.7) node {\tiny$x^{n}_{5}$};

\draw[color=black] (3.3,1) node {\tiny$x^{n}_{6}$};

\draw[color=black] (4,-0.7) node {\tiny$x^{n}_{2}$};

\draw[color=black] (3.3,0) node {\tiny$x^{n}_{1}$};

\draw[color=black] (4.8,1) node {\tiny$x^{n}_{4}$};

\draw[color=black] (4.8,0) node {\tiny$x^{n}_{3}$};

\draw[color=black] (5.5,0.5) node {\tiny$x^{n+1}_{2}$};

\draw[color=black] (5,-1) node {\tiny$x^{n+1}_{6}$};

\draw[color=black] (6,0) node {\tiny$x^{n+1}_{3}$};

\draw[color=black] (6,-1) node {\tiny$x^{n+1}_{4}$};

\end{scriptsize}
\end{tikzpicture}

 }%
%\vskip -1 cm
\caption{The graph $\tilde{H}(n)$.}
\label{mhtilde}
\end{figure}
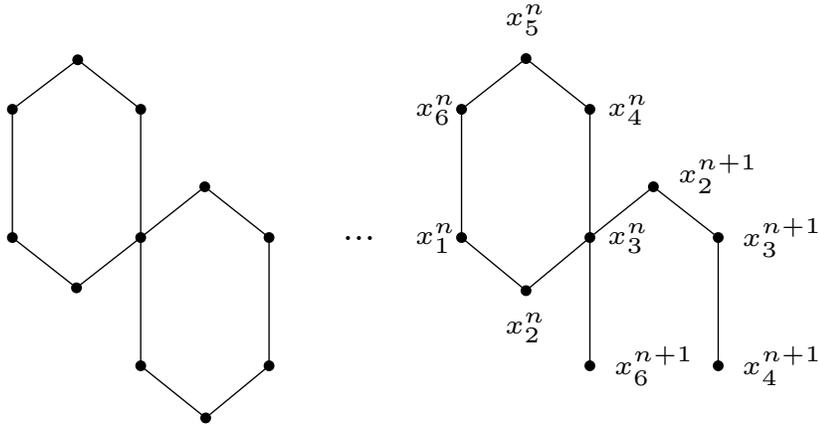

\noindent Next, we let $\overline{H}(n)$ be constructed from $H(n)$ by joining two vertices to a vertex at distance two from the cut vertex of the $n^{th}$ hexagon. Further, we let $\overline{H}(n)$ be constructed from $\overline{H}(n)$ by joining an end vertex of a path of length one to a vertex of degree one of $\overline{H}(n)$. The graphs $\overline{H}(n)$ and $\tilde{H}(n)$ are shown in Figures \ref{hbar} and \ref{mhtilde}, respectively.

%\vskip -2 cm
\vskip 5 pt

%\vskip -2 cm
\vskip 5 pt

\noindent Then, we let
\vskip 5 pt

\indent $\overline{h}(n, k)$ = the number of maximal independent sets containing $k$ vertices of $\overline{H}(n)$,
\vskip 5 pt

\indent $\tilde{h}(n, k)$ = the number of maximal independent sets containing $k$ vertices of $\tilde{H}(n)$,
\vskip 5 pt

\noindent and let
\begin{align}
\overline{H}(x, y) = \sum_{n \geq 0}\sum_{k \geq 0}\overline{h}(n, k)x^{n}y^{k},\notag
\end{align}
\begin{align}
\tilde{H}(x, y) = \sum_{n \geq 0}\sum_{k \geq 0}\tilde{h}(n, k)x^{n}y^{k}.\notag
\end{align}
\noindent Now, we are ready to prove Theorem \ref{bivariate-Hh}.

\noindent \emph{Proof of Theorem \ref{bivariate-Hh}} Let $D$ be a maximal independent set of $H(n)$ containing $k$ vertices. We distinguish $2$ cases. 
\vskip 5 pt

\noindent \textbf{Case 1:} $x^{n}_{4} \in D$\\
\indent Thus, $x^{n}_{3}, x^{n}_{5} \notin D$. Removing $x^{n}_{3}, x^{n}_{4} , x^{n}_{5}$ from $H(n)$ result in $\overline{H}(n - 1)$. There are $\overline{h}(n - 1, k - 1)$ possibilities of $D$.
\vskip 5 pt

\noindent \textbf{Case 2:} $x^{n}_{4} \notin D$\\
\indent We further have the following $3$ subcases.
\vskip 5 pt

\noindent \textbf{Subcase 2.1:} $n^{n}_{3} \notin D$ and $x^{n}_{5} \in D$.\\
\indent Thus, $x^{n}_{6} \notin D$ and $x^{n}_{2} \in D$. This implies that $x^{n}_{1} \notin D$. Removing all the vertices of the $n^{th}$ hexagon results in $\tilde{H}(n - 2)$. We have that there are $\tilde{h}(n - 2, k - 2)$ possibilities of $D$.
\vskip 5 pt

\noindent \textbf{Subcase 2.2:} $n^{n}_{3} \in D$ and $x^{n}_{5} \notin D$.\\
\indent Thus, $x^{n}_{2} \notin D$ and $x^{n}_{6} \in D$. This implies that $x^{n}_{1} \notin D$. Removing all the vertices of the $n^{th}$ hexagon results in $\tilde{H}(n - 2)$. We have that there are $\tilde{h}(n - 2, k - 2)$ possibilities of $D$.
\vskip 5 pt

\noindent \textbf{Subcase 2.3:} $n^{n}_{3} \in D$ and $x^{n}_{5} \in D$.\\
\indent Thus, $x^{n}_{2}, x^{n}_{6} \notin D$. Removing the vertices $x^{n}_{2}, ..., x^{n}_{6}$ results in $H(n - 1)$. We have that there are $h(n - 1, k - 2)$ possibilities of $D$.
\vskip 5 pt

\indent From, the two cases, we have that
\begin{align}\label{mh11}
h(n, k) = \overline{h}(n - 1, k - 1) + 2\Tilde{h}(n - 2, k - 2) +  h(n - 1, k - 2).
\end{align}
\noindent For $n \geq 2 $ and $k \geq 2 $, we multiply $x^{n}y^{k}$ throughout (\ref{m11}) and sum over all $x^{n}y^{k}$. Thus, we have that

\begin{align}\label{hhhhh}
\sum_{n \geq 2}\sum_{k \geq 2}h(n, k)x^{n}y^{k} &= \sum_{n \geq 2}\sum_{k \geq 2}\overline{h}(n - 1, k - 1)x^{n}y^{k} + 2\sum_{n \geq 2}\sum_{k \geq 2}\Tilde{h}(n - 2, k - 2)x^{n}y^{k}\notag\\
& \text{  }+ \sum_{n \geq 2}\sum_{k \geq 2}{h}(n - 1, k - 2)x^{n}y^{k}.
\end{align}

\noindent We first consider the term $\displaystyle{\sum_{n \geq 2}\sum_{k \geq 2}h(n, k)x^{n}y^{k}}$. It can be checked that

\begin{center}
    $h(0, 0) = 1, h(0, 1) = 0, h(0, 2) = 0$, and $ h(0, k) = 0$ for all $k \geq 3$
    \vskip 5 pt
    $h(1, 0) = 0, h(1, 1) = 0, h(1, 2) = 3 ,h(1, 3) = 2$, and $ h(1, k) = 0$ for all $k \geq 4$
    \vskip 5 pt
\end{center}

\begin{align}\label{mh1}
\sum_{n \geq 2}\sum_{k \geq 2}h(n, k)x^{n}y^{k} &= \sum_{n \geq 2}\sum_{k \geq 2}h(n, k)x^{n}y^{k}  + h(0, 0)  + h(1, 2)xy^{2} + h(1, 3)xy^{3} -h(0, 0)\notag\\
&\text{  }-h(1, 2)xy^{2} - h(1, 3)xy^{3}\notag\\
                  &= \sum_{n \geq 0}\sum_{k \geq 0}h(n, k)x^{n}y^{k} - 1 - 3xy^{2} - 2xy^{3} \notag\\
                  &= H(x,y) - 1 - 3xy^2 - 2xy^{3}
\end{align}

\noindent Now, we consider the term $\displaystyle{\sum_{n \geq 2}\sum_{k \geq 2}\overline{h}(n - 1, k - 1)x^{n}y^{k}}$. clearly,

\begin{center}
    $\overline{h}(0, 0) = 0, \overline{h}(0, 1) = 1 , \overline{h}(0, 2) = 1$ and $\overline{h}(0, k) = 0$ for all $k \geq 2$
    \vskip 5 pt
    $\overline{h}(1, 0) = 0, \overline{h}(1, 1) = 0$ and $\overline{h}(1, k) = 0$ for all $k \geq 2$
    \vskip 5 pt
\end{center}

\begin{align}\label{mhbar1}
\sum_{n \geq 2}\sum_{k \geq 2}\overline{h}(n - 1, k - 1)x^{n}y^{k}
    &= xy(\sum_{n \geq 2}\sum_{k \geq 2}\overline{h}(n - 1, k - 1)x^{n - 1}y^{k - 1})\notag\\
    &= xy(\sum_{n \geq 2}\sum_{k \geq 2}\overline{h}(n - 1, k - 1)x^{n - 1}y^{k - 1} + \overline{h}(0,1)y + \overline{h}(0,2)y^{2}\notag\\
    & \text{  }- \overline{h}(0,1)y - \overline{h}(0,2)y^{2})\notag\\
    &= xy(\sum_{n \geq 0}\sum_{k \geq 0}\overline{h}(n, k)x^{n}y^{k} - y - y^{2})\notag\\
    &= xy\overline{H}(x, y) - xy^{2} - xy^{3}.
\end{align}

\noindent Next, we consider the term $\displaystyle{\sum_{n \geq 2}\sum_{k \geq 2}\Tilde{2h}(n - 2, k - 2)x^{n}y^{k}}$. It can be checked that, 
\begin{center}
    $\Tilde{h}(0, 0) = 0, \Tilde{h}(0, 1) = 0 , \Tilde{h}(0, 2) = 3, \Tilde{h}(0,3) =1$ and $\Tilde{h}(0, k) = 0$ for all $k \geq 4$
    \vskip 5 pt
    $\Tilde{h}(1, 0) = 0, \Tilde{h}(1, 1) = 0, \Tilde{h}(1,2) = 0,\Tilde{h}(1, 3) = 2 ,\Tilde{h}(1, 4) = 5 , \Tilde{h}(1, 5) = 4 ,\Tilde{h}(1, 6) = 1$ and
    $\Tilde{h}(1, k) = 0$ for all $k \geq 7$
    \vskip 5 pt
\end{center}
\noindent Thus,

\begin{align}\label{mhtil1}
\sum_{n \geq 2}\sum_{k \geq 2}2\Tilde{h}(n - 2, k - 2)x^{n}y^{k}
    &= 2x^{2}y^{2}(\sum_{n \geq 2}\sum_{k \geq 2}\Tilde{h}(n - 2, k - 2)x^{n - 2}y^{k - 2})\notag\\
    &= 2x^{2}y^{2}(\sum_{n \geq 0}\sum_{k \geq 0}\Tilde{h}(n , k)x^{n}y^{k}\notag\\
    &= 2x^{2}y^{2}\Tilde{H}(x, y).
\end{align}

\noindent Finally, we consider the term $\displaystyle{\sum_{n \geq 2}\sum_{k \geq 2}h(n - 1, k - 2)x^{n}y^{k}}$.
\begin{align}\label{mhh}
\sum_{n \geq 2}\sum_{k \geq 2}h(n - 1, k - 2)x^{n}y^{k}
    &= xy^{2}(\sum_{n \geq 2}\sum_{k \geq 2}h(n - 1, k - 2)x^{n - 1}y^{k - 2})\notag\\
    &= xy^{2}(\sum_{n \geq 1}\sum_{k \geq 0}h(n , k)x^{n}y^{k} + h(0,0) - h(0,0))\notag\\
    &= xy^{2}H(x, y) - xy^{2}.
\end{align}

\noindent Plugging (\ref{mh1}), (\ref{mhbar1}), (\ref{mhtil1}) and (\ref{mhh}) to (\ref{hhhhh}), we have 
\begin{align}
H(x,y) - 1- 3xy^{2} - 2xy^{3} &= xy\overline{H}(x, y) - xy^{2} - xy^{3} + 2x^{2}y^{2}\Tilde{H}(x, y) + xy^{2}H(x, y) - xy^{2}.\notag
\end{align}

\noindent which can be solved that

\begin{align}\label{semifinalmh4}
H(x,y) &= 1 + xy^2 + xy^{3} + xy\overline{H}(x, y) + 2x^{2}y^{2}\Tilde{H}(x, y) + xy^{2}H(x, y).
\end{align}
\vskip 5 pt

\indent Next, we will find the recurrence relation of $\overline{h}(n, k)$. Let $D$ be a maximal independent set of $\overline{H}(n)$ containing $k$ vertices. We distinguish $2$ cases. 
\vskip 5 pt

\noindent \textbf{Case 1:} $x^{n}_{3} \in D$\\
\indent Thus, $x^{n}_{2}, x^{n}_{4}, x^{n+1}_{2}, x^{n+1}_{6} \notin D$. We further distinguish $2$ subcases.
\vskip 5 pt

\noindent \textbf{Subcase 1.1:} $x^{n}_{5} \in D$\\
\indent Thus, $x^{n}_{6} \notin D$. Removing $x^{n}_{2}, ..., x^{n}_{6}$ and $x^{n+1}_{2}, x^{n+1}_{6}$ from $\overline{H}(n)$ results in $H(n - 1)$. There are $h(n - 1, k - 2)$ possibilities of $D$.
\vskip 5 pt

\noindent \textbf{Subcase 1.2:} $x^{n}_{5} \notin D$\\
\indent By the maximality of $D$, $x^{n}_{6} \in D$. Thus, $x^{n}_{1} \notin D$. Removing all vertices of the $n^{th}$ hexagon and $x^{n+1}_{2}, x^{n+1}_{6}$ from $\overline{H}(n)$ result in $\tilde{H}(n - 2)$. There are $\tilde{h}(n - 2, k - 2)$ possibilities of $D$.
\vskip 5 pt

\noindent \textbf{Case 2:} $x^{n}_{3} \notin D$\\
\indent By the maximality of $D$, $x^{n+1}_{2}, x^{n+1}_{6} \in D$. Removing the vertices $x^{n}_{3}, x^{n+1}_{2}, x^{n+1}_{6}$ from $\overline{H}(n)$ result in $\tilde{H}(n - 1)$. There are $\tilde{h}(n - 1, k - 2)$ possibilities of $D$.
\vskip 5 pt

\indent From the two cases, we have that
\begin{align}\label{mhbar2}
\overline{h}(n, k) = h(n - 1, k - 2)
+ \Tilde{h}(n - 2, k - 2) + \Tilde{h}(n - 1, k - 2)
\end{align}
\vskip 5 pt
\noindent For $n \geq 2, k \geq 2$, we multiply $x^{n}y^{k}$ throughout (\ref{mhbar2}) and sum over all $n$ and $k$. We have that
\begin{align}\label{mhbar22}
\displaystyle{
\sum_{n \geq 2}\sum_{k \geq 2}\overline{h}(n, k) = \sum_{n \geq 2}\sum_{k \geq 2}h(n - 1, k - 2) + \sum_{n \geq 2}\sum_{k \geq 2}\Tilde{h}(n - 2, k - 2) + \sum_{n \geq 2}\sum_{k \geq 2}\Tilde{h}(n - 1, k - 2).}
\end{align}

\noindent We first consider the term $\displaystyle{\sum_{n \geq 2}\sum_{k \geq 2}\overline{h}(n , k )x^{n}y^{k}}$. It can be checked that, 

\begin{center}
    $\overline{h}(0, 0) = 0, \overline{h}(0, 1) = 1, \overline{h}(0, 2) =1$ and $\overline{h}(0, k) = 0$ for all $k \geq 3$
    \vskip 5 pt
    $\overline{h}(1, 0) = 0, \overline{h}(1, 1) = 0, \overline{h}(1, 2) =1, \overline{h}(1, 3) = 1, \overline{h}(1, 4) = 3, \overline{h}(1, 5) = 1$ and $\overline{h}(1, k) = 0$ for all $k \geq 6$
    \vskip 5 pt
\end{center}

\begin{align}\label{mhbar3}
\sum_{n \geq 2}\sum_{k \geq 2}\overline{h}(n, k)x^{n}y^{k}
    &= \sum_{n \geq 2}\sum_{k \geq 2}\overline{h}(n,k)x^{n}y^{k} + \overline{h}(0, 1)y +\overline{h}(0, 2)y^{2} + \overline{h}(1, 2)xy^{2}\notag\\
    &+ \overline{h}(1, 3)xy^{3} + \overline{h}(1,4)xy^{4} +\overline{h}(1, 5)xy^{5} - \overline{h}(0, 1)y - \overline{h}(0, 2)y^{2}\notag\\
    &- \overline{h}(1, 2)xy^{2} - \overline{h}(1, 3)xy^{3} - \overline{h}(1,4)xy^{4} - \overline{h}(1, 5)xy^{5}\notag\\
    &= \sum_{n \geq 0}\sum_{k \geq 0}\overline{h}(n , k )x^{n}y^{k} - y - y^{2} - xy^{2} - xy^{3} - 3xy^{4} - xy^{5}\notag\\
    &= \overline{H}(x, y) - y - y^{2} - xy^{2} - xy^{3} - 3xy^{4} - xy^{5}
\end{align}

\noindent We next consider the term $\displaystyle{\sum_{n \geq 2}\sum_{k \geq 2}h(n - 1 , k - 2)x^{n}y^{k}}$, clearly

\begin{align}\label{mh2}
\sum_{n \geq 2}\sum_{k \geq 2}h(n-1, k-2)x^{n}y^{k} 
&= xy^2(\sum_{n \geq 2}\sum_{k \geq 2}h(n-1, k-2)x^{n-1}y^{k-2} + h(0, 0) -h(0,0))\notag\\
                     &= xy^2(\sum_{n \geq 1}\sum_{k \geq 0}h(n, k)x^{n}y^{k} - 1)\notag\\
                  &= xy^2H(x,y) - xy^2
\end{align}

\noindent We next consider the term $\displaystyle{\sum_{n \geq 2}\sum_{k \geq 2}\Tilde{h}(n - 2 , k - 2 )x^{n}y^{k}}$.

\begin{align}\label{mhtil2}
\sum_{n \geq 2}\sum_{k \geq 2}\Tilde{h}(n - 2 , k - 2)x^{n}y^{k}
    &= x^{2}y^{2}(\sum_{n \geq 2}\sum_{k \geq 2}\Tilde{h}(n - 2, k - 2 )x^{n - 2}y^{k - 2})\notag\\
    &= x^{2}y^{2}\Tilde{H}(x, y)
\end{align}    
\vskip 5 pt

\noindent Finally, We next consider the term $\displaystyle{\sum_{n \geq 2}\sum_{k \geq 2}\Tilde{h}(n - 1 , k - 2)x^{n}y^{k}}$ , Clearly
\begin{align}\label{mhtil3}
\sum_{n \geq 2}\sum_{k \geq 2}\Tilde{h}(n - 1 , k - 2)x^{n}y^{k}
    &= xy^{2}(\sum_{n \geq 2}\sum_{k \geq 2}\Tilde{h}(n - 1 , k -2 )x^{n - 1}y^{k - 2} + \Tilde{h}(0, 2) + \Tilde{h}(0, 3)\notag\\ 
    &- \Tilde{h}(0, 2) - \Tilde{h}(0, 3))\notag\\
    &= xy^{2}(\sum_{n \geq 0}\sum_{k \geq 0}\Tilde{h}(n , k )x^{n}y^{k} - 3y^{2} - y^{3})\notag\\
    &= xy^{2}\Tilde{H}(x, y) - 3xy^{4} - xy^{5}
\end{align}    
\vskip 5 pt

\noindent Plugging (\ref{mhbar3}) , (\ref{mh2}) , (\ref{mhtil2}) and (\ref{mhtil3}) to (\ref{mhbar22}), we have 
\begin{align}
\overline{H}(x, y) - y &- y^{2} - xy^{2} - xy^{3} - 3xy^{4} - xy^{5} \notag\\
&= xy^2H(x,y) - xy^2 + x^{2}y^{2}\Tilde{H}(x, y) + xy^{2}\Tilde{H}(x, y) - 3xy^{4} - xy^{5}.
\end{align}
\noindent which can be solved that
\begin{align}\label{semifinalmhbar4}
\overline{H}(x, y) &= xy^2H(x,y) + (x^{2}y^{2}+xy^{2})\tilde{H}(x,y) + y + y^{2} +xy^{3}.
\end{align}
\vskip 5 pt

\indent Next, we will find the recurrence relation of $\tilde{h}(n, k)$. Let $D$ be a maximal independent set of $\tilde{H}(n)$ containing $k$ vertices. We distinguish $2$ cases. 
\vskip 5 pt

\noindent \textbf{Case 1:} $x^{n}_{3} \in D$\\
\indent Thus, $x^{n}_{2}, x^{n}_{4}, x^{n+1}_{2}, x^{n+1}_{6} \notin D$. Clearly, either $x^{n+1}_{3} \in D$ or $x^{n+1}_{4} \in D$. Further, we distinguish $2$ subcases.
\vskip 5 pt

\noindent \textbf{Subcase 1.1:} $x^{n}_{5} \in D$\\
\indent Thus, $x^{n}_{6} \notin D$. Removing $x^{n}_{2}, ..., x^{n}_{6}$ and $x^{n+1}_{2}, x^{n+1}_{3}, x^{n+1}_{4}, x^{n+1}_{6}$ from $tilde{H}(n)$ result in $H(n - 1)$. Since either $x^{n+1}_{3}$ or $x^{n+1}_{4}$ is in $D$, there are $2h(n - 1, k - 3)$ possibilities of $D$.
\vskip 5 pt

\noindent \textbf{Subcase 1.2:} $x^{n}_{5} \notin D$\\
\indent By the maximality of $D$, $x^{n}_{6} \in D$. THus, $x^{n}_{1} \notin D$. Removing all vertices of the $n^{th}$ hexagon and $x^{n+1}_{2}, x^{n+1}_{3}, x^{n+1}_{4}, x^{n+1}_{6}$ from $\tilde{H}(n)$ result in $\tilde{H}(n - 2)$. Since either $x^{n+1}_{3}$ or $x^{n+1}_{4}$ is in $D$, , there are $2\tilde{h}(n - 2, k - 3)$ possibilities of $D$.
\vskip 5 pt

\noindent \textbf{Case 2:} $x^{n}_{3} \notin D$\\
\indent By the maximality of $D$, $x^{n+1}_{6} \in D$. We further distinguish $2$ subcases.
\vskip 5 pt

\noindent \textbf{Subcase 2.1:} $x^{n+1}_{2} \in D$\\
\indent Thus, $x^{n+1}_{3} \notin D$. By maximality of $D$, $x^{n+1}_{4} \in D$. Removing the vertices $x^{n}_{3}, x^{n+1}_{2}, x^{n+1}_{3},$ $x^{n+1}_{4}, x^{n+1}_{6}$ from $\tilde{H}(n)$ result in $\tilde{H}(n - 1)$. There are $\tilde{h}(n - 1, k - 3)$ possibilities of $D$.
\vskip 5 pt

\noindent \textbf{Subcase 2.2:} $x^{n+1}_{2} \notin D$\\
\indent By maximality of $D$, $x^{n+1}_{3} \in D$. Thus, $x^{n+1}_{4} \notin D$. Removing the vertices $x^{n}_{3}, x^{n+1}_{2}, x^{n+1}_{3},$ $x^{n+1}_{4}, x^{n+1}_{6}$ from $\tilde{H}(n)$ result in $\tilde{H}(n - 1)$. There are $\tilde{h}(n - 1, k - 2)$ possibilities of $D$.

\indent From all the cases, we have that
\begin{align}\label{mhtil4}
\Tilde{h}(n, k) = 2h(n - 1,k - 3) + 2\Tilde{h}(n - 2, k - 3) + \Tilde{h}(n - 1, k - 3) + \Tilde{h}(n - 1,k - 2)
\end{align}
\vskip 5 pt
\noindent For $n \geq 2, k \geq 3$, we multiply $x^{n}y^{k}$ throughout (\ref{mhtil4}) and sum over all $n, k$. We have that

\begin{align}\label{mhtil44}
\sum_{n \geq 2}\sum_{k \geq 3}\Tilde{h}(n, k) = \sum_{n \geq 2}\sum_{k \geq 3}2h(n - 1,k - 3) + \sum_{n \geq 2}\sum_{k \geq 3}2\Tilde{h}(n - 2, k - 3) \notag\\
+ \sum_{n \geq 2}\sum_{k \geq 3}\Tilde{h}(n - 1, k - 3) + \sum_{n \geq 2}\sum_{k \geq 3}\Tilde{h}(n - 1,k - 2).
\end{align}

\noindent We first consider the term $\displaystyle{\sum_{n \geq 2}\sum_{k \geq 3}\Tilde{h}(n , k )x^{n}y^{k}}$. Recall that, 
\begin{center}
    $\Tilde{h}(0, 0) = 0, \Tilde{h}(0, 1) = 0 , \Tilde{h}(0, 2) = 3, \Tilde{h}(0,3) =1$ and $\Tilde{h}(0, k) = 0$ for all $k \geq 4$
    \vskip 5 pt
    $\Tilde{h}(1, 0) = 0, \Tilde{h}(1, 1) = 0, \Tilde{h}(1,2) = 0,\Tilde{h}(1, 3) = 2 ,\Tilde{h}(1, 4) = 5 , \Tilde{h}(1, 5) = 4 ,\Tilde{h}(1, 6) = 1$ and
    $\Tilde{h}(1, k) = 0$ for all $k \geq 7$.
    \vskip 5 pt
\end{center}
\noindent Thus,
\begin{align}\label{mhtil5}
\sum_{n \geq 2}\sum_{k \geq 3}\Tilde{h}(n , k)x^{n}y^{k}
    &= \sum_{n \geq 2}\sum_{k \geq 3}\Tilde{h}(n , k)x^{n}y^{k}
    + \Tilde{h}(0,2)y^{2} + \Tilde{h}(0,3)y^{3} + \Tilde{h}(1,3)xy^{3}\notag\\
    &+ \Tilde{h}(1,4)xy^{4} + \Tilde{h}(1,5)xy^{5} + \Tilde{h}(1,6)xy^{6} -\Tilde{h}(0,2)y^{2} - \Tilde{h}(0,3)y^{3}\notag\\
    &- \Tilde{h}(1,3)xy^{3} - \Tilde{h}(1,4)xy^{4} - \Tilde{h}(1,5)xy^{5} - \Tilde{h}(1,6)xy^{6}\notag\\
    &= \sum_{n \geq 0}\sum_{k \geq 0}\Tilde{h}(n , k )x^{n}y^{k} - 3y^{2} - y^{3} - 2xy^{3} - 5xy^{4} - 4xy^{5} - xy^{6}\notag\\
    &= \Tilde{H}(x, y) - 3y^{2} - y^{3} - 2xy^{3} - 5xy^{4} - 4xy^{5} - xy^{6}
\end{align}    
\vskip 5 pt

\noindent We next consider the term $\displaystyle{\sum_{n \geq 2}\sum_{k \geq 3}2h(n - 1 , k - 3 )x^{n}y^{k}}$.

\begin{align}\label{mh3}
\sum_{n \geq 2}\sum_{k \geq 3}2h(n-1, k-3)x^{n}y^{k} &= 2xy^{3}(\sum_{n \geq 2}\sum_{k \geq 3}h(n-1, k-3)x^{n-1}y^{k-3} + h(0, 0) -h(0,0))\notag\\
                  &= 2xy^{3}(\sum_{n \geq 0}\sum_{k \geq 0}h(n, k)x^{n}y^{k} - 1)\notag\\
                  &= 2xy^{3}H(x,y) - 2xy^{3}
\end{align}

\noindent Now, we consider the term $\displaystyle{\sum_{n \geq 2}\sum_{k \geq 3}2\Tilde{h}(n - 2 , k - 3)x^{n}y^{k}}$.

\begin{align}\label{mhtil6}
\sum_{n \geq 2}\sum_{k \geq 3}2\Tilde{h}(n - 2 , k - 3)x^{n}y^{k}
    &= 2x^{2}y^{3}(\sum_{n \geq 2}\sum_{k \geq 3}\Tilde{h}(n - 2 , k - 3)x^{n - 2}y^{k - 3})\notag\\
    &= 2x^{2}y^{3}\Tilde{H}(x, y)
\end{align}      
\vskip 5 pt

\noindent Next, we consider the term
$\displaystyle{\sum_{n \geq 2}\sum_{k \geq 3}\Tilde{h}(n - 1 , k - 3)x^{n}y^{k}}$.

\begin{align}\label{mhtil7}
\sum_{n \geq 2}\sum_{k \geq 3}\Tilde{h}(n - 1 , k - 3)x^{n}y^{k}
    &= xy^{3}(\sum_{n \geq 2}\sum_{k \geq 3}\Tilde{h}(n - 1 , k - 3 )x^{n - 1}y^{k - 3}\notag\\
    &+ \Tilde{h}(0,2)y^{2} + \Tilde{h}(0,3)y^{3} - \Tilde{h}(0,2)y^{2} - \Tilde{h}(0,3)y^{3})\notag\\
    &= xy^{3}(\sum_{n \geq 0}\sum_{k \geq 0}\Tilde{h}(n , k )x^{n}y^{k} - 3y^{2} - y^{3})\notag\\
    &= xy^{3}\Tilde{H}(x, y) - 3xy^{5} - xy^{6}
\end{align}

\noindent Finally, we consider the term
$\displaystyle{\sum_{n \geq 2}\sum_{k \geq 3}\Tilde{h}(n - 1 , k - 2)x^{n}y^{k}}$. It can be checked that, 

\begin{align}\label{mhtil8}
\sum_{n \geq 2}\sum_{k \geq 3}\Tilde{h}(n - 1 , k - 2)x^{n}y^{k}
    &= xy^{2}(\sum_{n \geq 2}\sum_{k \geq 3}\Tilde{h}(n - 1 , k - 2 )x^{n - 1}y^{k - 2}\notag\\
    &+ \Tilde{h}(0,2)y^{2} + \Tilde{h}(0,3)y^{3} - \Tilde{h}(0,2)y^{2} - \Tilde{h}(0,3)y^{3})\notag\\
    &= xy^{2}(\sum_{n \geq 0}\sum_{k \geq 0}\Tilde{h}(n , k )x^{n}y^{k} - 3y^{2} -y^{3})\notag\\
    &= xy^{2}\Tilde{H}(x, y) - 3xy^{4} - xy^{5}
\end{align}

\noindent Plugging (\ref{mhtil5}) , (\ref{mh3}) , (\ref{mhtil6}), (\ref{mhtil7}) , and (\ref{mhtil8}) to (\ref{mhtil44}), we have 
\begin{align}
&\Tilde{H}(x, y) - 3y^{2} - y^{3} - 2xy^{3} - 5xy^{4} - 4xy^{5} - xy^{6} =\notag\\
& 2xy^{3}H(x,y) - 2xy^{3} + 2x^{2}y^{3}\Tilde{H}(x, y) + xy^{3}\Tilde{H}(x, y) - 3xy^{5} - xy^{6} + xy^{2}\Tilde{H}(x, y) - 3xy^{4} - xy^{5}
\end{align}
\noindent which can be solved that

\begin{align}\label{semifinalmhtil9}
\Tilde{H}(x, y) = \frac{1}{1 - 2x^{2}y^{3} - xy^{3} - xy^{2}}(2xy^{3}H(x, y) + 3y^{2} + y^{3} + 2xy^{4})
\end{align}
\vskip 5 pt

\noindent By plugging (\ref{semifinalmh4}) to (\ref{semifinalmhbar4}) to (\ref{semifinalmhtil9}) , we have
\begin{align}
H(x, y) &= \frac{1-2x^{2}y^{3}  + xy^{2} - x^{3}y^{5} + x^{2}y^{5} + 5x^{2}y^{4} + xy^{3} - x^{2}y^{6} + x^{3}y^{7}}{1 - 3x^{2}y^{3} - xy^{3} - 2xy^{2} - x^{3}y^{5} + x^{2}y^{5} + x^{2}y^{4} - x^{3}y^{6}}\notag\\
\end{align}
\noindent as required.
\qed

\indent Now, we are ready to prove Theorem \ref{recurrence-h}. By Theorem \ref{bivariate-Hh} with $y = 1$, we have that
\vskip 5 pt

\noindent \emph{Proof of Theorem \ref{recurrence-h}}
\vskip 5 pt

\begin{align}
\sum_{n \geq 0}h(n)x^{n} &= H(x, 1)\notag\\
                         &= \frac{3x^{2} + 2x + 1}
                         {-2x^{3} - x^{2} - 3x + 1}
\end{align}
 
 \noindent which can  be solved that

\begin{align}\label{recursived-mh}
3x^{2} + 2x + 1\notag\\   
&= (-2x^{3} - x^{2} - 3x + 1)\sum_{n \geq 0}h(n)x^{n}\notag\\
                  &= -2\sum_{n \geq 0}h(n)x^{n + 3} - \sum_{n \geq 0}h(n)x^{n + 2} - 3\sum_{n \geq 0}h(n)x^{n + 1} + \sum_{n \geq 0}h(n)x^{n} \notag\\
                  &= -2\sum_{n \geq 3}h(n-3)x^{n} - \sum_{n \geq 2}h(n - 2)x^{n} - 3\sum_{n \geq 1}h(n - 1)x^{n} + \sum_{n \geq 0}h(n)x^{n} \notag\\
                  &= h(0) + h(1)x + h(2)x^{2} - 3h(0)x -3h(1)x^{2} - h(0)x^{2}\notag\\
                  &+ \sum_{n \geq 3}(-2h(n-3) - 3h(n-1) - h(n-2) + h(n))x^{n}
\end{align}
\noindent Because the order of the polynomial on the left hand side of (\ref{recursived-mh}) is two, the coefficients of $x^{n}$ for all $n \geq 3$ must be $0$. Thus, $h(n) - 3h(n - 1) - h(n - 2) - 2h(n - 3) = 0$ implying that

\begin{align*}
    h(n) = 3h(n - 1) + h(n - 2) + 2h(n - 3)
\end{align*}

and this completes the proof.
\qed

\subsection{Para-Hexagonal Cacti}
\begin{figure}[H]
\centering
\definecolor{ududff}{rgb}{0.30196078431372547,0.30196078431372547,1}
\resizebox{1\textwidth}{!}{%

\begin{tikzpicture}[line cap=round,line join=round,>=triangle 45,x=1cm,y=1cm]
\clip(-0.8083969725917267,-1.5837349250118342) rectangle (8.231945270112451,3.999462124674588);
\draw [line width=0.3pt] (-0.02405899207354767,0.9828713769549792)-- (0.3745474533070522,1.5702914017263894);
\draw [line width=0.3pt] (0.3745474533070522,1.5702914017263894)-- (0.9829467646774412,1.5807810450258788);
\draw [line width=0.3pt] (0.9829467646774412,1.5807810450258788)-- (1.3916045138741162,0.9767710067152169);
\draw [line width=0.3pt] (1.3916045138741162,0.9767710067152169)-- (1.8011389420376196,0.395451352183569);
\draw [line width=0.3pt] (1.3916045138741162,0.9767710067152169)-- (1.7906492987381306,1.5702914017263894);
\draw [line width=0.3pt] (1.3916045138741162,0.9767710067152169)-- (0.9934364079769307,0.37447206558459006);
\draw [line width=0.3pt] (-0.02405899207354767,0.9828713769549792)-- (0.38503709660654173,0.37447206558459006);
\draw [line width=0.3pt] (0.38503709660654173,0.37447206558459006)-- (0.9934364079769307,0.37447206558459006);
\draw [line width=0.3pt] (1.7906492987381306,1.5702914017263894)-- (2.388558966809031,1.5702914017263894);
\draw [line width=0.3pt] (2.388558966809031,1.5702914017263894)-- (2.7976550554891206,0.9723817336554896);
\draw [line width=0.3pt] (2.7976550554891206,0.9723817336554896)-- (2.3990486101085207,0.37447206558459006);
\draw [line width=0.3pt] (1.8011389420376196,0.395451352183569)-- (2.3990486101085207,0.37447206558459006);
\draw [line width=0.3pt] (4.3920808370115205,1.5807810450258788)-- (3.9829847483314316,0.9618920903560004);
\draw [line width=0.3pt] (4.3920808370115205,1.5807810450258788)-- (4.97950086178293,1.5702914017263894);
\draw [line width=0.3pt] (4.97950086178293,1.5702914017263894)-- (5.391604513874117,0.9767710067152169);
\draw [line width=0.3pt] (5.391604513874117,0.9767710067152169)-- (5.00048014838191,0.38496170888407955);
\draw [line width=0.3pt] (3.9829847483314316,0.9618920903560004)-- (4.381591193712031,0.38496170888407955);
\draw [line width=0.3pt] (4.381591193712031,0.38496170888407955)-- (5.00048014838191,0.38496170888407955);
\begin{scriptsize}
\draw [fill=black] (-0.02405899207354767,0.9828713769549792) circle (1pt);
\draw [fill=black] (1.3916045138741162,0.9767710067152169) circle (1pt);
\draw [fill=black] (2.7976550554891206,0.9723817336554896) circle (1pt);
\draw [fill=black] (3.9829847483314316,0.9618920903560004) circle (1pt);
\draw [fill=black] (5.391604513874117,0.9767710067152169) circle (1pt);
\draw [fill=black] (0.3745474533070522,1.5702914017263894) circle (1pt);
\draw [fill=black] (0.9829467646774412,1.5807810450258788) circle (1pt);
\draw [fill=black] (0.38503709660654173,0.37447206558459006) circle (1pt);
\draw [fill=black] (0.9934364079769307,0.37447206558459006) circle (1pt);
\draw [fill=black] (1.7906492987381306,1.5702914017263894) circle (1pt);
\draw [fill=black] (2.388558966809031,1.5702914017263894) circle (1pt);
\draw [fill=black] (1.8011389420376196,0.395451352183569) circle (1pt);
\draw [fill=black] (2.3990486101085207,0.37447206558459006) circle (1pt);
\draw [fill=black] (4.3920808370115205,1.5807810450258788) circle (1pt);
\draw [fill=black] (4.97950086178293,1.5702914017263894) circle (1pt);
\draw [fill=black] (4.381591193712031,0.38496170888407955) circle (1pt);
\draw [fill=black] (5.00048014838191,0.38496170888407955) circle (1pt);

\draw[color=black] (3.3,1) node {$...$};

%%hex1

\draw[color=black] (-0.2,1) node {\tiny$x^{1}_{1}$};

\draw[color=black] (1.7,1) node {\tiny$x^{2}_{1}$};

\draw[color=black] (3.1,1) node {\tiny$x^{2}_{4}$};

\draw[color=black] (1.1,1) node {\tiny$x^{1}_{4}$};

\draw[color=black] (1.1,1.8) node {\tiny$x^{1}_{3}$};

\draw[color=black] (0.5,1.8) node {\tiny$x^{1}_{2}$};

\draw[color=black] (2.4,1.8) node {\tiny$x^{2}_{3}$};

\draw[color=black] (1.8,1.8) node {\tiny$x^{2}_{2}$};

\draw[color=black] (1.1,0.1) node {\tiny$x^{1}_{5}$};

\draw[color=black] (0.5,0.1) node {\tiny$x^{1}_{6}$};

\draw[color=black] (2.4,0.1) node {\tiny$x^{2}_{5}$};

\draw[color=black] (1.8,0.1) node {\tiny$x^{2}_{6}$};

%hexn
\draw[color=black] (3.8,1) node {\tiny$x^{n}_{1}$};

\draw[color=black] (5.1,1.8) node {\tiny$x^{n}_{3}$};

\draw[color=black] (4.5,1.8) node {\tiny$x^{n}_{2}$};

\draw[color=black] (5.7,1) node {\tiny$x^{n}_{4}$};

\draw[color=black] (5.1,0.1) node {\tiny$x^{n}_{5}$};

\draw[color=black] (4.5,0.1) node {\tiny$x^{n}_{6}$};
\end{scriptsize}
\end{tikzpicture}

 }%
\vskip -3 cm
\caption{The para-hexagonal cactus $G(n)$ of $n$ hexagons.}
\label{labelphn}
\end{figure}
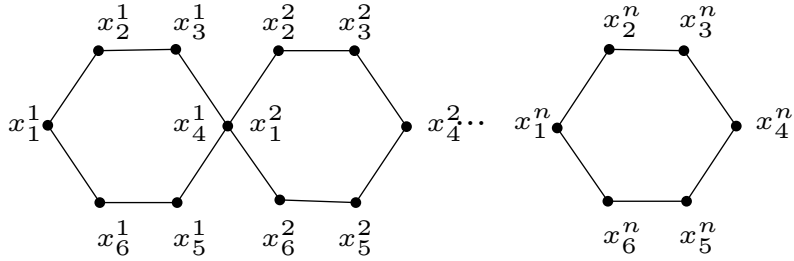

First, we may name all the vertices of $G(n)$ as shown in Figure \ref{labelphn}.

\vskip -3 cm
\vskip 5 pt

\noindent Then, we recall that
\vskip 5 pt

\indent $g(n)$ = the number of all maximal independent sets of $G(n)$

\noindent and

\indent $G(n, k)$ = the number of maximal independent sets containing $k$ vertices of $G(n)$.
\vskip 5 pt

\noindent Thus,
\begin{align}\label{Qn}
G(n) = \sum_{k \geq 0}g(n, k).    
\end{align}
\vskip 5 pt

\noindent Further, we let 
\begin{align}
G(x) = \sum_{n \geq 0}g(n)x^{n}\notag
\end{align}
\noindent be the generating function of $g(n)$ and we let
\begin{align}
G(x, y) = \sum_{n \geq 0}\sum_{k \geq 0}g(n, k)x^{n}y^{k}\notag
\end{align}
\noindent be the bi-variate generating function of $g(n, k)$. It is worth noting that, when $y = 1$, we have
\begin{align}\label{qxytoqn}
G(x, 1) = \sum_{n \geq 0}(\sum_{k \geq 0}g(n, k)(1)^{k})x^{n} = \sum_{n \geq 0}g(n)x^{n} = G(x).
\end{align}

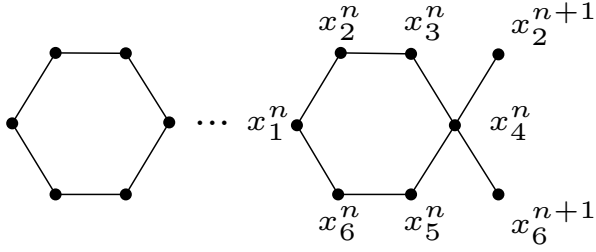
\begin{figure}[H]
\centering
\definecolor{ududff}{rgb}{0.30196078431372547,0.30196078431372547,1}
\resizebox{0.8\textwidth}{!}{%

\begin{tikzpicture}[line cap=round,line join=round,>=triangle 45,x=1cm,y=1cm]
\clip(-1.3044193314100485,-1.0954962890926483) rectangle (5.106902020096289,2.864051866643371);
\draw [line width=0.3pt] (-0.004836677106818454,0.990579241506063)-- (0.2996682258189778,0.5054414834987486);
\draw [line width=0.3pt] (0.2996682258189778,0.5054414834987486)-- (-0.004836677106818454,0);
\draw [line width=0.3pt] (1.1953037164397389,0.4855384725960651)-- (1.5024247274329101,0.9933656630324388);
\draw [line width=0.3pt] (1.5024247274329101,0.9933656630324388)-- (1.9957903904653493,0.9867313260648776);
\draw [line width=0.3pt] (1.9957903904653493,0.9867313260648776)-- (2.303237990022459,0.4855384725960651);
\draw [line width=0.3pt] (2.303237990022459,0.4855384725960651)-- (1.9957903904653493,0);
\draw [line width=0.3pt] (1.1953037164397389,0.4855384725960651)-- (1.4843193763909703,0);
\draw [line width=0.3pt] (1.4843193763909703,0)-- (1.9957903904653493,0);
\draw [line width=0.3pt] (2.303237990022459,0.4855384725960651)-- (2.6084174905302744,0.9831137451631534);
\draw [line width=0.3pt] (2.303237990022459,0.4855384725960651)-- (2.6084174905302744,0);
\draw [line width=0.3pt] (-0.8016317107961806,0.4988071465311874)-- (-0.4975752725670899,0.9933656630324388);
\draw [line width=0.3pt] (-0.8016317107961806,0.4988071465311874)-- (-0.4975752725670899,0);
\draw [line width=0.3pt] (-0.4975752725670899,0.9933656630324388)-- (-0.004836677106818454,0.990579241506063);
\draw [line width=0.3pt] (-0.4975752725670899,0)-- (-0.004836677106818454,0);
\begin{scriptsize}
\draw [fill=black] (-0.8016317107961806,0.4988071465311874) circle (1pt);
\draw [fill=black] (-0.004836677106818454,0.990579241506063) circle (1pt);
\draw [fill=black] (-0.004836677106818454,0) circle (1pt);
\draw [fill=black] (0.2996682258189778,0.5054414834987486) circle (1pt);
\draw [fill=black] (1.1953037164397389,0.4855384725960651) circle (1pt);
\draw [fill=black] (1.5024247274329101,0.9933656630324388) circle (1pt);
\draw [fill=black] (1.4843193763909703,0) circle (1pt);
\draw [fill=black] (1.9957903904653493,0) circle (1pt);
\draw [fill=black] (2.303237990022459,0.4855384725960651) circle (1pt);
\draw [fill=black] (1.9957903904653493,0.9867313260648776) circle (1pt);
\draw [fill=black] (2.6084174905302744,0.9831137451631534) circle (1pt);
\draw [fill=black] (2.6084174905302744,0) circle (1pt);
\draw [fill=black] (-0.4975752725670899,0.9933656630324388) circle (1pt);
\draw [fill=black] (-0.4975752725670899,0) circle (1pt);

\draw[color=black] (0.6,0.5) node {$...$};
%hexn
\draw[color=black] (1,0.5) node {\tiny$x^{n}_{1}$};

\draw[color=black] (1.5,1.2) node {\tiny$x^{n}_{2}$};

\draw[color=black] (2.1,1.2) node {\tiny$x^{n}_{3}$};

\draw[color=black] (2.7,0.5) node {\tiny$x^{n}_{4}$};

\draw[color=black] (1.5,-0.2) node {\tiny$x^{n}_{6}$};

\draw[color=black] (2.1,-0.2) node {\tiny$x^{n}_{5}$};

\draw[color=black] (3,-0.2) node {\tiny$x^{n+1}_{6}$};

\draw[color=black] (3,1.2) node {\tiny$x^{n+1}_{2}$};
\end{scriptsize}

\end{tikzpicture}

 }%
\vskip -1 cm
\caption{The graph $\overline{G}(n).$}
\label{qbar}
\end{figure}

\begin{figure}[H]
\centering
\definecolor{ududff}{rgb}{0.30196078431372547,0.30196078431372547,1}
\resizebox{0.9\textwidth}{!}{%

\begin{tikzpicture}[line cap=round,line join=round,>=triangle 45,x=1cm,y=1cm]
\clip(-1.248987944662694,-1.5540924731902557) rectangle (5.775500024935567,2.7841393607947578);
\draw [line width=0.3pt] (-0.004836677106818454,0.990579241506063)-- (0.2996682258189778,0.5054414834987486);
\draw [line width=0.3pt] (0.2996682258189778,0.5054414834987486)-- (-0.004836677106818454,0);
\draw [line width=0.3pt] (0.6629515060126447,0.4898164476493253)-- (0.970072517005816,0.997643638085699);
\draw [line width=0.3pt] (0.970072517005816,0.997643638085699)-- (1.4634381800382552,0.9910093011181378);
\draw [line width=0.3pt] (1.4634381800382552,0.9910093011181378)-- (1.7708857795953656,0.4898164476493253);
\draw [line width=0.3pt] (1.7708857795953656,0.4898164476493253)-- (1.4634381800382552,0);
\draw [line width=0.3pt] (0.6629515060126447,0.4898164476493253)-- (0.9519671659638762,0);
\draw [line width=0.3pt] (0.9519671659638762,0)-- (1.4634381800382552,0);
\draw [line width=0.3pt] (1.7708857795953656,0.4898164476493253)-- (2.076065280103183,0.9873917202164137);
\draw [line width=0.3pt] (1.7708857795953656,0.4898164476493253)-- (2.076065280103183,0);
\draw [line width=0.3pt] (-0.8016317107961806,0.4988071465311874)-- (-0.4975752725670899,0.9933656630324388);
\draw [line width=0.3pt] (-0.8016317107961806,0.4988071465311874)-- (-0.4975752725670899,0);
\draw [line width=0.3pt] (-0.4975752725670899,0.9933656630324388)-- (-0.004836677106818454,0.990579241506063);
\draw [line width=0.3pt] (-0.4975752725670899,0)-- (-0.004836677106818454,0);
\draw [line width=0.3pt] (2.076065280103183,0.9873917202164137)-- (2.5699067613713895,0.9891664816425414);
\draw [line width=0.3pt] (2.076065280103183,0)-- (2.5830178479797867,-0.00727610059566227);
\begin{scriptsize}
\draw [fill=black] (-0.8016317107961806,0.4988071465311874) circle (1pt);
\draw [fill=black] (-0.004836677106818454,0.990579241506063) circle (1pt);
\draw [fill=black] (-0.004836677106818454,0) circle (1pt);
\draw [fill=black] (0.2996682258189778,0.5054414834987486) circle (1pt);
\draw [fill=black] (0.6629515060126447,0.4898164476493253) circle (1pt);
\draw [fill=black] (0.970072517005816,0.997643638085699) circle (1pt);
\draw [fill=black] (0.9519671659638762,0) circle (1pt);
\draw [fill=black] (1.4634381800382552,0) circle (1pt);
\draw [fill=black] (1.7708857795953656,0.4898164476493253) circle (1pt);
\draw [fill=black] (1.4634381800382552,0.9910093011181378) circle (1pt);
\draw [fill=black] (2.076065280103183,0.9873917202164137) circle (1pt);
\draw [fill=black] (2.076065280103183,0) circle (1pt);
\draw [fill=black] (-0.4975752725670899,0.9933656630324388) circle (1pt);
\draw [fill=black] (-0.4975752725670899,0) circle (1pt);
\draw [fill=black] (2.5699067613713895,0.9891664816425414) circle (1pt);
\draw [fill=black] (2.5830178479797867,-0.00727610059566227) circle (1pt);

\draw[color=black] (0.5,0.5) node {$...$};
%hexn
\draw[color=black] (1,0.5) node {\tiny$x^{n}_{1}$};

\draw[color=black] (1,1.2) node {\tiny$x^{n}_{2}$};

\draw[color=black] (1.5,1.2) node {\tiny$x^{n}_{3}$};

\draw[color=black] (2.1,0.5) node {\tiny$x^{n}_{4}$};

\draw[color=black] (1,-0.2) node {\tiny$x^{n}_{6}$};

\draw[color=black] (1.51,-0.2) node {\tiny$x^{n}_{5}$};

\draw[color=black] (2.2,-0.2) node {\tiny$x^{n+1}_{6}$};

\draw[color=black] (2.3,1.2) node {\tiny$x^{n+1}_{2}$};

\draw[color=black] (3.1,0) node {\tiny$x^{n+1}_{5}$};

\draw[color=black] (3.1,1) node {\tiny$x^{n+1}_{3}$};

\end{scriptsize}
\end{tikzpicture}

 }%
\vskip -2 cm
\caption{The graph $\tilde{G}(n)$.}
\label{qtilde}
\end{figure}
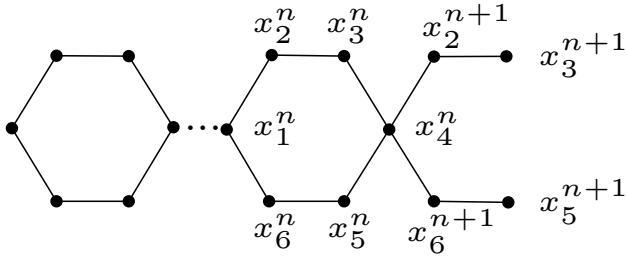

\noindent Next, we let $\overline{G}(n)$ be constructed from $G(n)$ by joining two vertices to a vertex at distance three from the cut vertex of the $n^{th}$ hexagon. Further, we let $\tilde{G}(n)$ be constructed from $\overline{G}(n)$ by joining two vertices to the end vertices of $\overline{G}(n)$, one vertex each. The graphs $\overline{G}(n)$ and $\tilde{G}(n)$ are shown in Figures \ref{qbar} and \ref{qtilde}, respectively.

\vskip -3 cm
\vskip 5 pt

\vskip -3 cm
\vskip 5 pt

\noindent Then, we let
\vskip 5 pt

\indent $\overline{g}(n, k)$ = the number of maximal independent sets containing $k$ vertices of $\overline{G}(n)$,
\vskip 5 pt

\indent $\tilde{g}(n, k)$ = the number of maximal independent sets containing $k$ vertices of $\tilde{G}(n)$,
\vskip 5 pt

\noindent and let
\begin{align}
\overline{G}(x, y) = \sum_{n \geq 0}\sum_{k \geq 0}\overline{g}(n, k)x^{n}y^{k}\notag
\end{align}
\begin{align}
\tilde{G}(x, y) = \sum_{n \geq 0}\sum_{k \geq 0}\tilde{g}(n, k)x^{n}y^{k}\notag
\end{align}
\noindent Now, we are ready to prove Theorem \ref{bivariate-ph}.
\vskip 5 pt

\noindent \emph{Proof of Theorem \ref{bivariate-ph}} Let $D$ be a maximal independent set of $G(n)$ containing $k$ vertices. We distinguish $2$ cases. 
\vskip 5 pt

\noindent \textbf{Case 1:} $x^{n}_{4} \in D$\\
\indent Thus, $x^{n}_{3}, x^{n}_{5} \notin D$. Removing the vertices $x^{n}_{3}, x^{n}_{4}, x^{n}_{5}$ from $G(n)$ result in $\overline{G}(n - 1)$. There are $\overline{g}(n - 1, k - 1)$ possibilities of $D$.
\vskip 5 pt

\noindent \textbf{Case 2:} $x^{n}_{4} \notin D$\\
\indent We further have the following $3$ subcases.
\vskip 5 pt

\noindent \textbf{Subcase 2.1:} $n^{n}_{3} \notin D$ and $x^{n}_{5} \in D$.\\
\indent Thus, $x^{n}_{6} \notin D$ and $x^{n}_{2} \in D$. This implies that $x^{n}_{1} \notin D$. Removing all the vertices of the $n^{th}$ hexagon results in $\tilde{G}(n - 2)$. We have that there are $\tilde{g}(n - 2, k - 2)$ possibilities of $D$.
\vskip 5 pt

\noindent \textbf{Subcase 2.2:} $n^{n}_{3} \in D$ and $x^{n}_{5} \notin D$.\\
\indent Thus, $x^{n}_{2} \notin D$ and $x^{n}_{6} \in D$. This implies that $x^{n}_{1} \notin D$. Removing all the vertices of the $n^{th}$ hexagon results in $\tilde{G}(n - 2)$. We have that there are $\tilde{g}(n - 2, k - 2)$ possibilities of $D$.
\vskip 5 pt

\noindent \textbf{Subcase 2.3:} $n^{n}_{3} \in D$ and $x^{n}_{5} \in D$.\\
\indent Thus, $x^{n}_{2}, x^{n}_{6} \notin D$. Removing the vertices $x^{n}_{2}, ..., x^{n}_{6}$ results in $G(n - 1)$. We have that there are $g(n - 1, k - 2)$ possibilities of $D$.
\vskip 5 pt
\vskip 5 pt

\indent From, the two cases, we have that
\begin{align}\label{ph11}
g(n, k) = \overline{g}(n - 1, k - 1) + 2\Tilde{g}(n - 2, k - 2) + g(n - 1, k - 2)
\end{align}
\noindent For $n \geq 2 $ and $k \geq 2 $, we multiply $x^{n}y^{k}$ throughout (\ref{ph11}) and sum over all $x^{n}y^{k}$. Thus, we have that
\begin{align}\label{ph22}
\sum_{n \geq 2}\sum_{k \geq 2}g(n, k)x^{n}y^{k} &= \sum_{n \geq 2}\sum_{k \geq 2}\overline{g}(n - 1, k - 1)x^{n}y^{k} + 2\sum_{n \geq 2}\sum_{k \geq 2}\Tilde{g}(n - 2, k - 2)x^{n}y^{k} \notag\\
&+ \sum_{n \geq 2}\sum_{k \geq 2}g(n - 1, k - 2)x^{n}y^{k}
\end{align}

\noindent We first consider the term $\displaystyle{\sum_{n \geq 2}\sum_{k \geq 2}g(n, k)x^{n}y^{k}}$. It can be checked that

\begin{center}
    $g(0, 0) = 1, g(0, 1) = 0, g(0, 2) = 0$, and $ g(0, k) = 0$ for all $k \geq 3$
    \vskip 5 pt
    $g(1, 0) = 0, g(1, 1) = 0, g(1, 2) = 3 ,g(1, 3) = 2$, and $ g(1, k) = 0$ for all $k \geq 4$
    \vskip 5 pt
\end{center}

\begin{align}\label{ph1}
\sum_{n \geq 2}\sum_{k \geq 2}g(n, k)x^{n}y^{k} &= \sum_{n \geq 2}\sum_{k \geq 2}g(n, k)x^{n}y^{k}  + g(0, 0)\notag\\
&+ g(1, 2)xy^{2} + g(1, 3)xy^{3} -g(0, 0) -g(1, 2)xy^{2} - g(1, 3)xy^{3}\notag\\
                  &= \sum_{n \geq 0}\sum_{k \geq 0}g(n, k)x^{n}y^{k} - 1 - 3xy^{2} - 2xy^{3} \notag\\
                  &= G(x,y) - 1 - 3xy^2 - 2xy^{3}
\end{align}

\noindent Now, we consider the term $\displaystyle{\sum_{n \geq 2}\sum_{k \geq 2}\overline{g}(n - 1, k - 1)x^{n}y^{k}}$. clearly,

\begin{center}
    $\overline{g}(0, 0) = 0, \overline{g}(0, 1) = 1 , \overline{g}(0, 2) = 1$ and $\overline{g}(0, k) = 0$ for all $k \geq 2$
    \vskip 5 pt
\end{center}

\begin{align}\label{phbar1}
\sum_{n \geq 2}\sum_{k \geq 2}\overline{g}(n - 1, k - 1)x^{n}y^{k}
    &= xy(\sum_{n \geq 2}\sum_{k \geq 2}\overline{g}(n - 1, k - 1)x^{n - 1}y^{k - 1})\notag\\
    &= xy(\sum_{n \geq 2}\sum_{k \geq 2}\overline{g}(n - 1, k - 1)x^{n - 1}y^{k - 1} + \overline{g}(0,1)y + \overline{g}(0,2)y^{2}\notag\\
    &- \overline{g}(0,1)y - \overline{g}(0,2)y^{2})\notag\\
    &= xy(\sum_{n \geq 0}\sum_{k \geq 0}\overline{g}(n, k)x^{n}y^{k} - y - y^{2})\notag\\
    &= xy\overline{G}(x, y) - xy^{2} - xy^{3}.
\end{align}

\noindent Next, we consider the term $\displaystyle{\sum_{n \geq 2}\sum_{k \geq 2}2\Tilde{g}(n - 2, k - 2)x^{n}y^{k}}$. It can be check that,

\begin{align}\label{phtil1}
\sum_{n \geq 2}\sum_{k \geq 2}2\Tilde{g}(n - 2, k - 2)x^{n}y^{k}
    &= 2x^{2}y^{2}(\sum_{n \geq 2}\sum_{k \geq 2}\Tilde{g}(n - 2, k - 2)x^{n - 2}y^{k - 2})\notag\\
    &= 2x^{2}y^{2}(\sum_{n \geq 0}\sum_{k \geq 0}\Tilde{g}(n , k)x^{n}y^{k}\notag\\
    &= 2x^{2}y^{2}\Tilde{G}(x, y).
\end{align}

\noindent Finally, we consider the term
$\displaystyle{\sum_{n \geq 2}\sum_{k \geq 2}g(n - 1 , k - 2)x^{n}y^{k}}$. It can be checked that, 

\begin{align}\label{ph2}
\sum_{n \geq 2}\sum_{k \geq 2}g(n - 1, k - 2)x^{n}y^{k}
    &= xy^{2}(\sum_{n \geq 2}\sum_{k \geq 2}g(n - 1, k - 2)x^{n - 1}y^{k - 2})\notag\\
    &= xy^{2}(\sum_{n \geq 2}\sum_{k \geq 2}g(n - 1, k - 2)x^{n - 1}y^{k - 2} + g(0,0)-g(0,0))\notag\\
    &= xy^{2}(\sum_{n \geq 0}\sum_{k \geq 0}g(n, k)x^{n}y^{k} - 1)\notag\\
    &= xy^{2}G(x, y) - xy^{2}
\end{align}

\noindent Plugging (\ref{ph1}), (\ref{phbar1}) , (\ref{phtil1}) and (\ref{ph2}) to (\ref{ph22}), we have 
\begin{align}
G(x,y) - 1- 3xy^{2} - 2xy^{3} &= xy\overline{G}(x, y) - xy^{2} - xy^{3} + 2x^{2}y^{2}\Tilde{G}(x, y) + xy^{2}G(x, y) - xy^{2} \notag
\end{align}

\noindent which can be solved that

\begin{align}\label{semifinalph4}
G(x,y) &= 1 + 3xy^2 + 2xy^{3} + xy\overline{G}(x, y) - xy^{2} - xy^{3} + 2x^{2}y^{2}\Tilde{G}(x, y) + xy^{2}G(x, y) - xy^{2}.
\end{align}
\vskip 5 pt

\indent Next, we will find the recurrence relation of $\overline{G}(n)$. Let $D$ be a maximal independent set of $\overline{G}(n)$ containing $k$ vertices. We distinguish $2$ cases. 
\vskip 5 pt

\noindent \textbf{Case 1:} $x^{n}_{4} \in D$\\
\indent Thus, $x^{n}_{3}, x^{n}_{5}, x^{n+1}_{2}, x^{n+1}_{6} \notin D$. Removing $x^{n}_{3}, x^{n}_{4}, x^{n}_{5}, x^{n+1}_{2}, x^{n+1}_{6}$ from $\overline{G}(n)$ results in $\overline{G}(n - 1)$. There are $\overline{g}(n - 1, k - 1)$ possibilities of $D$.
\vskip 5 pt

\noindent \textbf{Case 2:} $x^{n}_{4} \notin D$\\
\indent By the maximality of $D$, $x^{n+1}_{2}, x^{n+1}_{6} \in D$. Removing $x^{n}_{4}, x^{n+1}_{2}, x^{n+1}_{6}$ from $\overline{G}(n)$ results in $\tilde{G}(n - 1)$. There are $\tilde{g}(n - 1, k - 2)$ possibilities of $D$.
\vskip 5 pt

\indent From the two cases, we have that
\begin{align}\label{phbar2}
\overline{g}(n, k) = \overline{g}(n - 1, k - 1) + \Tilde{g}(n - 1, k - 2)
\end{align}
\vskip 5 pt

\noindent For $n \geq 1$ and $k \geq 2$, we multiply $x^{n}y^{k}$ throughout (\ref{phbar2}) and sum over all $x^{n}y^{k}$. Thus, we have that

\begin{align}\label{phbar22222}
\sum_{n \geq 1}\sum_{k \geq 2}\overline{g}(n, k)x^{n}y^{k} = \sum_{n \geq 1}\sum_{k \geq 2}\overline{g}(n - 1, k - 1)x^{n}y^{k} + \sum_{n \geq 1}\sum_{k \geq 2}\Tilde{g}(n - 1, k - 2)x^{n}y^{k}
\end{align}

\noindent We first consider the term $\displaystyle{\sum_{n \geq 1}\sum_{k \geq 2}\overline{g}(n , k )x^{n}y^{k}}$. It can be checked that, 

\begin{center}
    $\overline{g}(0, 0) = 0, \overline{g}(0, 1) = 1, \overline{g}(0, 2) =1$ and $\overline{g}(0, k) = 0$ for all $k \geq 3$
    \vskip 5 pt
\end{center}

\begin{align}\label{phbar3}
\sum_{n \geq 1}\sum_{k \geq 2}\overline{g}(n, k)x^{n}y^{k}
    &= \sum_{n \geq 1}\sum_{k \geq 2}\overline{g}(n,k)x^{n}y^{k} + \overline{g}(0, 1)y +\overline{g}(0, 2)y^{2}
    -\overline{g}(0, 1)y - \overline{g}(0, 2)y^{2}\notag\\
    &= \sum_{n \geq 0}\sum_{k \geq 0}\overline{g}(n , k )x^{n}y^{k} - y - y^{2}\notag\\
    &= \overline{G}(x, y) - y - y^{2}
\end{align}

\noindent We next consider the term $\displaystyle{\sum_{n \geq 1}\sum_{k \geq 2}\overline{g}(n - 1 , k - 1 )x^{n}y^{k}}$, clearly

\begin{align}\label{phbar4}
\sum_{n \geq 1}\sum_{k \geq 2}\overline{g}(n-1, k-1)x^{n}y^{k} &= xy(\sum_{n \geq 1}\sum_{k \geq 2}g(n-1, k-1)x^{n-1}y^{k-1}\notag\\
                  &= xy(\sum_{n \geq 0}\sum_{k \geq 1}\overline{g}(n , k)x^{n}y^{k})\notag\\
                  &= xy(\sum_{n \geq 0}\sum_{k \geq 0}\overline{g}(n , k)x^{n}y^{k})\notag\\
                  &= xy\overline{G}(x,y) 
\end{align}

\noindent Finally, we consider the term $\displaystyle{\sum_{n \geq 1}\sum_{k \geq 2}\Tilde{g}(n - 1 , k - 2 )x^{n}y^{k}}$. It can be checked that,

\begin{align}\label{phtil2}
\sum_{n \geq 1}\sum_{k \geq 2}\Tilde{g}(n - 1 , k - 2)x^{n}y^{k}
    &= xy^{2}(\sum_{n \geq 1}\sum_{k \geq 2}\Tilde{g}(n - 1, k - 2 )x^{n - 1}y^{k - 2})\notag\\
    &= xy^{2}\Tilde{G}(x, y)
\end{align}    
\vskip 5 pt

\noindent Plugging (\ref{phbar3}) , (\ref{phbar4}) and (\ref{phtil2}) to (\ref{phbar22222}), we have 

\begin{align}
\overline{G}(x, y) - y - y^{2} &= xy\overline{G}(x,y) + xy^{2}\Tilde{G}(x, y).
\end{align}

\noindent which can be solved that

\begin{align}\label{semifinalphbar4}
\overline{G}(x, y) &= y + y^{2} + xy\overline{G}(x,y) + xy^{2}\Tilde{G}(x, y).
\end{align}
\vskip 5 pt

\indent Next, we will find the recurrence relation of $\tilde{G}$. Let $D$ be a maximal independent set of $\tilde{G}(n)$ containing $k$ vertices. We distinguish $4$ cases. 
\vskip 5 pt

\noindent \textbf{Case 1:} $x^{n+1}_{3} \in D$ and $x^{n+1}_{5} \in D$\\
\indent Thus, $x^{n+1}_{2}, x^{n+1}_{6} \notin D$. Removing $x^{n+1}_{2}, x^{n+1}_{3}, x^{n+1}_{5}, x^{n+1}_{6}$ from $\tilde{G}(n)$ results in $G(n)$. There are $g(n, k - 2)$ possibilities of $D$.
\vskip 5 pt

\noindent \textbf{Case 2:} $x^{n+1}_{3} \notin D$ and $x^{n+1}_{5} \in D$\\
\indent Thus, $x^{n+1}_{2} \in D$ but $x^{n+1}_{6} \notin D$. Since $D$ is independent, $x^{n+1}_{1} \notin D$. Removing $x^{n+1}_{1}, x^{n+1}_{2},$ $x^{n+1}_{3}, x^{n+1}_{5}, x^{n+1}_{6}$ from $\tilde{G}(n)$ results in $\tilde{G}(n-1)$. There are $\tilde{g}(n-1, k - 2)$ possibilities of $D$.
\vskip 5 pt

\noindent \textbf{Case 3:} $x^{n+1}_{3} \in D$ and $x^{n+1}_{5} \notin D$\\
\indent By similar arguments in Case 2, there are $\tilde{g}(n-1, k - 2)$ possibilities of $D$.
\vskip 5 pt

\noindent \textbf{Case 4:} $x^{n+1}_{3} \notin D$ and $x^{n+1}_{5} \notin D$\\
\indent Thus, $x^{n+1}_{2}, x^{n+1}_{6} \in D$. Since $D$ is independent, $x^{n+1}_{1} \notin D$. Removing $x^{n+1}_{1}, x^{n+1}_{2}, x^{n+1}_{3}, x^{n+1}_{5},$ $x^{n+1}_{6}$ from $\tilde{G}(n)$ results in $\tilde{G}(n-1)$. There are $\tilde{g}(n-1, k - 2)$ possibilities of $D$.
\vskip 5 pt

\indent From the three cases, we have that
\begin{align}\label{phtil4}
    \Tilde{g}(n, k) = g(n ,k - 2) + 3\Tilde{g}(n - 1, k - 2).
\end{align}
\vskip 5 pt

\noindent For $n \geq 1$ and $k \geq 2$, we multiply $x^{n}y^{k}$ throughout (\ref{phtil4}) and sum over all $x^{n}y^{k}$. Thus, we have that

\begin{align}\label{phtil44444}
\sum_{n \geq 1}\sum_{k \geq 2}\Tilde{g}(n, k)x^{n}y^{k} = \sum_{n \geq 1}\sum_{k \geq 2}g(n ,k - 2)x^{n}y^{k} + \sum_{n \geq 1}\sum_{k \geq 2}3\Tilde{g}(n - 1, k - 2)x^{n}y^{k}.
\end{align}

\noindent We first consider the term $\displaystyle{\sum_{n \geq 1}\sum_{k \geq 2}\Tilde{g}(n , k )x^{n}y^{k}}$. It can be checked that, 

\begin{center}
    $\Tilde{g}(0, 0) = 0, \Tilde{g}(0, 1) = 0 , \Tilde{g}(0, 2) = 3, \Tilde{g}(0,3) =1$ and $\Tilde{g}(0, k) = 0$ for all $k \geq 4$
    \vskip 5 pt
\end{center}

\begin{align}\label{phtil5}
\sum_{n \geq 1}\sum_{k \geq 2}\Tilde{g}(n , k)x^{n}y^{k}
    &= \sum_{n \geq 1}\sum_{k \geq 2}\Tilde{g}(n , k)x^{n}y^{k}
    + \Tilde{g}(0,2)y^{2} + \Tilde{g}(0,3)y^{3}\notag\\ &-\Tilde{g}(0,2)y^{2} - \Tilde{g}(0,3)y^{3}\notag\\
    &= \sum_{n \geq 0}\sum_{k \geq 0}\Tilde{g}(n , k )x^{n}y^{k} - 3y^{2} - y^{3}\notag\\
    &= \Tilde{G}(x, y) - 3y^{2} - y^{3}
\end{align}    
\vskip 5 pt

\noindent We next consider the term $\displaystyle{\sum_{n \geq 1}\sum_{k \geq 2}g(n, k - 2)x^{n}y^{k}}$, clearly

\begin{align}\label{ph3}
\sum_{n \geq 1}\sum_{k \geq 2}g(n, k-2)x^{n}y^{k} &= y^{2}(\sum_{n \geq 1}\sum_{k \geq 2}g(n, k-2)x^{n}y^{k-2} + g(0, 0) -g(0,0))\notag\\
                  &= y^{2}(\sum_{n \geq 0}\sum_{k \geq 0}g(n, k)x^{n}y^{k} - 1)\notag\\
                  &= y^{2}G(x,y) - y^{2}
\end{align}

\noindent Finally, we consider the term
$\displaystyle{\sum_{n \geq 1}\sum_{k \geq 2}3\Tilde{g}(n - 1 , k - 2)x^{n}y^{k}}$. It can be checked that, 

\begin{align}\label{phtil6}
\sum_{n \geq 1}\sum_{k \geq 2}3\Tilde{g}(n - 1 , k - 2)x^{n}y^{k}
    &= 3xy^{2}\sum_{n \geq 1}\sum_{k \geq 2}\Tilde{g}(n - 1 , k - 2 )x^{n - 1}y^{k - 2}\notag\\
    &= 3xy^{2}(\sum_{n \geq 0}\sum_{k \geq 0}\Tilde{g}(n , k )x^{n}y^{k})\notag\\
    &= 3xy^{2}\Tilde{G}(x, y) 
\end{align}

\noindent Plugging (\ref{phtil5}) , (\ref{ph3}) and (\ref{phtil6}) to (\ref{phtil44444}), we have

\begin{align}
\Tilde{G}(x, y) - 3y^{2} - y^{3} &= y^{2}G(x,y) - y^{2} + 3xy^{2}\Tilde{G}(x, y) 
\end{align}

\noindent which can be solved that

\begin{align}\label{semifinalphtil7}
\Tilde{G}(x, y) = 3y^{2} + y^{3} + y^{2}G(x,y) - y^{2} + 3xy^{2}\Tilde{G}(x, y) 
\end{align}
\vskip 5 pt

\noindent By plugging (\ref{semifinalph4}) to (\ref{semifinalphbar4}) to (\ref{semifinalphtil7}) , we have
\begin{align}
G(x, y) = \frac{x^{3}y^{6}- x^{3}y^{5}+ x^{2}y^{6}-2x^{2}y^{5} - 3x^{2}y^{4}   +2x^2y^3+2xy^3 -xy^{2} -xy + 1}{1 - xy - 4xy^{2} +4x^{2}y^{3} - x^{2}y^{5} + x^{2}y^{4} - x^{3}y^{5} }.\notag
\end{align}

\noindent as required.
\qed

\noindent We will prove Theorem \ref{recurrence-ph}.
\vskip 5 pt

\noindent \emph{Proof of Theorem \ref{recurrence-ph}}. By Theorem \ref{bivariate-ph} with $y = 1$, we have that
\vskip 5 pt

\begin{align}
\sum_{n \geq 0}g(n)x^{n} &= G(x, 1)\notag\\
                         &= \frac{1 -2x^{2}}{1 - 5x + 4x^{2} - x^{3}}
\end{align}
 
 \noindent which can  be solved that

\begin{align}\label{recursived-ph}
1 -2x^{2}   &= (1 - 5x + 4x^{2} - x^{3})\sum_{n \geq 0}g(n)x^{n}\notag\\
                  &= -\sum_{n \geq 0}g(n)x^{n + 3} + 4\sum_{n \geq 0}g(n)x^{n + 2} - 5\sum_{n \geq 0}g(n)x^{n + 1} + \sum_{n \geq 0}g(n)x^{n} \notag\\
                  &= -\sum_{n \geq 3}g(n-3)x^{n} + 4\sum_{n \geq 2}g(n - 2)x^{n} - 5\sum_{n \geq 1}g(n - 1)x^{n} + \sum_{n \geq 0}g(n)x^{n} \notag\\
                  &= g(0) + g(1)x + g(2)x^{2} - 5g(0)x - 5g(1)x^{2} + 4g(0)x^{2}\notag\\
                  &+ \sum_{n \geq 3}(-g(n-3) + 4g(n-2) - 5g(n-1) + g(n))x^{n}
\end{align}
\noindent Because the order of the polynomial on the left hand side of (\ref{recursived-ph}) is two, the coefficients of $x^{n}$ for all $n \geq 3$ must be $0$. Thus, $g(n) - 5g(n - 1) + 4g(n - 2) - g(n - 3) = 0$ implying that

\begin{align*}
    g(n) = 5g(n - 1) - 4g(n - 2) + g(n - 3)
\end{align*}
\noindent and this proves Theorem \ref{recurrence-ph}.
\qed

\subsection{Ortho-Hexagonal Cacti}
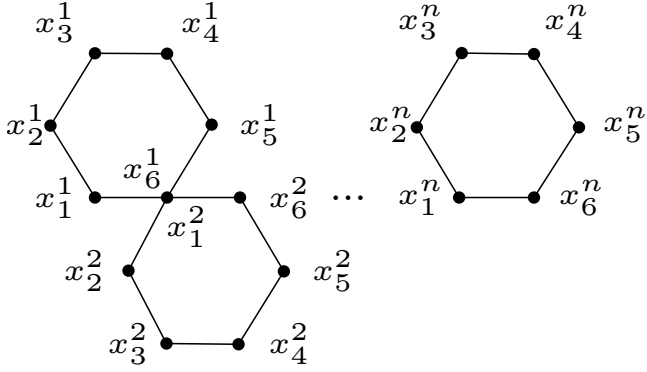
\begin{figure}[H]
\centering
\definecolor{ududff}{rgb}{0.30196078431372547,0.30196078431372547,1}
\resizebox{0.8\textwidth}{!}{%

\begin{tikzpicture}[line cap=round,line join=round,>=triangle 45,x=1cm,y=1cm]
\clip(-0.32537133026435466,-1.5238107651131987) rectangle (5.916913749858403,2.33134280175197);
\draw [line width=0.3pt] (0.7636911477878787,0.983396738469664)-- (1.068196050713675,0.4982589804623496);
\draw [line width=0.3pt] (1.068196050713675,0.4982589804623496)-- (0.7636911477878787,0);
\draw [line width=0.3pt] (2.469888418299665,0.4797553368451276)-- (2.7770094292928365,0.9875825272815013);
\draw [line width=0.3pt] (2.7770094292928365,0.9875825272815013)-- (3.270375092325276,0.9809481903139401);
\draw [line width=0.3pt] (3.270375092325276,0.9809481903139401)-- (3.5778226918823863,0.4797553368451276);
\draw [line width=0.3pt] (3.5778226918823863,0.4797553368451276)-- (3.270375092325276,0);
\draw [line width=0.3pt] (2.469888418299665,0.4797553368451276)-- (2.7589040782508967,0);
\draw [line width=0.3pt] (2.7589040782508967,0)-- (3.270375092325276,0);
\draw [line width=0.3pt] (-0.03310388590148339,0.4916246434947884)-- (0.2709525523276073,0.9861831599960398);
\draw [line width=0.3pt] (-0.03310388590148339,0.4916246434947884)-- (0.2709525523276073,0);
\draw [line width=0.3pt] (0.2709525523276073,0.9861831599960398)-- (0.7636911477878787,0.983396738469664);
\draw [line width=0.3pt] (0.2709525523276073,0)-- (0.7636911477878787,0);
\draw [line width=0.3pt] (0.5,-0.5)-- (0.7636911477878787,0);
\draw [line width=0.3pt] (1.2619236772230777,0)-- (1.561584020328892,-0.5050167892307563);
\draw [line width=0.3pt] (0.5,-0.5)-- (0.7594444433566199,-0.9983600370269141);
\draw [line width=0.3pt] (1.2518740925457486,-1.0033848293655787)-- (1.561584020328892,-0.5050167892307563);
\draw [line width=0.3pt] (0.7636911477878787,0)-- (1.2619236772230777,0);
\draw [line width=0.3pt] (0.7594444433566199,-0.9983600370269141)-- (1.2518740925457486,-1.0033848293655787);
\begin{scriptsize}
\draw [fill=black] (-0.03310388590148339,0.4916246434947884) circle (1pt);
\draw [fill=black] (0.7636911477878787,0.983396738469664) circle (1pt);
\draw [fill=black] (0.7636911477878787,0) circle (1pt);
\draw [fill=black] (1.068196050713675,0.4982589804623496) circle (1pt);
\draw [fill=black] (2.469888418299665,0.4797553368451276) circle (1pt);
\draw [fill=black] (2.7770094292928365,0.9875825272815013) circle (1pt);
\draw [fill=black] (2.7589040782508967,0) circle (1pt);
\draw [fill=black] (3.270375092325276,0) circle (1pt);
\draw [fill=black] (3.5778226918823863,0.4797553368451276) circle (1pt);
\draw [fill=black] (3.270375092325276,0.9809481903139401) circle (1pt);
\draw [fill=black] (0.2709525523276073,0.9861831599960398) circle (1pt);
\draw [fill=black] (0.2709525523276073,0) circle (1pt);
\draw [fill=black] (1.2619236772230777,0) circle (1pt);
\draw [fill=black] (0.7594444433566199,-0.9983600370269141) circle (1pt);
\draw [fill=black] (1.2518740925457486,-1.0033848293655787) circle (1pt);
\draw [fill=black] (1.561584020328892,-0.5050167892307563) circle (1pt);
\draw [fill=black] (0.5,-0.5) circle (1pt);

\draw[color=black] (2,0) node {$...$};

%%%%hex1
\draw[color=black] (0,1.2) node {\tiny$x^{1}_{3}$};

\draw[color=black] (-0.2,0.5) node {\tiny$x^{1}_{2}$};

\draw[color=black] (0,0) node {\tiny$x^{1}_{1}$};

\draw[color=black] (0.6,0.2) node {\tiny$x^{1}_{6}$};

\draw[color=black] (1,1.2) node {\tiny$x^{1}_{4}$};

\draw[color=black] (1.4,0.5) node {\tiny$x^{1}_{5}$};

%%%hex2
\draw[color=black] (1.6,0) node {\tiny$x^{2}_{6}$};

\draw[color=black] (1.9,-0.5) node {\tiny$x^{2}_{5}$};

\draw[color=black] (0.2,-0.5) node {\tiny$x^{2}_{2}$};

\draw[color=black] (0.9,-0.2) node {\tiny$x^{2}_{1}$};

\draw[color=black] (1.6,-1) node {\tiny$x^{2}_{4}$};

\draw[color=black] (0.5,-1) node {\tiny$x^{2}_{3}$};

%%hexn

\draw[color=black] (2.5,1.2) node {\tiny$x^{n}_{3}$};

\draw[color=black] (2.3,0.5) node {\tiny$x^{n}_{2}$};

\draw[color=black] (2.5,0) node {\tiny$x^{n}_{1}$};

\draw[color=black] (3.6,0) node {\tiny$x^{n}_{6}$};

\draw[color=black] (3.5,1.2) node {\tiny$x^{n}_{4}$};

\draw[color=black] (3.9,0.5) node {\tiny$x^{n}_{5}$};

\end{scriptsize}
\end{tikzpicture}

 }%
%\vskip -1 cm
\caption{The graph $Q(n)$ with labelled vertices.}
\label{ohofnh}
\end{figure}

First, we may name all the vertices of the ortho-hexagonal cactus of $n$ hexagons as shown by Figure \ref{ohofnh}.

\vskip -1 cm
\vskip 5 pt

\noindent Then, we recall that
\vskip 5 pt

\indent $q(n)$ = the number of all maximal independent sets of $Q(n)$

\noindent and

\indent $q(n, k)$ = the number of maximal independent sets containing $k$ vertices of $Q(n)$.
\vskip 5 pt

\noindent Thus,
\begin{align}\label{Pn}
Q(n) = \sum_{k \geq 0}q(n, k).    
\end{align}
\vskip 5 pt

\noindent Further, we let 
\begin{align}
Q(x) = \sum_{n \geq 0}q(n)x^{n}\notag
\end{align}
\noindent be the generating function of $q(n)$ and we let
\begin{align}
Q(x, y) = \sum_{n \geq 0}\sum_{k \geq 0}q(n, k)x^{n}y^{k}\notag
\end{align}
\noindent be the bi-variate generating function of $q(n, k)$. It is worth noting that, when $y = 1$, we have
\begin{align}\label{pxytopn}
Q(x, 1) = \sum_{n \geq 0}(\sum_{k \geq 0}q(n, k)(1)^{k})x^{n} = \sum_{n \geq 0}q(n)x^{n} = Q(x).
\end{align}

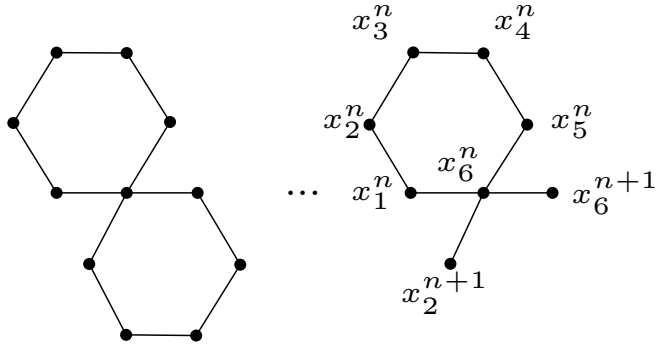
\begin{figure}[H]
\centering
\definecolor{ududff}{rgb}{0.30196078431372547,0.30196078431372547,1}
\resizebox{0.8\textwidth}{!}{%

\begin{tikzpicture}[line cap=round,line join=round,>=triangle 45,x=1cm,y=1cm]
\clip(-0.39058887004316967,-1.598224723822629) rectangle (6.033226076009011,2.3690393123895466);
\draw [line width=0.3pt] (0.7636911477878787,0.983396738469664)-- (1.068196050713675,0.4982589804623496);
\draw [line width=0.3pt] (1.068196050713675,0.4982589804623496)-- (0.7636911477878787,0);
\draw [line width=0.3pt] (2.469888418299665,0.4797553368451276)-- (2.7770094292928365,0.9875825272815013);
\draw [line width=0.3pt] (2.7770094292928365,0.9875825272815013)-- (3.270375092325276,0.9809481903139401);
\draw [line width=0.3pt] (3.270375092325276,0.9809481903139401)-- (3.5778226918823863,0.4797553368451276);
\draw [line width=0.3pt] (3.5778226918823863,0.4797553368451276)-- (3.270375092325276,0);
\draw [line width=0.3pt] (2.469888418299665,0.4797553368451276)-- (2.7589040782508967,0);
\draw [line width=0.3pt] (2.7589040782508967,0)-- (3.270375092325276,0);
\draw [line width=0.3pt] (-0.03310388590148339,0.4916246434947884)-- (0.2709525523276073,0.9861831599960398);
\draw [line width=0.3pt] (-0.03310388590148339,0.4916246434947884)-- (0.2709525523276073,0);
\draw [line width=0.3pt] (0.2709525523276073,0.9861831599960398)-- (0.7636911477878787,0.983396738469664);
\draw [line width=0.3pt] (0.2709525523276073,0)-- (0.7636911477878787,0);
\draw [line width=0.3pt] (0.5,-0.5)-- (0.7636911477878787,0);
\draw [line width=0.3pt] (1.2619236772230777,0)-- (1.561584020328892,-0.5050167892307563);
\draw [line width=0.3pt] (0.5,-0.5)-- (0.7594444433566199,-0.9983600370269141);
\draw [line width=0.3pt] (1.2518740925457486,-1.0033848293655787)-- (1.561584020328892,-0.5050167892307563);
\draw [line width=0.3pt] (0.7636911477878787,0)-- (1.2619236772230777,0);
\draw [line width=0.3pt] (0.7594444433566199,-0.9983600370269141)-- (1.2518740925457486,-1.0033848293655787);
\draw [line width=0.3pt] (3.270375092325276,0)-- (3.755591075111256,0);
\draw [line width=0.3pt] (3.270375092325276,0)-- (3.0383080049123388,-0.4986518900368979);
\begin{scriptsize}
\draw [fill=black] (-0.03310388590148339,0.4916246434947884) circle (1pt);
\draw [fill=black] (0.7636911477878787,0.983396738469664) circle (1pt);
\draw [fill=black] (0.7636911477878787,0) circle (1pt);
\draw [fill=black] (1.068196050713675,0.4982589804623496) circle (1pt);
\draw [fill=black] (2.469888418299665,0.4797553368451276) circle (1pt);
\draw [fill=black] (2.7770094292928365,0.9875825272815013) circle (1pt);
\draw [fill=black] (2.7589040782508967,0) circle (1pt);
\draw [fill=black] (3.270375092325276,0) circle (1pt);
\draw [fill=black] (3.5778226918823863,0.4797553368451276) circle (1pt);
\draw [fill=black] (3.270375092325276,0.9809481903139401) circle (1pt);
\draw [fill=black] (0.2709525523276073,0.9861831599960398) circle (1pt);
\draw [fill=black] (0.2709525523276073,0) circle (1pt);
\draw [fill=black] (1.2619236772230777,0) circle (1pt);
\draw [fill=black] (0.7594444433566199,-0.9983600370269141) circle (1pt);
\draw [fill=black] (1.2518740925457486,-1.0033848293655787) circle (1pt);
\draw [fill=black] (1.561584020328892,-0.5050167892307563) circle (1pt);
\draw [fill=black] (0.5,-0.5) circle (1pt);
\draw [fill=black] (3.755591075111256,0) circle (1pt);
\draw [fill=black] (3.0383080049123388,-0.4986518900368979) circle (1pt);

\draw[color=black] (2,0) node {$...$};

%%hexn

\draw[color=black] (2.5,1.2) node {\tiny$x^{n}_{3}$};

\draw[color=black] (2.3,0.5) node {\tiny$x^{n}_{2}$};

\draw[color=black] (2.5,0) node {\tiny$x^{n}_{1}$};

\draw[color=black] (3.1,0.2) node {\tiny$x^{n}_{6}$};

\draw[color=black] (3.5,1.2) node {\tiny$x^{n}_{4}$};

\draw[color=black] (3.9,0.5) node {\tiny$x^{n}_{5}$};

\draw[color=black] (3,-0.7) node {\tiny$x^{n+1}_{2}$};

\draw[color=black] (4.2,0) node {\tiny$x^{n+1}_{6}$};
\end{scriptsize}
\end{tikzpicture}

 }%
%\vskip -1 cm
\caption{The graph $\overline{Q}(n)$}
\label{ohbar}
\end{figure}

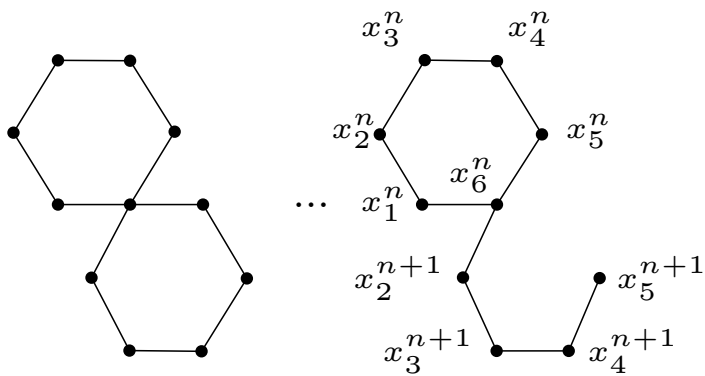
\begin{figure}[H]
\centering
\definecolor{ududff}{rgb}{0.30196078431372547,0.30196078431372547,1}
\resizebox{0.8\textwidth}{!}{%

\begin{tikzpicture}[line cap=round,line join=round,>=triangle 45,x=1cm,y=1cm]
\clip(-0.32537133026435466,-1.5238107651131987) rectangle (5.916913749858403,2.33134280175197);
\draw [line width=0.3pt] (0.7636911477878787,0.983396738469664)-- (1.068196050713675,0.4982589804623496);
\draw [line width=0.3pt] (1.068196050713675,0.4982589804623496)-- (0.7636911477878787,0);
\draw [line width=0.3pt] (2.469888418299665,0.4797553368451276)-- (2.7770094292928365,0.9875825272815013);
\draw [line width=0.3pt] (2.7770094292928365,0.9875825272815013)-- (3.270375092325276,0.9809481903139401);
\draw [line width=0.3pt] (3.270375092325276,0.9809481903139401)-- (3.5778226918823863,0.4797553368451276);
\draw [line width=0.3pt] (3.5778226918823863,0.4797553368451276)-- (3.270375092325276,0);
\draw [line width=0.3pt] (2.469888418299665,0.4797553368451276)-- (2.7589040782508967,0);
\draw [line width=0.3pt] (2.7589040782508967,0)-- (3.270375092325276,0);
\draw [line width=0.3pt] (-0.03310388590148339,0.4916246434947884)-- (0.2709525523276073,0.9861831599960398);
\draw [line width=0.3pt] (-0.03310388590148339,0.4916246434947884)-- (0.2709525523276073,0);
\draw [line width=0.3pt] (0.2709525523276073,0.9861831599960398)-- (0.7636911477878787,0.983396738469664);
\draw [line width=0.3pt] (0.2709525523276073,0)-- (0.7636911477878787,0);
\draw [line width=0.3pt] (0.5,-0.5)-- (0.7636911477878787,0);
\draw [line width=0.3pt] (1.2619236772230777,0)-- (1.561584020328892,-0.5050167892307563);
\draw [line width=0.3pt] (0.5,-0.5)-- (0.7594444433566199,-0.9983600370269141);
\draw [line width=0.3pt] (1.2518740925457486,-1.0033848293655787)-- (1.561584020328892,-0.5050167892307563);
\draw [line width=0.3pt] (0.7636911477878787,0)-- (1.2619236772230777,0);
\draw [line width=0.3pt] (0.7594444433566199,-0.9983600370269141)-- (1.2518740925457486,-1.0033848293655787);
\draw [line width=0.3pt] (3.270375092325276,0)-- (3.0383080049123388,-0.4986518900368979);
\draw [line width=0.3pt] (3.0383080049123388,-0.4986518900368979)-- (3.267964222953829,-1.0005077912145033);
\draw [line width=0.3pt] (3.267964222953829,-1.0005077912145033)-- (3.7605880685099775,-1.0005077912145033);
\draw [line width=0.3pt] (3.7605880685099775,-1.0005077912145033)-- (3.9718356882360695,-0.5028131662957285);
\begin{scriptsize}
\draw [fill=black] (-0.03310388590148339,0.4916246434947884) circle (1pt);
\draw [fill=black] (0.7636911477878787,0.983396738469664) circle (1pt);
\draw [fill=black] (0.7636911477878787,0) circle (1pt);
\draw [fill=black] (1.068196050713675,0.4982589804623496) circle (1pt);
\draw [fill=black] (2.469888418299665,0.4797553368451276) circle (1pt);
\draw [fill=black] (2.7770094292928365,0.9875825272815013) circle (1pt);
\draw [fill=black] (2.7589040782508967,0) circle (1pt);
\draw [fill=black] (3.270375092325276,0) circle (1pt);
\draw [fill=black] (3.5778226918823863,0.4797553368451276) circle (1pt);
\draw [fill=black] (3.270375092325276,0.9809481903139401) circle (1pt);
\draw [fill=black] (0.2709525523276073,0.9861831599960398) circle (1pt);
\draw [fill=black] (0.2709525523276073,0) circle (1pt);
\draw [fill=black] (1.2619236772230777,0) circle (1pt);
\draw [fill=black] (0.7594444433566199,-0.9983600370269141) circle (1pt);
\draw [fill=black] (1.2518740925457486,-1.0033848293655787) circle (1pt);
\draw [fill=black] (1.561584020328892,-0.5050167892307563) circle (1pt);
\draw [fill=black] (0.5,-0.5) circle (1pt);
\draw [fill=black] (3.0383080049123388,-0.4986518900368979) circle (1pt);
\draw [fill=black] (3.267964222953829,-1.0005077912145033) circle (1pt);
\draw [fill=black] (3.7605880685099775,-1.0005077912145033) circle (1pt);
\draw [fill=black] (3.9718356882360695,-0.5028131662957285) circle (1pt);

\draw[color=black] (2,0) node {$...$};

%%hexn

\draw[color=black] (2.5,1.2) node {\tiny$x^{n}_{3}$};

\draw[color=black] (2.3,0.5) node {\tiny$x^{n}_{2}$};

\draw[color=black] (2.5,0) node {\tiny$x^{n}_{1}$};

\draw[color=black] (3.1,0.2) node {\tiny$x^{n}_{6}$};

\draw[color=black] (3.5,1.2) node {\tiny$x^{n}_{4}$};

\draw[color=black] (3.9,0.5) node {\tiny$x^{n}_{5}$};

\draw[color=black] (2.6,-0.5) node {\tiny$x^{n+1}_{2}$};

\draw[color=black] (4.4,-0.5) node {\tiny$x^{n+1}_{5}$};

\draw[color=black] (2.8,-1) node {\tiny$x^{n+1}_{3}$};

\draw[color=black] (4.2,-1) node {\tiny$x^{n+1}_{4}$};
\end{scriptsize}
\end{tikzpicture}

 }%
%\vskip -1 cm
\caption{The graph $\tilde{Q}(n)$.}
\label{ohtilde}
\end{figure}

\noindent Next, we let $\overline{Q}(n)$ be constructed from $Q(n)$ by joining two vertices to a vertex at distance one from the cut vertex of the $n^{th}$ hexagons. Further, we let $\tilde{Q}(n)$ be constructed from $Q(n)$ by joining an end vertex of a path of length three to a vertex at distance one from the cut vertex of the $n^{th}$ hexagons. The graphs $\overline{Q}(n)$ and $\tilde{Q}(n)$ are shown in Figures \ref{ohbar} and \ref{ohtilde}, respectively.
\vskip 5 pt

\vskip -1 cm
\vskip 5 pt

\vskip -1 cm
\vskip 5 pt

\noindent Then, we let
\vskip 5 pt

\indent $\overline{Q}(n, k)$ = the number of maximal independent sets containing $k$ vertices of $\overline{Q}(n)$,
\vskip 5 pt

\indent $\tilde{Q}(n, k)$ = the number of maximal independent sets containing $k$ vertices of $\tilde{Q}(n)$
\vskip 5 pt

\noindent and let
\begin{align}
\overline{Q}(x, y) = \sum_{n \geq 0}\sum_{k \geq 0}\overline{q}(n, k)x^{n}y^{k}\notag\\
\tilde{Q}(x, y) = \sum_{n \geq 0}\sum_{k \geq 0}\tilde{q}(n, k)x^{n}y^{k}\notag
\end{align}
\noindent be the bi-variate generating functions of $\overline{q}(n, k)$ and $\tilde{q}(n, k)$, respectively. Now, we are ready to prove Theorem \ref{bivariate-qh}.
\vskip 5 pt

\noindent \emph{Proof of Theorem \ref{bivariate-qh}} Let $D$ be a maximal independent set of $Q(n)$ containing $k$ vertices. We distinguish $2$ cases. 
\vskip 5 pt

\noindent \textbf{Case 1:} $x^{n}_{4} \in D$\\
\indent Thus, $x^{n}_{3}, x^{n}_{5} \notin D$. Removing the vertices $x^{n}_{3}, x^{n}_{4}, x^{n}_{5}$ from $Q(n)$ result in $\overline{Q}(n - 1)$. There are $\overline{q}(n - 1, k - 1)$ possibilities of $D$.
\vskip 5 pt

\noindent \textbf{Case 2:} $x^{n}_{4} \notin D$\\
\indent We further have the following $3$ subcases.
\vskip 5 pt

\noindent \textbf{Subcase 2.1:} $n^{n}_{3} \notin D$ and $x^{n}_{5} \in D$.\\
\indent Thus, $x^{n}_{6} \notin D$ and $x^{n}_{2} \in D$. This implies that $x^{n}_{1} \notin D$. Removing all the vertices of the $n^{th}$ hexagon results in $\tilde{Q}(n - 2)$. We have that there are $\tilde{q}(n - 2, k - 2)$ possibilities of $D$.
\vskip 5 pt

\noindent \textbf{Subcase 2.2:} $n^{n}_{3} \in D$ and $x^{n}_{5} \notin D$.\\
\indent Thus, $x^{n}_{2} \notin D$ and $x^{n}_{6} \in D$. This implies that $x^{n}_{1} \notin D$. Removing all the vertices of the $n^{th}$ hexagon results in $\tilde{Q}(n - 2)$. We have that there are $\tilde{q}(n - 2, k - 2)$ possibilities of $D$.
\vskip 5 pt

\noindent \textbf{Subcase 2.3:} $n^{n}_{3} \in D$ and $x^{n}_{5} \in D$.\\
\indent Thus, $x^{n}_{2}, x^{n}_{6} \notin D$. Removing the vertices $x^{n}_{2}, ..., x^{n}_{6}$ results in $Q(n - 1)$. We have that there are $q(n - 1, k - 2)$ possibilities of $D$.
\vskip 5 pt

\indent From, the two cases, we have that
\begin{align}\label{qh11}
q(n, k) = \overline{q}(n - 1, k - 1) + 2\Tilde{q}(n - 2, k - 2) + q(n - 1, k - 2)
\end{align}
\noindent For $n \geq 2 $ and $k \geq 2 $, we multiply $x^{n}y^{k}$ throughout (\ref{qh11}) and sum over all $x^{n}y^{k}$. Thus, we have that
\begin{align}\label{qh22}
\sum_{n \geq 2}\sum_{k \geq 2}q(n, k)x^{n}y^{k} &= \sum_{n \geq 2}\sum_{k \geq 2}\overline{q}(n - 1, k - 1)x^{n}y^{k} + 2\sum_{n \geq 2}\sum_{k \geq 2}\Tilde{q}(n - 2, k - 2)x^{n}y^{k}\notag\\
&+ \sum_{n \geq 2}\sum_{k \geq 2}q(n - 1, k - 2)x^{n}y^{k}.
\end{align}

\noindent We first consider the term $\displaystyle{\sum_{n \geq 2}\sum_{k \geq 2}q(n, k)x^{n}y^{k}}$. It can be checked that

\begin{center}
    $q(0, 0) = 1, q(0, k) = 0$ for all $k \geq 1$
    \vskip 5 pt
    $q(1, 0) = 0, q(1, 1) = 0, q(1, 2) = 3 ,q(1, 3) = 2$, and $ q(1, k) = 0$ for all $k \geq 4$
    \vskip 5 pt
\end{center}

\begin{align}\label{qh1}
\sum_{n \geq 2}\sum_{k \geq 2}q(n, k)x^{n}y^{k} &= \sum_{n \geq 2}\sum_{k \geq 2}q(n, k)x^{n}y^{k}  + q(0, 0)  + q(1, 2)xy^{2} + q(1, 3)xy^{3} -q(0, 0)\notag\\
&-q(1, 2)xy^{2} - q(1, 3)xy^{3}\notag\\
                  &= \sum_{n \geq 0}\sum_{k \geq 0}q(n, k)x^{n}y^{k} - 1 - 3xy^{2} - 2xy^{3} \notag\\
                  &= Q(x,y) - 1 - 3xy^2 - 2xy^{3}
\end{align}

\noindent Now, we consider the term $\displaystyle{\sum_{n \geq 2}\sum_{k \geq 2}\overline{q}(n - 1, k - 1)x^{n}y^{k}}$. Clearly,

\begin{center}
    $\overline{q}(0, 0) = 0, \overline{q}(0, 1) = 1 , \overline{q}(0, 2) = 1$ and $\overline{q}(0, k) = 0$ for all $k \geq 3$
    \vskip 5 pt
\end{center}

\begin{align}\label{qhbar1}
\sum_{n \geq 2}\sum_{k \geq 2}\overline{q}(n - 1, k - 1)x^{n}y^{k}
    &= xy(\sum_{n \geq 2}\sum_{k \geq 2}\overline{q}(n - 1, k - 1)x^{n - 1}y^{k - 1})\notag\\
    &= xy(\sum_{n \geq 2}\sum_{k \geq 2}\overline{q}(n - 1, k - 1)x^{n - 1}y^{k - 1} + \overline{q}(0,1)y + \overline{q}(0,2)y^{2}\notag\\
    &- \overline{q}(0,1)y - \overline{q}(0,2)y^{2})\notag\\
    &= xy(\sum_{n \geq 0}\sum_{k \geq 0}\overline{q}(n, k)x^{n}y^{k} - y - y^{2})\notag\\
    &= xy\overline{Q}(x, y) - xy^{2} - xy^{3}.
\end{align}

\noindent Next, we consider the term $\displaystyle{\sum_{n \geq 2}\sum_{k \geq 2}2\Tilde{q}(n - 2, k - 2)x^{n}y^{k}}$. It can be check that,

\begin{align}\label{qhtil1}
\sum_{n \geq 2}\sum_{k \geq 2}2\Tilde{q}(n - 2, k - 2)x^{n}y^{k}
    &= 2x^{2}y^{2}(\sum_{n \geq 2}\sum_{k \geq 2}\Tilde{q}(n - 2, k - 2)x^{n - 2}y^{k - 2})\notag\\
    &= 2x^{2}y^{2}\sum_{n \geq 0}\sum_{k \geq 0}\Tilde{q}(n , k)x^{n}y^{k}\notag\\
    &= 2x^{2}y^{2}\Tilde{Q}(x, y).
\end{align}

\noindent Finally, we consider the term
$\displaystyle{\sum_{n \geq 2}\sum_{k \geq 2}q(n - 1 , k - 2)x^{n}y^{k}}$. It can be checked that, 

\begin{align}\label{qh2}
\sum_{n \geq 2}\sum_{k \geq 2}q(n - 1, k - 2)x^{n}y^{k}
    &= xy^{2}(\sum_{n \geq 2}\sum_{k \geq 2}q(n - 1, k - 2)x^{n - 1}y^{k - 2})\notag\\
    &= xy^{2}(\sum_{n \geq 2}\sum_{k \geq 2}q(n - 1, k - 2)x^{n - 1}y^{k - 2} + q(0,0)-q(0,0))\notag\\
    &= xy^{2}(\sum_{n \geq 0}\sum_{k \geq 0}q(n, k)x^{n}y^{k} - 1)\notag\\
    &= xy^{2}Q(x, y) - xy^{2}
\end{align}

\noindent Plugging (\ref{qh1}), (\ref{qhbar1}) , (\ref{qhtil1}) and (\ref{qh2}) to (\ref{qh22}), we have 
\begin{align}
Q(x,y) - 1- 3xy^{2} - 2xy^{3} &= xy\overline{Q}(x, y) - xy^{2} - xy^{3} + 2x^{2}y^{2}\Tilde{Q}(x, y) + xy^{2}Q(x, y) - xy^{2}\notag\\
\end{align}

\noindent which can be solved that

\begin{align}\label{semifinalqh4}
Q(x,y) &= 1 + xy^2 + xy^{3} + xy\overline{Q}(x, y) + 2x^{2}y^{2}\Tilde{Q}(x, y) + xy^{2}Q(x, y).
\end{align}
\vskip 5 pt

\indent Next, we will find the recurrence relation of $\overline{Q}(n)$. Let $D$ be a maximal independent set of $\overline{Q}(n)$ containing $k$ vertices. We distinguish $2$ cases. 
\vskip 5 pt

\noindent \textbf{Case 1:} $x^{n}_{6} \in D$\\
\indent Thus, $x^{n}_{1}, x^{n}_{5}, x^{n+1}_{2}, x^{n+1}_{6} \notin D$. We further distinguish $2$ cases.
\vskip 5 pt

\noindent \textbf{Case 1.1:} $x^{n}_{2}, x^{n}_{4} \in D$\\
\indent Removing the vertices of the $n^{th}$ hexagon and $x^{n+1}_{2}, x^{n+1}_{6}$ from $\overline{Q}(n)$ results in $\tilde{Q}(n - 2)$. There are $\tilde{q}(n - 2, k - 3)$ possibilities of $D$.
\vskip 5 pt

\noindent \textbf{Case 1.2:} $x^{n}_{3} \in D$\\
\indent Thus, $x^{n}_{2}, x^{n}_{4} \notin D$. Removing the vertices of the $n^{th}$ hexagon and $x^{n+1}_{2}, x^{n+1}_{6}$ from $\overline{Q}(n)$ results in $\tilde{Q}(n - 2)$. There are $\tilde{q}(n - 2, k - 2)$ possibilities of $D$.
\vskip 5 pt

\noindent \textbf{Case 2:} $x^{n}_{6} \notin D$\\
\indent By the maximality of $D$, $x^{n+1}_{2}, x^{n+1}_{6} \in D$. Removing $x^{n}_{6}, x^{n+1}_{2}, x^{n+1}_{6}$ from $\overline{Q}(n)$ results in $\tilde{Q}(n - 1)$. There are $\tilde{q}(n - 1, k - 2)$ possibilities of $D$.
\vskip 5 pt

\indent From all the cases, we have that
\begin{align}\label{qhbar2}
\overline{q}(n, k) = \Tilde{q}(n - 2, k - 3) + \Tilde{q}(n - 1, k - 2) + \Tilde{q}(n - 2, k - 2)
\end{align}
\vskip 5 pt
\noindent For $n \geq 2 $ and $k \geq 3 $, we multiply $x^{n}y^{k}$ throughout (\ref{qh11}) and sum over all $x^{n}y^{k}$. Thus, we have that
\begin{align}\label{qhbar22}
\sum_{n \geq 2}\sum_{k \geq 3}\overline{q}(n, k)x^{n}y^{k} = &\sum_{n \geq 2}\sum_{k \geq 3}\Tilde{q}(n - 2, k - 3)x^{n}y^{k} + \sum_{n \geq 2}\sum_{k \geq 3}\Tilde{q}(n - 1, k - 2)x^{n}y^{k} \notag\\
&+ \sum_{n \geq 2}\sum_{k \geq 3}\Tilde{q}(n - 2, k - 2)x^{n}y^{k}
\end{align}

\noindent We first consider the term $\displaystyle{\sum_{n \geq 2}\sum_{k \geq 3}\overline{q}(n , k )x^{n}y^{k}}$. It can be checked that, 

\begin{center}
    $\overline{q}(0, 0) = 0, \overline{q}(0, 1) = 1, \overline{q}(0, 2) =1$ and $\overline{q}(0, k) = 0$ for all $k \geq 3$
    \vskip 5 pt
    $\overline{q}(1, 0) = 0, \overline{q}(1, 1) = 0, \overline{q}(1, 2) =1 , \overline{q}(1, 3) = 1$ 
    \vskip 5 pt
     $\overline{q}(1, 4) = 3, \overline{q}(1, 5) =1$ and $\overline{q}(1, k) = 0$ for all $k \geq 6$
    \vskip 5 pt
\end{center}

\begin{align}\label{qhbar3}
\sum_{n \geq 2}\sum_{k \geq 3}\overline{q}(n, k)x^{n}y^{k}
    &= \sum_{n \geq 2}\sum_{k \geq 3}\overline{q}(n,k)x^{n}y^{k} + \overline{q}(0, 1)y +\overline{q}(0, 2)y^{2}
    +\overline{q}(1, 2)xy^{2} +\overline{q}(1, 3)xy^{3}\notag\\
    &+\overline{q}(1, 4)xy^{4} +\overline{q}(1, 5)xy^{5} -\overline{q}(0, 1)y - \overline{q}(0, 2)y^{2}\notag\\
    &-\overline{q}(1, 2)xy^{2} -\overline{q}(1, 3)xy^{3} -\overline{q}(1, 4)xy^{4} -\overline{q}(1, 5)xy^{5}\notag\\
    &= \sum_{n \geq 0}\sum_{k \geq 0}\overline{q}(n , k )x^{n}y^{k} - y - y^{2} -xy^{2} - xy^{3} - 3xy^{4} - xy^{5}\notag\\
    &= \overline{Q}(x, y) - y - y^{2} - xy^{2} - xy^{3} - 3xy^{4} - xy^{5}
\end{align}

It can be checked that, 

\begin{center}
    $\Tilde{q}(0, 0) = 0, \Tilde{q}(0, 1) = 0 , \Tilde{q}(0, 2) = 3, \Tilde{q}(0,3) =1$ and $\Tilde{q}(0, k) = 0$ for all $k \geq 4$
    \vskip 5 pt
    $\Tilde{q}(1, 0) = 0, \Tilde{q}(1, 1) = 0 , \Tilde{q}(1, 2) = 0, \Tilde{q}(1,3) =1$
    \vskip 5 pt
    $\Tilde{q}(1, 4) = 10, \Tilde{q}(1, 5) = 4$ and $\Tilde{q}(1, k) = 0$ for all $k \geq 6$
    \vskip 5 pt
\end{center}

\noindent We next consider the term $\displaystyle{\sum_{n \geq 2}\sum_{k \geq 3}\Tilde{q}(n - 2 , k - 3)x^{n}y^{k}}$, clearly
\begin{align}\label{qhtil2}
\sum_{n \geq 2}\sum_{k \geq 3}\Tilde{q}(n - 2 , k - 3)x^{n}y^{k}
    &= x^{2}y^{3}(\sum_{n \geq 2}\sum_{k \geq 3}\Tilde{q}(n - 2, k - 3 )x^{n - 2}y^{k - 3})\notag\\
    &= x^{2}y^{3}\Tilde{Q}(x, y)
\end{align}    
\vskip 5 pt

\noindent We next consider the term $\displaystyle{\sum_{n \geq 2}\sum_{k \geq 3}\Tilde{q}(n - 1 , k - 2)x^{n}y^{k}}$, clearly
\begin{align}\label{qhtil3}
\sum_{n \geq 2}\sum_{k \geq 3}\Tilde{q}(n - 1 , k - 2)x^{n}y^{k}
    &= xy^{2}(\sum_{n \geq 2}\sum_{k \geq 3}\Tilde{q}(n - 1, k - 2)x^{n - 1}y^{k - 2} - \Tilde{q}(0,2)y^{2} - \Tilde{q}(0,3)y^{3} )\notag\\
    &= xy^{2}\Tilde{Q}(x, y) - 3xy^{4} - xy^5
\end{align}    
\vskip 5 pt

\noindent Finally, we consider the term $\displaystyle{\sum_{n \geq 2}\sum_{k \geq 3}\Tilde{q}(n - 2 , k - 2 )x^{n}y^{k}}$.

\begin{align}\label{qhtil4}
\sum_{n \geq 2}\sum_{k \geq 3}\Tilde{q}(n - 2 , k - 2)x^{n}y^{k}
    &= x^{2}y^{2}(\sum_{n \geq 2}\sum_{k \geq 3}\Tilde{q}(n - 2, k - 2 )x^{n - 2}y^{k - 2})\notag\\
    &= x^{2}y^{2}\Tilde{Q}(x, y)
\end{align}    
\vskip 5 pt

\noindent Plugging (\ref{qhbar3}) , (\ref{qhtil2}) , (\ref{qhtil3}) and (\ref{qhtil4}) to (\ref{qhbar22}), we have 

\begin{align}
\overline{Q}(x, y) - y - y^{2} - xy^{2} - xy^{3} - 3xy^{4} - xy^{5} &= x^{2}y^{3}\Tilde{Q}(x, y) + xy^{2}\Tilde{Q}(x, y) - 3xy^{4} - xy^5\notag\\
&+ x^{2}y^{2}\Tilde{Q}(x, y)
\end{align}

\noindent which can be solved that

\begin{align}\label{semifinalqhbar4}
\overline{Q}(x, y) &= y + y^{2} + xy^{2} + xy^{3} + (x^{2}y^{3} + xy^{2} + x^{2}y^{2})\Tilde{Q}(x, y)
\end{align}
\vskip 5 pt

\indent Next, we will find the recurrence relation of $\tilde{Q}(n)$. Let $D$ be a maximal independent set of $\tilde{Q}(n)$ containing $k$ vertices. We distinguish $2$ cases. 
\vskip 5 pt

\noindent \textbf{Case 1:} $x^{n+1}_{5} \in D$\\
\indent Thus, $x^{n+1}_{4} \notin D$. We further distinguish $2$ cases.
\vskip 5 pt

\noindent \textbf{Case 1.1:} $x^{n + 1}_{3} \in D$\\
\indent Thus, $x^{n + 1}_{2} \notin D$. Removing the vertices $x^{n+1}_{2}, ..., x^{n+1}_{5}$ from $\tilde{Q}(n)$ results in $Q(n)$. There are $q(n, k - 2)$ possibilities of $D$.
\vskip 5 pt

\noindent \textbf{Case 1.2:} $x^{n + 1}_{3} \notin D$\\
\indent By maximality of $D$, $x^{n + 1}_{2} \in D$. Thus, $x^{n}_{6} \notin D$. Removing the vertices $x^{n}_{6}, x^{n+1}_{2}, ..., x^{n+1}_{5}$ from $\tilde{Q}(n)$ results in $\tilde{Q}(n - 1)$. There are $\tilde{q}(n - 1, k - 2)$ possibilities of $D$.
\vskip 5 pt

\noindent \textbf{Case 2:} $x^{n+1}_{5} \notin D$\\
\indent By the maximality of $D$, $x^{n+1}_{4} \in D$. Thus, $x^{n+1}_{3} \notin D$. We further distinguish $2$ cases.
\vskip 5 pt

\noindent \textbf{Case 2.1:} $x^{n + 1}_{2} \in D$\\
Thus, $x^{n}_{6} \notin D$. Removing $x^{n}_{6}, x^{n+1}_{2}, x^{n+1}_{6}$ from $\tilde{Q}(n)$ results in $\tilde{Q}(n - 1)$. There are $\tilde{q}(n - 1, k - 2)$ possibilities of $D$.
\vskip 5 pt

\noindent \textbf{Case 2.2:} $x^{n + 1}_{2} \notin D$\\
By the maximality of $D$, $x^{n - 1}_{6} \in D$. Thus, $x^{n}_{5} \notin D$. 

\noindent \textbf{Case 2.2.1:} $x^{n}_{2}, x^{n}_{4} \in D$\\
\indent Removing the vertices of the $n^{th}$ hexagon and $x^{n+1}_{2}, ..., x^{n+1}_{6}$ from $\tilde{Q}(n)$ results in $\tilde{Q}(n - 2)$. There are $\tilde{q}(n - 2, k - 4)$ possibilities of $D$.
\vskip 5 pt

\noindent \textbf{Case 2.2.2:} $x^{n}_{3} \in D$\\
\indent Thus, $x^{n}_{2}, x^{n}_{4} \notin D$. Removing the vertices of the $n^{th}$ hexagon and $x^{n+1}_{2}, ..., x^{n+1}_{6}$ from $\tilde{Q}(n)$ results in $\tilde{Q}(n - 2)$. There are $\tilde{q}(n - 2, k - 3)$ possibilities of $D$.
\vskip 5 pt

\indent From all the cases, we have that
\begin{align}\label{qhtil44}
    \Tilde{q}(n, k) = q(n ,k - 2) + 2\Tilde{q}(n - 1, k - 2) + \Tilde{q}(n - 2,k - 3) + \Tilde{q}(n - 2,k - 4)
\end{align}
\vskip 5 pt

\noindent For $n \geq 2$ and $k \geq 4$, we multiply $x^{n}y^{k}$ throughout (\ref{qhtil44}) and sum over all $x^{n}y^{k}$. Thus, we have that

\begin{align}\label{qhtil44444}
\sum_{n \geq 2}\sum_{k \geq 4}\Tilde{q}(n, k)x^{n}y^{k} =& \sum_{n \geq 2}\sum_{k \geq 4}q(n ,k - 2)x^{n}y^{k} + \sum_{n \geq 2}\sum_{k \geq 4}2\Tilde{q}(n - 1, k - 2)x^{n}y^{k}  \notag\\ &+ \sum_{n \geq 2}\sum_{k \geq 4}\Tilde{q}(n - 2,k - 3)x^{n}y^{k} + \sum_{n \geq 2}\sum_{k \geq 4}\Tilde{q}(n - 2,k - 4)x^{n}y^{k}
\end{align}

\noindent We first consider the term $\displaystyle{\sum_{n \geq 2}\sum_{k \geq 4}\Tilde{q}(n , k )x^{n}y^{k}}$. 

\begin{align}\label{qhtil5}
\sum_{n \geq 2}\sum_{k \geq 4}\Tilde{q}(n , k)x^{n}y^{k}
    &= \sum_{n \geq 2}\sum_{k \geq 4}\Tilde{q}(n , k)x^{n}y^{k}
    + \Tilde{q}(0,2)y^{2} + \Tilde{q}(0,3)y^{3}\notag\\
    &+\Tilde{q}(1,3)xy^{3} + \Tilde{q}(1,4)xy^{4} + \Tilde{q}(1,5)xy^{5}\notag\\
    &-\Tilde{q}(0,2)y^{2} - \Tilde{q}(0,3)y^{3} - \Tilde{q}(1,3)xy^{3} - \Tilde{q}(1,4)xy^{4} - \Tilde{q}(1,5)xy^{5}\notag\\
    &= \sum_{n \geq 0}\sum_{k \geq 0}\Tilde{q}(n , k )x^{n}y^{k} - 3y^{2} - y^{3} - xy^{3} -10xy^{4} -4xy^{5}\notag\\
    &= \Tilde{Q}(x, y) - 3y^{2} - y^{3} - xy^{3} -10xy^{4} -4xy^{5}
\end{align}    
\vskip 5 pt

\noindent We next consider the term $\displaystyle{\sum_{n \geq 2}\sum_{k \geq 4}q(n, k - 2)x^{n}y^{k}}$, clearly

\begin{align}\label{qh3}
\sum_{n \geq 2}\sum_{k \geq 4}q(n, k-2)x^{n}y^{k} &= y^{2}(\sum_{n \geq 2}\sum_{k \geq 4}q(n, k-2)x^{n}y^{k-2} + q(0, 0) +q(1,2)xy^{2} +q(1,3)xy^{3}\notag\\
                  &-q(0,0)- (1,2)xy^{2} - q(1,3)xy^{3})\notag\\
                  &= y^{2}(\sum_{n \geq 0}\sum_{k \geq 0}q(n, k)x^{n}y^{k} - 1 - 3xy^{2} - 2xy^{3})\notag\\
                  &= y^{2}Q(x,y) - y^{2} - 3xy^{4} - 2xy^{5}
\end{align}
\vskip 5 pt

\noindent We next consider the term $\displaystyle{\sum_{n \geq 2}\sum_{k \geq 4}2\Tilde{q}(n -1, k - 2)x^{n}y^{k}}$, clearly

\begin{align}\label{qhtil6}
\sum_{n \geq 2}\sum_{k \geq 4}2\Tilde{q}(n - 1 , k - 2)x^{n}y^{k}
    &= 2xy^{2}(\sum_{n \geq 2}\sum_{k \geq 4}\Tilde{q}(n - 1 , k - 2 )x^{n - 1}y^{k - 2}\notag\\
    &+ \Tilde{q}(0,2)y^{2} + \Tilde{q}(0,3)y^{3} -\Tilde{q}(0,2)y^{2} - \Tilde{q}(0,3)y^{3})\notag\\
    &= 2xy^{2}(\sum_{n \geq 0}\sum_{k \geq 0}\Tilde{q}(n , k )x^{n}y^{k} - 3y^{2} - y^{3})\notag\\
    &= 2xy^{2}\Tilde{Q}(x, y) - 6xy^{4} - 2xy^{5}
\end{align}
\vskip 5 pt

\noindent We next consider the term $\displaystyle{\sum_{n \geq 2}\sum_{k \geq 4}\Tilde{q}(n - 2, k - 3)x^{n}y^{k}}$.

\begin{align}\label{qhtil7}
\sum_{n \geq 2}\sum_{k \geq 4}\Tilde{q}(n - 2 , k - 3)x^{n}y^{k}
    &= x^{2}y^{3}(\sum_{n \geq 2}\sum_{k \geq 4}\Tilde{q}(n - 2, k - 3 )x^{n - 2}y^{k - 3})\notag\\
    &= x^{2}y^{3}\Tilde{Q}(x, y)
\end{align}    
\vskip 5 pt

\noindent Finally, we consider the term
$\displaystyle{\sum_{n \geq 2}\sum_{k \geq 4}\Tilde{q}(n - 2 , k - 4)x^{n}y^{k}}$. 

\begin{align}\label{qhtil8}
\sum_{n \geq 2}\sum_{k \geq 4}\Tilde{q}(n - 2 , k - 4)x^{n}y^{k}
    &= x^{2}y^{4}(\sum_{n \geq 2}\sum_{k \geq 4}\Tilde{q}(n - 2, k - 4 )x^{n - 2}y^{k - 4})\notag\\
    &= x^{2}y^{4}\Tilde{Q}(x, y)
\end{align}    
\vskip 5 pt

\noindent Plugging (\ref{qhtil5}) , (\ref{qh3}) , (\ref{qhtil6}) , (\ref{qhtil7}) and (\ref{qhtil8}) to (\ref{qhtil44444}), we have

\begin{align}
\Tilde{Q}(x, y) - 3y^{2} - y^{3} - xy^{3} -10xy^{4} -4xy^{5} &= y^{2}Q(x,y) - y^{2} - 3xy^{4} - 2xy^{5}\notag\\
&+ 2xy^{2}\Tilde{Q}(x, y) - 6xy^{4} - 2xy^{5}\notag\\
&+ x^{2}y^{3}\Tilde{Q}(x, y) + x^{2}y^{4}\Tilde{Q}(x, y)
\end{align}

\noindent which can be solved that

\begin{align}\label{semifinalqhtil7}
\Tilde{Q}(x, y) = \frac{1}{1 - 2xy^{2} - x^{2}y^{3} - x^{2}y^{4}}(2y^{2} + y^{3} + y^{2}Q(x,y) + xy^{3} + xy^{4})
\end{align}
\vskip 5 pt

\noindent By plugging (\ref{semifinalqh4}) to (\ref{semifinalqhbar4}) to (\ref{semifinalqhtil7}) , we have
\begin{align}
Q(x, y) = \frac{x^{2}y^{6}+2xy^{3}+1}{1 - x^{2}y^{5} - x^{2}y^{4} - x^{2}y^{3} - 3xy^2}
\end{align}
\noindent as required.
\qed

\noindent We will prove Theorem \ref{recurrence-qh}
\vskip 5 pt

\noindent \emph{Proof of Theorem \ref{recurrence-qh}} By Theorem \ref{bivariate-qh} with $y = 1$, we have that
\vskip 5 pt

\begin{align}
\sum_{n \geq 0}q(n)x^{n} &= Q(x, 1)\notag\\
                         &= \frac{x^{2} + 2x + 1}{1 - 3x - 3x^{2}}
\end{align}
 
 \noindent which can  be solved that

\begin{align}\label{recursived-qh}
x^{2} + 2x + 1 &= (1 - 3x - 3x^{2})\sum_{n \geq 0}q(n)x^{n}\notag\\
                         &= \sum_{n \geq 0}q(n)x^{n} - 3\sum_{n \geq 0}q(n)x^{n + 1} - 3\sum_{n \geq 0}q(n)x^{n + 2}\notag\\
                         &= \sum_{n \geq 0}q(n)x^{n} - 3\sum_{n \geq 1}q(n - 1)x^{n} - 3\sum_{n \geq 2}q(n - 2)x^{n}\notag\\
                         &= q(0) + q(1)x - 3q(0)x + \sum_{n \geq 2}(q(n) - 3q(n - 1) - 3q(n - 2))x^{n}.
\end{align}

\noindent Because the order of the polynomial on the left hand side of (\ref{recursived-qh}) is two, the coefficients of $x^{n}$ for all $n \geq 3$ must be $0$. Thus, $q(n) - 3q(n - 1) - 3q(n - 2) = 0$ implying that

\begin{align*}
    q(n) = 3q(n - 1) + 3q(n - 2).
\end{align*}


\begin{thebibliography}{99}
\bibitem{AH1} {BAbrishami, A., Henning, M. A., 2018, \ ``Independent domination in subcubic graphs of girth at least six'', \emph{Discrete Mathematics,} Vol. 341(1), 155-164.}

\bibitem{AH} {Abrishami, G., Henning, M. A., 2022, \ ``An improved upper bound on the independent domination number in cubic graphs of girth at least six'', \emph{Graphs and Combinatorics,} Vol. 38(2):50.}


\bibitem{BH} {Babikir, A., Henning, M. A., 2020, \ ``Domination versus independent domination in graphs of small regularity'', \emph{Discrete Mathematics,} Vol. 343(7), 111727.}



\bibitem{BN} {Bak, J., Newman, D. J., 2010,\ ``Complex analysis'', Springer, New York.}


\bibitem{BM} {Bisdorff, R., Marichal, J. L., 2008, \ ``Counting non-isomorphic maximal independent sets of the $n$-cycle graph'', \emph{Journal of Integer Sequences,} Vol. 11, pp. 1-16.}





\bibitem{CXL} {Chen, A., Xiong, X., Lin, F., 2016, \ ``Distance-based topological indices of the tree-like polyphenyl systems'', \emph{Applied Mathematics and Computation,} Vol. 281, pp. 233-242.}


\bibitem{CGR} {Cruz, R., Gutman, I., Rada, J., 2021, \ ``Sombor index of chemical graphs'', \emph{Applied Mathematics and Computation,} Vol. 399, 126018.}


\bibitem{DG} {Das, K. C., Gutman, I., 2022, \ ``On Sombor index of trees'', \emph{Applied Mathematics and Computation,} Vol. 412, 126575.}


\bibitem{D} {Dobrynin, 2022, \ ``Wiener index of families of unicyclic graphs obtained from a tree'', \emph{MATCH, Communication in Mathematical and Computer Chemistry,} Vol. 88, pp. 461-470.}






 \bibitem{DL} {Doslic, T. and Lits, M. S., 2011, \ ``Matchings and independent  sets in polyphenylene chains'', \emph{Communications in Mathematical and in Computer Chemistry ,} Vol. 67, pp. 314-329.}




\bibitem{DM} {Doslic, T. and Maloy, F., 2010, \ ``Chain hexagonal cacti: Matching and independent sets'', \emph{Discrete Mathematics,} Vol. 310, pp. 1676-1690.}



\bibitem{DZ} {Doslic, T. and Zubac, I., 2015, \ ``Counting maximal matchings in linear polymers'', \emph{Ars Mathematica Contemporanea,} Vol. 11(2), pp. 256-265.}


\bibitem{DR} {Drmota, M., Ramos, L., Requile, C. and Rue, J., 2020, \ ``Maximal independent sets and maximal matchings in series-parallel and related graph classes'', \emph{Electronic Journal of Combinatorics,} Vol. 27(1), pp. 1-36.}


\bibitem{DM1} {Dyer, M. and Muller, H., 2019, \ ``Counting independent sets in cocomparability graphs'', \emph{Information Processing Letters,} Vol. 144, pp. 31-36.}



\bibitem{DFJ} {Dyer, M., Frieze, A. and Jerrum, M., 2002, \ ``On counting independent sets in sparse graphs'', \emph{SIAM Journal on Computing,} Vol. 31(5), pp. 1527-1541.}



\bibitem{F1} {Farrell, E.J., 1987, \ ``Matchings in hexagonal cacti'', \emph{International Journal of Mathematics and Mathematical Sciences}, Vol. 10, pp. 321-338.}


\bibitem{F2} {Farrell, E.J., 1998, \ ``Matchings in rectangular cacti'', \emph{International Journal of Mathematics and Mathematical Sciences}, Vol. 9, pp. 163-183.}


\bibitem{F3} {Farrell, E.J., 2000, \ ``Matchings in pentagonal cacti'', \emph{International Journal of Mathematics and Mathematical Sciences}, Vol. 11, pp. 109-126.}



\bibitem{F4} {Farrell, E.J., 2000, \ ``Matchings in triangular cacti'', \emph{International Journal of Mathematics and Mathematical Sciences}, Vol. 11, pp. 85-98.}





\bibitem{G} {Gao, W., 2022, \ ``Trees with maximum vertex-degree-based topological indices'', \emph{MATCH, Communication in Mathematical and Computer Chemistry,} Vol. 88, pp. 535-552.}



\bibitem{GGG} {Griggs, J. R., Grinstead C. M.  and Guichard, D. R., 1988, \ ``The number of maximal independent sets in a connected graph'', \emph{Discrete Mathematics,} Vol. 68, pp. 211-220.}


\bibitem{Gu1} {Gutman, I., 2017, \ ``On coindices of graphs and their complements'', \emph{Applied Mathematics and Computation,} Vol. 305, pp. 161-165.}


\bibitem{Gu2} {Gutman, I., 2020, \ ``Cycle energy and its size independence'', \emph{Discrete Applied Mathematics,} Vol. 284, pp. 534-537.}


\bibitem{HN} {Harary, F. and Norman, R., 1953, \ ``The dissimilarity characteristic of Husimi trees'', \emph{Annals of Mathematics}, Vol. 58, pp. 134-141.}


\bibitem{HP} {Harary, F. and Palmer, E. M., 1973, \ ``Graphical enumeration'', Academic Press, New York.}


\bibitem{HU} {Harary, F. and Uhlenbeck, G. E., 1953, \ ``On the number of Husimi trees'', \emph{Proceedings of the National Academy of Sciences}, Vol. 39, pp. 315-322.}



\bibitem{H} {Husimi, K., 1950, \ ``Note on Mayer’s theory of cluster integrals'', \emph{The Journal of Chemical Physics},  Vol. 18, pp. 682–684.}











\bibitem{JC} {Jou, M. J. and Chang, G. J., 2000, \ ``The number of maximum independent sets in graphs'', \emph{Taiwanese Journal of Mathematics,} Vol. 4(4), pp. 685-695.}


\bibitem{JL} {Jou, M. J. and Lin, J. J., 2016, \ ``An alternative proof of the largest number of maximal independent sets in connected graphs having at most two cycles'', \emph{Open Journal of Discrete Mathematics,} Vol. 6, pp. 227-237.}



\bibitem{LZZ} {Li, S., Zhang, H. and Zhang, X., 2013, \ ``Maximal independent sets in bipartite graphs with at least one cycle'', \emph{Discrete Mathematics and Theoretical Computer Science,} Vol. 15, pp. 243-258.}




\bibitem{LYH} {Liu, H., You, L., Huang, Y., 2022, \ ``Odering chemical graphs by Sombor indices and its applications'', \emph{MATCH, Communication in Mathematical and Computer Chemistry,} Vol. 88, pp. 5-22.}



\bibitem{Mayer} {Mayer, J.E. and Mayer, M.G., 1940, \ ``Statistical Mechanics'', John Wiley and Sons, New York.}


\bibitem{MS1} {Merrifield, R. M. and Simmons, H. E., 1980, \ ``The structure of molecular topological spaces'', \emph{Theoretica Chimica Acta}, Vol. 55, pp. 55-75.}


\bibitem{MS2} {Merrifield, R. M. and Simmons, H. E., 1981, \ ``Enumeration of structure-sensitive graphical subsets: Theory'', \emph{Proceedings of the National Academy of Sciences}, Vol. 78, pp. 692-695.}



\bibitem{MS3} {Merrifield, R. M. and Simmons, H. E., 1981, \ `` Enumeration of structure-sensitive graphical  subsets: Calculations'', \emph{Proceedings of the National Academy of Sciences}, Vol. 78, pp. 1329-1332.}



\bibitem{MS4} {Merrifield, R. M. and Simmons, H. E., 1981, \ ``Mathematical description of molecular  Structure: Molecular topology'', \emph{Proceedings of the National Academy of Sciences}, Vol. 78, pp. 2616-2619.}






\bibitem{MM} {Moon, J. W. and Mooser, L., 1965, \ ``On cliques in graphs'', \emph{Israel Journal of Mathematics,} Vol. 3, pp. 23-28.}





\bibitem{OC} {Oz, M. S., Cangul, I. N., 2022, \ ``Enumeration of independent sets in Benzenoid chains'', \emph{MATCH, Communication in Mathematical and Computer Chemistry,} Vol. 88, pp. 93-107.}





\bibitem{Ra} {Rada, J., 2017, \ ``Vertex-degree-based topological indices of hexagonal systems with equal number of edges'', \emph{Applied Mathematics and Computation,} Vol. 296, pp. 270-276.}






\bibitem{R} {Riddell, R.J., 1951, \ ``Contribution to the theory of condensation'', \emph{Doctoral Dissertation in Physics,} University of Michigan.}


\bibitem{SV} {Sagan, B. E. and Vatter, V. R., 2006, \ ``Maximal and maximum independent sets in
graphs with at most $r$ cycles'', \emph{Journal of Graph Theory,} Vol. 53, pp. 283-314.}


\bibitem{SY} {Shalu, M.A. and Yamini, S. D., 2016, \ ``Counting maximal independent sets in power set graphs'', \emph{Journal of Combinatorial Mathematics and Combinatorial Computing,} Vol. 96, pp. 283-291.}



\bibitem{SIA} {Siddiqui, M. K., Imran, M., Ahmad, A., 2016, \ ``On Zagreb indices, Zagreb polynomials of some nanostar dendrimers'', \emph{Applied Mathematics and Computation,} Vol. 280, pp. 132-139.}








\bibitem{SLH} {Song, X., Li, J., He, W., 2020, \ ``On Zagreb eccentricity indices of cacti'', \emph{Applied Mathematics and Computation,} Vol. 383, 125361.}



\bibitem{TBJK} {Tabassum, H., Bokhary, S. A. U. H., Jiarasuksakun, T., Kaemawichanurat, P., 2022, \ ``Counting the numbers of paths of all lengths in symmetric dendrimers and its applications'', \emph{MATCH, Communication in Mathematical and Computer Chemistry,} Vol. 88, pp. 659-681.}









\bibitem{U} {Uhlenbeck, G.E., 1950, \ ``Some Basic Problems of Statistical Mechanics'', \emph{American Mathematical Society}, Gibbs Lecture.}



\bibitem{W} {Wilf, H. S., 1986, \ ``The Number of maximal independent sets in a tree'', \emph{
SIAM Journal on Algebraic Discrete Methods,} Vol. 7, pp. 125-130.}

\bibitem{W2} {Wilf, H. S., 1990, \ ``Generatingfunctionology'', 
Academic Press, San Diego.}
\end{thebibliography}
\end{document}